# The Floyd-Warshall Algorithm and the Asymmetric TSP

## Howard Kleiman

**Section 3.1 Introduction**

Let $s$ and $t$ be two vertices of a connected weighted graph $G$ represented by the matrix $M_G$. The shortest path problem finds a path between $s$ and $t$ whose total edge weight is minimum. Generally, edge-weight is taken to mean distance but the word is used loosely and may represent some other measurable quantity. Dijkstra's algorithm [1] finds the distance between $s$ and all of the other vertices of $G$. However, it assumes that all edge weights are non-negative. The Floyd-Warshall algorithm [3], [6], finds the shortest paths between all pairs of nodes. Furthermore, unlike the Dijkstra algorithm, it allows arc weights to be negative. In this paper, we use the variant of the F-W algorithm given in [5]. Of particular interest to us, this algorithm allows us to find all cycles of smallest total arc-weight. The running time of the F-W algorithm is $O(n^3)$. In this paper, given a random $n$-cycle, $D$, obtained from the symmetric group, $S_n$, we apply $D^{-1}$ to permute the columns of $M_G$. Thus, if the arc $(a, D_0(b))$ has the weight $w(a, D_0(b))$ in $M_G$, its weight is transformed into $w(a,b)$ in $D^{-1}M_G$. It follows that a negative- weighted cycle, $C$, in $D^{-1}M_G$, say $(a_1, a_2, \dots, a_i)$, has the property that the total edge weight of $DC$ is less than that of $D$. Furthermore, as long as no pair of consecutive vertices of $C$ - $a_j a_{j+1}$ - has the property that $a_{j-1} a_j$ is an arc of $D$, $DC$ will be a derangement, i.e., a permutation that moves all points in $V = \{1, 2, \dots, n\}$. In chapter 1, we used $H$-admissible permutations to obtain hamilton circuits in graphs. In the algorithm given here, we use admissible permutations as defined in chapter 2. Starting with a random $n$-cycle, they were used there to obtain a spanning set of $[\log n]$ disjoint cycles. Here they are used to obtain a rough approximation to an optimal solution of the edge-distance assignment problem defined by the entries in $M_G$. This approximation yields a weight matrix in which the negative entries are generally fewer in number and smaller in absolute value than those in $M_G$. Theorem 3.1 and its corollary allow us to use considerably fewer trials than generally used in Floyd-Warshall. Once we have obtained an optimal solution, $\sigma_{APOPT}$, to the assignment problem, the only cycles we can obtain in $\sigma_{APOPT}^{-1} M_G$ are all positive. We then use the above theorems and the F-W algorithm to obtain the smallest positive cycles of total edge weight less than a fixed positive number, say $N$. As we proceed, we test cycles and sets of disjoint cycles to see if any product is an *n*-cycle. The *n*-cycle of smallest weight that we obtain is an approximate optimal solution to the TSP. The reason that it is not exact is that the F-W algorithm allows us to obtain cycles of smallest weight constructed by always using the shortest distance between any pair of vertices. It is possible that a cycle, $C$, may be a disjoint cycle of a permutation, $s$, such that $\sigma_{APOPT} s = \sigma_{TSPOPT}$ where $C$ contains a subpath *[a, ... ,b]* that is not the shortest path between $a$ and $b$. However, we do give conditions where our version of the F-W algorithm yields $\sigma_{TSPOPT}$. Furthermore, we give necessary conditions in theorem 3.6 that cycles of $s$ not obtainable by F-W must satisfy.

**Section 3.2 Theorems**

**Theorem 3.1** *Let $C = (a_1 a_2 ... a_n)$ be a cycle of length n. Assume that the weight $w(a_i, C(a_i))$ (i=1,2,...,n) corresponds to the arc $(a_i, C(a_i))$ of C. Then if*

$$W = \sum_{i=1}^{i=n} w(a_i, C(a_i)) \leq 0$$

*there exists at least one vertex $a_{i*}$ with $1 \leq a_{i*} \leq n$ such that*

$$S_m = \sum_{j=0}^{j=m} w(a_{i*+j}, C(a_{i*+j})) \leq 0 \qquad (A)$$

*where $m = 0,1,2,...,n-1$ and $i*+j$ is modulo $n$.*

**Proof.** We prove the theorem by induction. Let k = 2. We thus have a 2-cycle. If both arcs have non-positive value, then the theorem is proved. If the non-positive arc, $(a_1 a_2)$, has a smaller weight than a positive one, then the sum of the weights of the two arcs is positive. This can't be the case. Thus, the sum of the two weights is non-positive. Now let the theorem always be true when our cycle has k arcs. Suppose the cycle has k+1 arcs. In what follows, assume that if a value is 0, its sign is negative. Then one of the following is true: (a) there exists a pair of consecutive arcs both of whose values have the same sign or (b) the signs of values of the arcs consecutively alternate in sign. First, we consider (a). Without loss of generality, let the two arcs be $(a_1\ a_2)$ and $(a_2\ a_3)$. Assume that each one has a non-positive value. We now define the arc $(a_1\ a_3)$ where $w(a_1\ a_3) = w(a_1\ a_2) + w(a_2\ a_3)$. Now replace arcs $(a_1\ a_2)$ and $(a_2\ a_3)$ by $(a_1\ a_3)$. The result is a cycle $C'$ containing $k$ arcs. By induction, the theorem holds for $C'$. Now replace $(a_1\ a_3)$ by $(a_1\ a_2)$ and $(a_2\ a_3)$. Let $a_{i*}$ be a determining vertex of $C'$. Then the path, $P_{i*}$, from $a_{i*}$ to $a_1$ is non-positive. If both $w(a_1\ a_2)$ and $w(a_2\ a_3)$ are non-negative, then both $P_{i*} \cup w(a_1\ a_2)$ and $P_{i*} \cup w(a_1\ a_2) \cup w(a_2\ a_3)$ are non-positive. Thus, the theorem is valid in this case. Now assume that both $w(a_1\ a_2)$ and $w(a_2\ a_3)$ are positive. But, since the theorem is valid for $C'$, $P_{i*} \cup w(a_1\ a_3)$ is non-positive. This assures us that each of $P_{i*} \cup w(a_1\ a_2)$ and $P_{i*} \cup w(a_1\ a_2) \cup w(a_2\ a_3)$ is non-positive. Therefore, the theorem holds when the values of $(a_1\ a_2)$ and $(a_2\ a_3)$ are both the same. We now consider case (b). Here, the signs alternate between positive and negative. Since the arcs lie on a cycle, we may assume that the first sign is negative. If k+1 is odd, it is always true that at least one pair of consecutive arcs has the same sign. Thus, assume that k+1 is even. Then, starting with a non-positively-valued arc, we can arrange the arcs of the cycle in pairs where the first arc has a non-positive value, while the second one has a positive value. Suppose the sum of each pair of arcs is positive. Then the sum of the values of the cycle is positive. Therefore, there exists at least one pair of arcs the sum of whose values is non-positive. Remembering that the first arc of each pair is non-positive, we follow the same procedure as in (a): $C'$ is a non-positive cycle containing k arcs. Therefore, $C$ is also positive.

**Corollary 1**. *Suppose that C is a cycle such that*

$$W = \sum_{i=1}^{i=n} w(a_i, C(a_i)) \leq N \qquad (B)$$

*Then there exists a determining vertex $a_{i*}$ such that each partial sum, $S_m$, has the property that*

$$S_m = \sum_{j=0}^{j=m} w(a_{i*+j}, C(a_{i*+j})) \leq N \qquad (C)$$

*always holds. Here m =0,1,2, ... , n-1 while i\* + j is modulo n.*

**Proof**. Subtract $N$ from both sides of (B). Now let the weight of arc $(a_n\ C(a_n))$ become $w(a_n, C(a_n)) - N$. From theorem 3.1,

$$W^* = W - N = \sum_{i=1}^{i=n} w(a_i, C(a_i)) \leq 0 \quad (D)$$

Therefore, we can obtain a determining vertex $a_{i*}$ having the property that every partial sum having $a_{i*}$ as it initial vertex is non-positive. It follows that if we restore $w(a_n, C(a_n))$ to its original value, every partial sum with initial vertex $a_{i*}$ is less than or equal to $N$.

**Corollary 2**. *Let $C$ be a positively-valued cycle of length $n$ obtained from a cost matrix $W$ whose entries may be positive, negative or $0$. Then there exists at least one determining vertex, say $a_{i*}$, of $C$ such that each subpath $S_m$ having initial vertex $a_{i*}$ has the property that*

$$S_m = \sum_{j=0}^{j=m} w(a_{i*+j}, C(a_{i*+j})) \geq 0 \quad (E)$$

*where $m = 0, 1, 2, \ldots, n-1$ and $i*+j$ is modulo $n$.*

**Proof.** Let $N = 0$ in $(B)$. Then multiply each term of $(B)$ by $-1$. It follows that we obtain

$-W = \sum_{i=1}^{i=n} -w(a_i, C(a_i)) \leq 0$. The rest of the proof is similar to that of Corollary 1.

**Example 3.1** We now give an example of how to obtain $i = i'$.

Let $C = (a_1\ a_2\ \ldots\ a_n)$

$a_1 = -7$, $a_2 = -10$, $a_3 = +1$, $a_4 = +2$, $a_5 = -7$, $a_6 = +4$, $a_7 = -9$,

$a_8 = +11$, $a_9 = -2$, $a_{10} = -1$, $a_{11} = -4$, $a_{12} = -4$, $a_{13} = -8$, $a_{14} = +9$,

$a_{15} = +9$, $a_{16} = +21$, $a_{17} = +1$, $a_{18} = -2$, $a_{19} = -1$, $a_{20} = -3$,

$a_{21} = -3$, $a_{22} = -12$, $a_{23} = +6$, $a_{24} = +2$, $a_{25} = +3$

-7 −10 +1 +2 −7 +4 −9 +11 −2 −1 −4 −4 −8 +9 +9 + 21 +1

-2 -1 -3 -3 -12 +6 +2 +3

We now add terms with like signs going from left to right. We place the ordinal number of the first number in each sum above it.

   1    3    5    6   7     8    9    14    18    23

-17  +3  -7  +4  -9  +11  -19  +40   -21   +11

We next add the positive number to the right of each negative number to the negative number.

   1    5             18

-14  -3  +2  +21  -10

We now add terms with like signs going from left to right.

   1         18

-17  +23  -10

Finally, assuming that all points lie on a circle, we add like terms going from left to right. We thus obtain

  18         18

-27  +23  =  -4

This tells us that $i'$ is the eighteenth ordinal number, namely, $-2$.

Thus, the partial sums are: -2, -3, -6, -9, -21, -15, -13, -10, -17, -27,

-26, -24, -31, -27, -36, -25, -27, -28, -32 –36, -44, -35, -26, -5, -4

The following theorem is essential for the use of MIN(M) throughout the algorithm.

**Theorem 3.2 The Floyd-Warshall Algorithm** [3], [6]

*If we perform a triangle operation for successive values $j = 1,2,...,n$, each entry $d_{ik}$ of an n X n cost matrix M becomes equal to the value of the shortest path from i to k provided that M contains no negative cycles.*

The version given here is modeled on theorem 6.4 in [5].

**Proof.** We shall show by induction that that after the triangle operation for $j = j_0$ is executed, $d_{ik}$ is the value of the shortest path with intermediate vertices $v \leq j_0$, for all $i$ and $k$. The theorem holds for $j_0 = 1$ since $v = 0$. Assume that the inductive hypothesis is true for $j = j_0 - 1$ and consider the triangle operation for $j = j_0$:

$$d_{ik} = min\{d_{ik}, d_{ij_0} + d_{j_0k}\}.$$

If the shortest path from $i$ through $k$ with $v \leq j_0$ doesn't pass through $j_0$, $d_{ik}$ will be unchanged by this operation, the first argument in the min-operation will be selected, and $d_{ik}$ will still satisfy the inductive hypothesis. On the other hand if the shortest path from $i$ to $k$ with intermediate vertices $v \leq j_0$ does pass through $j_0$, $d_{ik}$ will be replaced by $d_{ij_0} + d_{j_0k}$. By the inductive hypothesis, $d_{ij_0}$ and $d_{j_0k}$ are both optimal values with intermediate vertices $v \leq j_0 - 1$. Therefore, $d_{ij_0} + d_{j_0k}$ is optimal with intermediate vertices v $v \leq j_0$.

We now give an example of how F-W works.

**Example 3.2.** Let d(1, 3) = 5, d(3, 7) = -2, d(1, 7) = 25. Then

d(1, 3) + d(3, 7) < d(1, 7). Note however, that the intermediate vertex 3 comes from the fact that we have reached column j = 3 in the algorithm. We now substitute d(1, 3) + d(3, 7) = 3 for the entry in (1, 7). Suppose now that d(1, 10) = 7 while d(7, 10) = -5.

$$d(1, 3) + d(3, 7) + d(7,10) = -2 < d(1, 10) = 7.$$

The intermediate vertices are now 3 and 7. Thus, -2 can be substituted for 7 in entry (1, 10). Again, we must have reached column 7 before this substitution could be made. Now suppose that d(10, 1) = 1. Then our negative path implies that the negative cycle (1 3 7 10) exists.

**Section 3.3 The Algorithm**

### PHASE 1

Let M be an $n \times n$ distance matrix. An arbitrary entry of M is d($i$, $j$j).

$V = \{1,2,\ldots,n\}$

**Step 1**. Sort each row of M in ascending order of numerical value. Call the matrix obtained MIN(M). An arbitrary entry ($i, j$) of MIN(M) is *ORDINAL((i, j))* which stands for the ordinal value of entry ($i, j$) in row $i$.

**Step 2.** Construct a random $n$-cycle, $D = (a_1\ a_2\ \ldots\ a_n)$ whose arcs are assigned the respective values of its arcs in M. Let $D^{-1}$ be the inverse of $D$.

**Step 3.** Let the function *ORD* with respect to $D$ be $ORD(i) = a_i$. Construct $ORD^{-1}(a_i) = i$, $i = 1, 2, \ldots, n$.

**Step 4.** Given each arc $(a_i, D(a_i))$ of $D$, obtain

$$DIFF(a_i) = d(a_i, D(a_i)) - d(a_i, \text{MIN(M)}(a_i, 1)).$$

(If MIN(M)($a_i$, 1) exists in $D$, we choose the second smallest value of row $a_i$, MIN(M)($a_i$, 2), etc. .

$MIN_1 = \min\{ DIFF(a_i) \in V \mid d(a_i, D(a_i)) - d(a_i, MIN(M)(a_i, j)), j \in V \}$. We then choose $MIN_2, \ldots, MIN_{[\log n]}$, corresponding to the second, third,…,([log n] + 1)-th smallest values of $DIFF(a_i)$ $(i=1,2,\ldots n)$. Generally, $MIN_1$ is negative.

**Step 5.** As we pointed out in the introduction, suppose that

$D = (1\ 7\ 12\ 17\ 20\ 15)$ is considered to be a permutation on the points in $\{1, 2, \ldots, 20\}$. Then, given arcs $(17, D(15))$, $(15, D(7))$, $(7, D(17))$ of a graph $G$, we obtain the permutation $(17\ 15\ 7)$ which, when applied to $D$, yields $D* = (1\ 7\ 20\ 15\ 12\ 17)$. Since $D(15) = 1$, $D(7) = 12$, $D(17) = 20$, each of the given vertices appears in $D*$. Continuing, we first work with the value d($a_1$, MIN(M)($a_1$, 1)) obtained in $MIN_1$. We assume that it is negative. We now use theorem 3.1 to obtain a cycle having a negative value such that each partial sum is negative. We note that we can never choose an arc of form $(a, D^{-1}(a))$: When we apply $D$ to the terminal vertex of such an arc, we obtain $(a, a)$, a loop. Thus, the new permutation being constructed would not be a derangement. Continuing, we obtain $(a_1, D^{-1}(MIN(M)(a_1, 1)))$ which is given the value d($a_1, D^{-1}(j_k)$). Letting $D^{-1}(j^k) = a_k$, we next obtain $(a_k, D^{-1}(\text{MIN(M)}(a_k + 1))) = (a_k, D^{-1}(j_m)) = (a_k, a_m)$. If the sum of the first two terms is negative, we continue. Otherwise, we stop. If we stop, we check to see if either the sum $(a_i\ a_k) + (a_k\ a_i)$ or $(a_i\ a_k) + (a_k\ a_m) + (a_m\ a_i)$ is negative. If at least one of them is negative, we choose the one with the smaller value and save it for use later on. If the path $[a_i, a_k, a_m]$ has a negative value, we continue with this procedure. Checking the value of $(a_i + a_k + a_m)$, if its value is smaller than that of the first two cycles, we substitute it for the 2-cycle saved earlier. At some point in this procedure, one of the following occurs:

(1)    The partial sum of the arcs of the path becomes positive.

(2)    The partial sum of the arcs remains negative but the path repeats a vertex therefore forming a cycle out of a subpath of arcs.

Assume (1) occurs with the path $[a_i\ a_k\ a_m\ ...\ a_z]$. In this case, we save the cycle, $C_1$, with smallest negative value. If (2) occurs, let our path be $[a_i\ a_k\ a_m\ ...\ a_p\ a_q\ a_r\ ...\ a_y\ a_m]$.

As in the first case, we have obtained the cycle, $C_2$, of smallest value from among the elements of the set

$$\{(a_i\ a_k), (a_i\ a_k\ a_m), ..., (a_i\ a_k\ a_m\ ...\ a_p), ..., (a_i\ a_k\ a_m\ ...\ a_p\ ...\ a_y)\}.$$

Without loss of generality, suppose the two cycles $(a_i\ a_k)$ and $(a_m\ a_p\ a_y)$ both have negative values. Call the product of the two cycles, $P_3$. Whichever of the three permutations, $C_1$, $C_2$, $P_3$ has the smallest (negative) value, gives us our best negative permutation to be applied to D to obtain $D_1$. On the other hand, suppose that only one of the cycles, $C_i$, has a negative value. Then we choose whichever of $C_1$ and $P_3$ has the smaller negative value. We follow the same procedure starting off with each value $i$ of MIN(M)$(i, j)$ for $j = 2, 3, ...[\log n]+1$. Once we have obtained the smallest (negative) cycle, say $C^*$, from among all cycles tested, we multiply it by $D$ to obtain $DC^* = D_1$, a *derangement*, i.e., a permutation which moves every point of $V = \{1,2,...,n\}$. We continue with this procedure until we no longer can obtain a negative cycle. At this stage, we obtain all negative entries from MIN(M)$(i, j)$ $(i = 1,2, ..., [\log n] + 1; j = 1,2, ... [\log n] + 1)$. We then apply the algorithm given above to each entry. If we obtain no further negative cycle, we go on to Phase 2.

**PHASE 2.**

Before going on, a *determining vertex* of a cycle of at most value $L$ is an initial vertex, $d$, of a path $P$ that traverses $C$ such that every subpath is of value at most $L$. From the corollary to theorem 3.1, at least one such vertex, $d$, always exists in a cycle $C$ of value $M$. Phase 2 uses a modified form of the Floyd-Warshall Algorithm to obtain an optimal solution to the Assignment Problem. For short, we call this algorithm the F-W n-vs algorithm, i.e., the Floyd-Warshall negatively valued subpath algorithm. Given $D_i$, we construct $D_i^{-1}M^-$ where $M$ is the value matrix. $D_i^{-1}$ permutes the columns of $M$. Thus, if $D_i = (a_1\ a_2\ ...\ a_n)$, $a_2 = D_i(a_1)$, etc. . Thus, if $(u\ D_i(u))$ is an arc of $D_i$, it is mapped onto $(u\ u)$ in $D_i^{-1}M^-$. It follows that the arcs of $D_i$ are mapped onto the diagonal elements of $D_i^{-1}M^-$. Therefore, the initial values of the diagonal elements of $D_i^{-1}M^-$ are all zero. Furthermore, if - using theorem 3.1 – we obtain a negative cycle, $C_1$, in $D_i^{-1}M^-$, $D_iC_1$ has a smaller value than $D_i$. The value of $C_1$ is represented by an underlined negative value lying on a diagonal entry. By permuting the columns of $D_i$, we obtain $C_1^{-1}D_i^{-1}M^-$. We continue this process until we no longer can obtain a negative cycle. This indicates that the last permutation obtained is $\sigma_{APOPT}$. We now discuss the details of the algorithm. We first obtain $\text{MIN}(D_i^{-1}M^-)$. It consists of sorting the row numbers of each row of $M$ from entries of smallest value to those of largest value. Using $\text{MIN}(D_i^{-1}M^-)$, we apply the F-W algorithm to each negative entry of $D_i^{-1}M^-$. As we leave a negative entry, we write its value in italics. If we are able to extend it to a negative path, we underline the value of the new path. We extend an underlined path by using $\text{MIN}(D_i^{-1}M^-)$ to obtain negatively-valued entries from $D_i^{-1}M^-$. By the nature of the F-W algorithm, we go from column 1 to column 2, column 3, ... , column $n$. An *iteration* of the algorithm is completed if we have gone through each column from 1 through *n* once. As an iteration proceeds, we place each of the underlined paths in a balanced, binary search tree, *NEGPATHS*. We delete underlined paths that become italicized; we add paths that have just become underlined. If an extension of a negative path has a smaller negative value than an italicized entry, we replace the italicized entry by a new underlined path. On the other hand, if we choose an element of *NEGPATHS* that isn't smaller than the current entry, we delete it from *NEGPATHS*. In general, if a new underlined path is at an entry *(a b)* where $b > a$, it will atill be underlined at the end of an iteration. As we go through the *j*-*th* iteration, we place the row and column numbers of each new underlined path entry in a matrix, $P_{in}$. If we were using the normal F-W algorithm and a negative cycle exists, we would require only one iteration to obtain the cycle, because we don't impose the conditions of theorem 3.1 upon it. Hopefully, using the modified algorithm, we will be able to lower the running time, $O(n^3)$, of the F-W algorithm. During an iteration, we stop at the moment we've obtained a negative cycle, $C_1$. We then construct $MIN(C_1^{-1} D_i^{-1} M^-)$. We continue the algorithm until we've either obtained a negative cycle or we have not obtained any underlined paths during an entire iteration.

We use $MIN(D_i^{-1} M^-)$ to extend an underlined path as well as path matrices, $P_{in}$ *(i = 1,2, … )* to keep track of the path corresponding to each negative entry *(i, j)*. We use theorem 3.1, starting with a (negative) determining vertex, to construct negative paths. Unlike the F-W algorithm, we don't piece together paths to construct new smaller paths. We add only one arc from $D_i^{-1} M^-$ at a time to extend a negative path. When we have gone through columns 1, 2, … , n, we have completed one *iteration* of this portion of the algorithm. Our goal is to obtain a negative cycle if one exists. If we are able to extend a path *(a b)* to *(a b)(b c)= (a c)*, we underline the new entry (*a c*). Once we have gone through the entry – either extending the path whose value is given at *(a b)* or not doing so - we write the entry at *(a b)* in italics. If no negative cycle is obtained using only *n* columns (as is the case in the normal F-W algorithm), we continue the algorithm into a second iteration using $D_i^{-1}(n)$. As we proceed with the iteration, after underlining the value of an extended negative path in $D_i^{-1}M^-(n)$, we place the value in the ordered list *NEGPATHS*. If a negative path is extended further, we delete the old entry and place the new value in *NEGPATHS*. Also, as we proceed, we use the entries in $P_0$ to construct $P_1$. If we cannot obtain a negative cycle after using $2n$ columns, we construct $D_i^{-1}M^-(2n)$ and $P_2$. We continue this procedure of iteration until either we obtain a negative cycle or we are unable to extend any negative path further. Once a negative cycle, s, is obtained, we use it to construct $D_{i+1}^{-1}M^-$. When we no longer can obtain a negative cycle, the last derangement obtained must be the smallest-valued derangement obtainable, namely, an optimal solution to the Assignment Problem, $\sigma_{APOPT}$.

## PHASE 3.

Before going on, we denote the smallest-valued tour obtainable using the modified F-W algorithm as $\sigma_{FWTSPOPT}$.

**Step 1.** We use the last matrix obtained - $\sigma_{APOPT}^{-1}M^-$ - at the start of the algorithm. Our ultimate goal is to obtain an optimal tour of $M$, $\sigma_{TSPOPT}$. Our immediate one is to obtain $\sigma_{FWTSPOPT}$.

**Step 2.** A *tour* is an $n$ cycle each of whose arcs has a finite value. We first use Karp's patching algorithm to obtain a tour, $T_{UPPERBOUND}$. We then check to see if any of the permutations obtained earlier using the F-W algorithm yield derangements smaller in value than $T_{UPPERBOUND}$ are tours. If they are, we call the smallest of them $T_{UPPERBOUND}$  In the worked-out example 3.3, we've used the difference of the total values of $\sigma_{APOPT}$ and $T_{UPPERBOUND}$. $|T_{UPPERBOUND}| = 161$. Thus, $m_0 = |T_{UPPERBOUND}| - |\sigma_{APOPT}| = 6$. Using the corollary to theorem 3.1, we can assume that each value in $D_3^{-1}M^-$ that belongs to a cycle is no larger than $m_0 - 1 = 5$.

*Note.* As in PHASE 1, we use MIN(M) to obtain possible candidates for arcs to extend paths. Also, as we proceed with PHASE 2, we could place the arcs and their values on a balanced, binary search tree as

we proceed, sorting each branch. Once we have gone through the first $n$ columns, we could use the search tree for choice of paths. If it is possible to permute branches of the tree to simulate permuting columns of a matrix, we might be able to save running time.

**Step 3a.** Given a tour with an initial upper bound, $T_{UPPERBOUND}$, we first check to see if any of the derangements obtained earlier that were of smaller value than $T_{UPPERBOUND}$ were tours. If so, we choose the smallest-valued tour and denote it by $T_{UPPERBOUND}$. To obtain $\sigma_{FWTSPOPT}$, we use the modified F-W algorithm to obtain all lowest-valued non-negative cycles obtainable smaller-valued than $T_{UPPERBOUND}$. Since we no longer have to worry about negative cycles, non-negative cycles may be obtainable provided $|T_{UPPERBOUND}| \neq |\sigma_{APOPT}|$.

*Note.* In general, to simplify our problem, We first obtain an ordinal value and its inverse for each point of $T_{UPPERBOUND}$. Then, if we obtain a set of pair-wise disjoint cycles, $p$, we shrink $T_{UPPERBOUND}$ so that it contains only the points of $p$ before applying $p$ to it.

**Theorem 3.3** Using the modified F-W algorithm, we can obtain a negative cycle, $C$, of length $L$ in $D_i^{-1}M^-$ in at most $L$ iterations.

**Proof.** We can obtain at least one arc in each iteration. Thus, we can obtain $C$ in at most $L$ iterations. In general, when using our specialized version of the F-W algorithm, we choose the first negative cycle we obtain. The reason for this becomes evident when we consider the following theorem:

**Theorem 3.4** Let $D_i^{-1}M^-$ contain both positive and negative entries. Suppose that $P$ is a path obtained using the modified F-W. algorithm with initial vertex $b_1$. Assume that $C = (a_1\ a_2\ ...\ a_m)$ is a negatively-valued cycle where $a_1 \neq b_1$. Then if $C$ is constructed as arcs are extended to $P$, thus making $P$ a non-simple path, our algorithm obtains $C$ as an independent cycle before it is constructed in $P$.

**Proof.** Using the modified F-W algorithm, no matter which vertex $a_i$, of $C$ is first obtained in $P$, it requires the same number of columns to retrun to $a_i$ and make $P$ a non-simple path. Since $b_l \neq a_l$, It requires at least one column to go from $b_i$ to $a_l$, concluding the proof.

In the following example, we construct $P$ and $C$.

**Example 3.3**

$P = [1^{-20}3^{5}7^{-5}13^{12}15^{1}19^{3}20^{-18}18^{1}14^{1}6^{3}7^{-5}]$,

$C = (20\ 18\ 14\ 6\ 7\ 13\ 15\ 19)$.

In what follows, Roman numerals represent numbers of iterations.

| P | | C | |
|---|---|---|---|
| **I.** | [1 3] + 3 | **I.** | [20 18] + 18 |
| | [1 7] + 4 | | ------------------------ |
| | [1 13] + 6 | **II.** | [20 14] + 16 |
| | [1 15] + 2 | | ------------------------ |
| | [1 19] + 4 | **III.** | [20 6] + 12 |
| | [1 20] + 1 | | [20 7] + 1 |
| ------------------------ | | | [20 13] + 6 |
| **II.** | [1 18] + 18 | | [20 15] + 2 |
| ------------------------ | | | |
| **III.** | [1 14] + 16 | | [20 19] + 4 |
| ------------------------ | | | [20 20] + 1 |
| **IV.** | [1 6] + 12 | | ------------------------ |
| | [1 7] + 1 | | + 60 |
| ------------------------ | | | |
| | + 67 | | |

**Theorem 3.5** *Let c be any real number, $S_a = \{ a_i \mid i = 1, 2, \ldots, n\}$, a set of real numbers in increasing order of value. For $i = 1, 2, \ldots, n$, let $b_i = a_i + c$. Then $S_b = \{ b_i \mid i = 1, 2, \ldots, n\}$ preserves the ordering of $S_a$.*

**Proof**. The theorem merely states that adding a fixed number to a set ordered according to the value of its elements retains the same ordering as that of the original set.

*Comment*. This is very useful when we are dealing with entries of the value matrix which have not been changed during the algorithm. However, as the algorithm goes on, if we have entry (i, j) and d(i, j), we must go through all values of j = j' where the (i, j')-th entry's value in the current $\sigma_a^{-1} M^-(k)$ is less than the value of the (i, j')-th entry in $D^{-1} M^-$ (The terminology will become clearer in the examples given.)

**Theorem 3.6** *Suppose that $m_j$ is the last upper bound obtained in Step 3a. If $\sigma_{TSPOPT} = \sigma_{APOPT} s$, let $C = (a\ a_1\ a_2\ \ldots\ a_s\ b)$ be an arbitrary disjoint cycle of s where a is a determining vertex of C. Then either (1) C will be obtained in Step 3a or (2) there always exists a cycle obtained in Step 3a that is of the form $C' = (a\ b_1\ b_2\ \ldots\ b_r\ b)$ and of value no greater than C. Furthermore, a is a determining vertex of C'.*

**Proof.** $\sigma_{TSPOPT} = \sigma_{APOPT} s$. Here s is a permutation in $S_n$. The value of s is no greater than $m_j$ since we're assuming that $|\sigma_{TSP}| \leq |\sigma_{FWTSPOPT}|$. Let $C = (a\ a_1\ a_2\ \ldots\ a_s\ b)$ be a disjoint cycle of s. We can obtain s (as well as *all* cycles of value less than $m_i$) by starting a path in $D_i^{-1} M^-$ whose first vertex is a determining vertex. Let $d(a, a_1, a_2, \ldots, a_s, b) = m$ and $m \cup d(b, a) = m' \leq m_j$. Since a is a determining vertex of C, there must exist a smallest-valued simple path, P, in $D_i^{-1} M^-$ from a to b of value no greater than m. This path is obtainable using our modified F-W algorithm where the first vertex is always a determining vertex. Each minimal path obtained using our modified F-W algorithm must have an initial vertex that is a *determining* vertex. Since the initial vertex of P, a, is a determining vertex, $C' = P \cup (b\ a)$.

**Corollary 3.6** *If we can obtain no cycle in Step 3a, then $\sigma_{TSPOPT} = \sigma_{FWTSPOPT}$.*

**Proof**. The corollary follows directly from the fact that if we could obtain a cycle (a … b) in Step 3b, then a corresponding cycle must exist in Step 3a.

Our algorithm always obtains $\sigma_{TSPOPT}$ in (at most) Step 3b because there always exists some permutation, $p$, such that $(\sigma_{APOPT})p = \sigma_{TSPOPT}$. Whether this procedure can be done in every case in a finite amount of time is an open question.

**Example 3.4**

PHASE 1.

Let D = (1 2 3 4 5 6 7 8), while M is the following distance matrix:

M

|   | 1 | 2 | 3 | 4 | 5 | 6 | 7 | 8 |   |
|---|---|---|---|---|---|---|---|---|---|
| 1 | ∞ | 23 | 99 | 17 | 12 | 99 | 18 | 24 | 1 |
| 2 | 43 | ∞ | 2 | 73 | 15 | 100 | 53 | 28 | 2 |
| 3 | 1 | 84 | ∞ | 19 | 53 | 68 | 44 | 34 | 3 |
| 4 | 89 | 41 | 45 | ∞ | 40 | 71 | 79 | 51 | 4 |
| 5 | 83 | 62 | 94 | 88 | ∞ | 36 | 6 | 50 | 5 |
| 6 | 61 | 62 | 98 | 50 | 29 | ∞ | 52 | 40 | 6 |
| 7 | 50 | 21 | 53 | 68 | 39 | 26 | ∞ | 25 | 7 |
| 8 | 16 | 42 | 61 | 54 | 81 | 34 | 92 | ∞ | 8 |
|   | 1 | 2 | 3 | 4 | 5 | 6 | 7 | 8 |   |

In the next table, MIN(M), we give the column numbers of the entries in each row of M arranged in increasing order of value.

MIN(M)

|   | 1 | 2 | 3 | 4 | 5 | 6 | 7 | 8 |   |
|---|---|---|---|---|---|---|---|---|---|
| 1 | 5 | 4 | 7 | 2 | 8 | 3 | 6 | 1 | 1 |
| 2 | 3 | 5 | 8 | 1 | 7 | 4 | 6 | 2 | 2 |
| 3 | 1 | 4 | 8 | 7 | 5 | 6 | 2 | 3 | 3 |
| 4 | 5 | 2 | 3 | 8 | 6 | 7 | 1 | 4 | 4 |
| 5 | 7 | 6 | 8 | 2 | 1 | 4 | 3 | 5 | 5 |
| 6 | 5 | 8 | 4 | 7 | 1 | 2 | 3 | 6 | 6 |
| 7 | 2 | 8 | 6 | 5 | 1 | 3 | 4 | 7 | 7 |
| 8 | 1 | 6 | 2 | 4 | 3 | 5 | 7 | 8 | 8 |
|   | 1 | 2 | 3 | 4 | 5 | 6 | 7 | 8 |   |

MIN(M)(1, 1) = (1, 5): 12;  MIN(M)(2, 1) = (2, 3): 2;

MIN(M)(3, 1) = (3, 1): 1;  MIN(M)(4, 1) = (4, 5): 40;

MIN(M)(5, 1) = (5, 7): 6;  MIN(M)(6, 1) = (6, 5): 29;

MIN(M)(7, 1) = (7, 2): 31;  MIN(M)(8, 1) = (8, 1) = 16.

We now obtain the values of DIFF(i) (i = 1, 2, ..., 8).

*DIFF*(1) = d(1, 5) - d(1, 2) = -11.

*DIFF*(2) = d(2, 3) - d(2, 3) = 0.

*DIFF*(3) = d(3, 1) - d(3, 4) = -18.

*DIFF*(4) = d(4. 5) - d(4, 5) = 0.

*DIFF*(5) = d(5, 7) - d(5, 6) = -30.

*DIFF*(6) = d(6, 5) - d(6, 7) = -23.

*DIFF*(7) = d(7, 2) - d(7, 8) = -4.

*DIFF*(8) = d(8, 1) - d(8, 1)  = 0.

min { $DIFF(i)$ | $i$ = 1, 2, ..., 8} = DIFF(5) = -30.

We will start Phase 1 with 5.

Before doing so, we write $D$ and $D^{-1}$ in ROW FORM.

$D = $

| -11 | 0 | -18 | 0 | -30 | -23 | -4 | 0 |
|---|---|---|---|---|---|---|---|
| 1 | 2 | 3 | 4 | 5 | 6 | 7 | 8 |
| 2 | 3 | 4 | 5 | 6 | 7 | 8 | 1 |

$D^{-1} = $

| 1 | 2 | 3 | 4 | 5 | 6 | 7 | 8 |
|---|---|---|---|---|---|---|---|
| 8 | 1 | 2 | 3 | 4 | 5 | 6 | 7 |

TRIAL 1.

    (5, 7) → (5, 6)   - 30

    (6, 5) → (6, 4)   - 23

    (4, 5) → (4, 4)

(4, 5) is an arc of $D$.

    (4, MIN(M)(4, 2)) = (4, 2)    $d(4, 2) - d(4, 5) = 41 - 40 = 1$

    (4, 2) → (4, 1)

    (1, 5) → (1, 4)   -11

$P$ = [5, 6, 4, 1, 4].

(1, 5) is obtained from the arc (1, 6) in M.

$d(1, 6) - d(1, 2) = 99 - 23 = 76$.

(4, 5) is obtained from (4, 6).

d(4, 6) - d(4, 5) = 71 - 40 = 31

(6, 5) is obtained from the loop (6, 6).

$s_{11} = (5^{-30} 6^{-23} 4^1 1^{76})$, $s_{112} = (5^{-30} 6^{-23} 4^{31})$.

$s_{113} = (5\ 6)$ is not admissible since $Ds_{113}$ is not a derangement.

The value of $s_{11}$ is -22.

$s_{12} = (5\ 6)(4\ 1)$. As shown earlier, (5 6) is not admissible.

(1, 4) is obtained from the arc (1, 5) of M.

d(1, 5) - d(1, 2) = 12 - 23 = -11.

Thus, since the second cycle has a negative value, we rename it

$s_{12} = (4^1 1^{-11})$. A reasonable question is why we go to the bother of

saving it: $s_{11}$ has a smaller value than $s_{12}$. Our reason for doing so is that in Phase 3, we want to obtain an 8-cycle with as small as possible value.

We may thus want to check to see if $Ds_{11}$ is an 8-cycle. In any event,

$s_{112} = (5\ 6\ 4)$ has the smallest negative value so far: -22.

TRIAL 2. MIN(M)(5, 2) = (5, 6).

    (5, 6) → (5, 5)

(5, 6) is an arc of D. We thus can go no farther with 5.

Therefore, $s_1 = (5\ 6\ 4)$. $Ds_1 = D_1$. We obtained $D_1$ by exchanging those rows of D that have as initial arcs 5, 6, 4, with rows (5, 7), (6, 5), (4, 6), respectively. We must construct DIFF values for these arcs. All other DIFF values remain the same.

    -11  0  -18 -31  0   0  -4   0

$$D_1 = \begin{matrix} 1 & 2 & 3 & 4 & 5 & 6 & 7 & 8 \\ 2 & 3 & 4 & 6 & 7 & 5 & 8 & 1 \end{matrix}$$

$$D_1^{-1} = \begin{matrix} 1 & 2 & 3 & 4 & 5 & 6 & 7 & 8 \\ 8 & 1 & 2 & 3 & 6 & 4 & 5 & 7 \end{matrix}$$

d(5, 7) - d(5, 7) = 0,  d(6, 5) - d(6, 5) = 0,  d(4, 5) - d(4, 6) = -31.

We start with 4.

TRIAL 1.  MIN(M)(4, 1) = (4, 5).

$\quad$ (4, 5) → (4, 6)  $\quad$ -31

$\quad$ (6, 5) → (6, 6)

(6, 5) is an arc of $D_1$.

$\quad$ (6, MIN(6, 2)) = (6, 8)  $\quad$ d(6, 8) - d(6, 5) = 11.

$\quad$ (6, 8) → (6, 7)  $\quad$ 11

$\quad$ (7, 2) → (7, 1)  $\quad$ -4

$\quad$ (1, 5) → (1, 6)  $\quad$ -11

P = [4, 6, 7, 1, 6].

The arc (1, 4) is derived from the arc (1, 6) of M.

d(1, 6) - d(1, 2) = 99 - 23 = 76.

d(7, 6) - d(7, 8) = 26 - 25 = 1.

(6, 6) is a loop. Thus, the cycle (4 6) is not an admissible cycle since

$D_1(1, 6)$ is not a derangement.

$s_{II} = (4^{-31} 6^{11} 7^{-4} 1^{76})$, $s_{III} = (4^{-31} 6^{11} 7^1)$.

The value of $s_{II}$ is non-negative. The value of $s_{III}$ is -19.

The only possibility for $s_{21}$ is $(6^{11}7^{-4}1^{-11})$. The value of the latter cycle is

-4. We keep $s_{21} = (6\ 7\ 1)$ in our bag of negative cycles.

TRIAL 2.  MIN(M)(4, 2) = (4, 2).  d(4, 2) - d(4, 6) = 41 - 71 = -30.

  (4, 2) → (4, 1) -30

  (1, 5) → (1, 6) -11

  (6, 5) → (6, 6)

(6, 5) is an arc of $D_1$.

  (6, MIN(M)(6, 2)) = (6, 8) d(6, 8) - d(6, 5) = 40 - 29 = 11.

  (6, 8) → (6, 7) 11

  (7, 2) → (7, 1) -4

$P = [4, 1, 6, 7, 1]$.

(7, 4) is derived from (7, 6) in M.

d(7, 6) - d(7, 8) = 26 - 25 = 1.

(6, 4) is derived from the loop (6, 6).

(1, 4) is derived from (1, 6).

d(1, 6) - d(1, 2) = 99 - 23 = 76.

$s_{21} = (4^{-30}1^{-11}6^{11}7^1)$, $s_{211} = (4^{-30}1^{76})$.

$s_{21}$ has a negative value: -29. $s_{211}$ has a non-negative value.

$s_{11}$ has the smallest negative value thus far.

(7, 1) is derived from the arc (7, 2) of M.

d(7, 2) - d(7, 8) = 21 - 25 = -4.

Then $s_{22} = (1^{-11}6^{11}7^{-4})$. The value of $s_{22}$ is -4. We therefore keep it in our bag of negative cycles.

TRIAL 3.  MIN(M)(4, 3) = (4, 3).  d(4, 3) - d(4, 6) = 45 - 71 = -26.

  (4, 3) → (4, 2) -26

$(2, 3) \rightarrow (2, 2)$.

$(2, 3)$ is an arc of $D_1$.

$(2, MIN(M)(2, 2)) = (2, 5)$.   $d(2, 5) - d(2, 3) = 15 - 2 = 13$.

$(2, 5) \rightarrow (2, 6)$   13

$(6, 5) \rightarrow (6, 6)$.

$(6, 5)$ is an arc of $D_1$.

$(6, MIN(M)(6, 2)) = (6, 8)$.   $d(6, 8) - d(6, 5) = 11$

$(6, 8) \rightarrow (6, 7)$   11

$(7, 2) \rightarrow (7, 1)$   -4

$(1, 5) \rightarrow (1, 6)$

$P = [4, 2, 6, 7, 1, 6]$.

$(1, 4)$ is derived from $(1, 6)$.   $d(1, 6) - d(1, 2) = 76$.

$d(7, 6) - d(7, 8) = 1$

$(6, 6)$ is a loop.

$d(2, 6) - d(2, 3) = 100 - 2 = 98$.

$P = [4^{-26} 2^{13} 6^{11} 7^{-4} 1^{76}]$.

Perusing the values we obtained above, the only possible negative cycle obtainable is $s_{313} = (4\ 2\ 6\ 7)$ which has a value of -1.

$s_{32} = (6\ 7\ 1) = s_{22}$; its value is -1.

Thus, $s_2 = s_{21} = (4\ 1\ 6\ 7)$.

In order to obtain $D_2$, we respectively exchange rows $(4, 6)$, $(1, 2)$,

$(6, 5)$, $(7, 8)$ of $D_1$ with rows $(4, 2)$, $(1, 5)$, $(6, 8)$, $(7, 6)$.

Our new rows have the following *DIFF* values:

$d(4, 5) - d(4, 2) = -1$,  $d(1, 5) - d(1, 5) = 0$,  $d(6, 5) - d(6, 8) = -11$,

d(7, 2) - d(7, 6) = -5.

$$D_2 = \begin{array}{ccccccc} 0 & 0 & -18 & -1 & 0 & -11 & -5 & 0 \\ 1 & 2 & 3 & 4 & 5 & 6 & 7 & 8 \\ \\ 5 & 3 & 4 & 2 & 7 & 8 & 6 & 1 \end{array}$$

$$D_2^{-1} = \begin{array}{cccccccc} 1 & 2 & 3 & 4 & 5 & 6 & 7 & 8 \\ \\ 8 & 4 & 2 & 3 & 1 & 7 & 5 & 6 \end{array}$$

We start with 3.

TRIAL 1.

    (3, 1) → (3, 18)  -18

    (8, 1) → (8, 8)

(8, 1) is an arc of $D_2$.

    (8, MIN(M)(8, 2)) = (8, 6).   d(8, 6) - d(8, 1) = 18.

The path is non-negative. If we can go no further using the smallest negative *DIFF,* we are allowed to apply the procedure to at most the next [*log n*] smallest negative *DIFF* values.

Thus, we start with 6.

TRIAL 1.

    (6, 5) → (6, 1)  -11

    (1, 5) → (1, 1)

(1, 5) is an arc of $D_2$.

(1, MIN(M)(1, 2)) = (1, 4).   d(1, 4) - d(1,5) = 5.

(1, 4) → (1, 3)   5

(3, 1) → (3, 8)   -18

(8, 1) → (8, 8)

(8, 1) is an arc of $D_2$.

(8, MIN(M)(8, 2)) = (8, 6).   d(8, 6) - d(8, 1) = 18.

(8, 6) → (8, 7)   18

(7, 2) → (7, 4)   -5

(4, 5) → (4, 1)

$P = [6\ 1\ 3\ 8\ 7\ 4\ 1]$.

(4, 6) is derived from (4, 8).

d(4, 8) - d(4, 2) = 51 - 41 = 10.

d(7, 8) - d(7, 6) = 25 - 26 = -1.

(8, 8) is a loop.

d(3, 8) - d(3, 4) = 34 - 19 = 15.

d(1, 8) - d(1, 5) = 24 - 12 = 12.

We thus obtain the following cycles from $P$:

$s_{l1} = (6^{-11}1^5 3^{-18} 8^{18} 7^{-5} 4^{10})$, $s_{l11} = (6^{-11} 1^5 3^{-18} 8^{18} 7^{-1})$, $s_{l13} = (6^{-11} 1^5 3^{15})$,
$s_{l14} = (6^{-11} 1^{12})$.

$s_{l1}$ has a value of -1.  $s_{l11}$ has a value of -7.  The other two cycles have non-negative values.

(4, 1) is derived from the arc (4, 5) in M.

d(4, 5) - d(4, 2) = -1.

$s_{l2} = (1^5 3^{-18} 8^{18} 7^{-5} 4^{-1})$. The value of $s_{l2}$ is -1.

TRIAL 2. MIN(M)(6, 2)) = (6, 8).   d(6, 8) - d(6, 8) = 0.

We can proceed no further using 6.

Since we have obtained a negative cycle using 6, we need not go to 7.

Thus, let $s_3 = s_{III} = (6\ 1\ 3\ 8\ 7)$. $D_3 = D_2 s_3$.

Rows (6, 8), (1, 5), (3, 4), *, 1), (7, 6) are replaced, respectively, by rows (6, 5), (1, 4), (3, 1), (8, 6), 7, 8) to obtain $D_3$.

The following are the new *DIFF* values:

d(6, 5) - d(6, 5) = 0.

d(1, 5) - d(1, 4) = -5.

d(3, 1) - d(3, 1) = 0.

d(8, 1) - d(8, 6) = -18.

d(7, 2) - d(7, 8) = -4.

$$D_3 = \begin{matrix} -5 & 0 & 0 & -1 & 0 & 0 & -4 & -15 \\ 1 & 2 & 3 & 4 & 5 & 6 & 7 & 8 \\ 4 & 3 & 1 & 2 & 7 & 5 & 8 & 6 \end{matrix}$$

$$D_3^{-1} = \begin{matrix} 1 & 2 & 3 & 4 & 5 & 6 & 7 & 8 \\ 3 & 4 & 2 & 1 & 6 & 8 & 5 & 7 \end{matrix}$$

We now start with 8.

TRIAL 1.

(8, 1) → (8, 3)  -18

(3, 1) → (3, 3).

(3, 1) is an arc of $D_3$.

  (3, MIN(M)(3, 1) = (3, 4).

  (3, 4) → (3, 1)  18

The path $P = [8, 3, 1]$ is non-negative.

TRIAL 2. MIN(M)(8, 2) = (8, 6).  d(8, 6) - d(8, 6) = 0. Thus, (8, 6) is an arc of $D_3$. We can go no further with 8.

We now start with 1.

TRIAL 1.

  (1, 5) → (1, 6)  -5

  (6, 5) → (6, 6).

(6, 5) is an arc of $D_3$.

  (6, MIN(M)(6, 2)) = (6, 8).  d(6, 8) - d(6, 5) = 11.

$P = [1, 6, 8]$ is non-negative.

(6, 1) is derived from (6, 4) of M. d(6, 4) - d(6, 5) = 50 - 29 = 21.

The cycle (1 6) has a non-negative value.

TRIAL 2. MIN(M)(1, 2) = (1, 4). d(1, 4) - d(1, 4) = 0.

(1, 4) is an arc of $D_3$. We can go no farther with 1.

The third smallest negative *DIFF* value is -4. It occurs at 7.

We start with 7.

TRIAL 1.

  (7, 2) → (7, 4)  -4

  (4, 5) → (4, 6)  -1

  (6, 5) → (6, 6).

(6, 5) is an arc of $D_3$.

  (6, MIN(M)(6, 2)) = (6, 8).  d(6, 8) - d(6, 5) = 11.

We can go no further since P = [7, 4, 6, 8] is non-negative.

(6, 7) is derived from (6, 8) of M.

d(6, 8) - d(6, 5) = 11.

d(4, 8) - d(4, 5) = 51 - 40 = 11.

$s_{II} = (7^{-4}4^{-1}6^{11})$, $s_{III} = (7^{-4}4^{11})$.

Both cycles have non-negative values.

TRIAL 2.  MIN(M)(7, 2)) = (7, 8).    d(7, 8) - d(7, 8) = 0.

(7, 8) is an arc of $D_3$. We thus conclude Phase 1 with

$D_3 = (1^{17}4^{41}2^{2}3^{1})(5^{6}7^{25}8^{34}6^{29})$. We next go to Phase 2 to try to obtain a derangement with a smaller value than $D_3$.

PHASE 2.

$$D_3^{-1}M$$

|   | 4 | 3 | 1 | 2 | 7 | 5 | 8 | 6 |   |
|---|---|---|---|---|---|---|---|---|---|
|   | 1 | 2 | 3 | 4 | 5 | 6 | 7 | 8 |   |
| 1 | 17 | 99 | ∞ | 23 | 16 | 12 | 24 | 99 | 1 |
| 2 | 73 | 2 | 43 | ∞ | 53 | 15 | 28 | 100 | 2 |
| 3 | 19 | ∞ | 1 | 84 | 44 | 53 | 34 | 68 | 3 |
| 4 | ∞ | 45 | 89 | 41 | 79 | 40 | 51 | 71 | 4 |
| 5 | 88 | 94 | 83 | 62 | 6 | ∞ | 50 | 36 | 5 |
| 6 | 50 | 98 | 61 | 62 | 52 | 29 | 40 | ∞ | 6 |
| 7 | 68 | 53 | 50 | 21 | ∞ | 39 | 25 | 26 | 7 |
| 8 | 54 | 61 | 16 | 42 | 92 | 81 | ∞ | 34 | 8 |
|   | 1 | 2 | 3 | 4 | 5 | 6 | 7 | 8 |   |

$$D_3^{-1}M^-$$

|   |   | 4 | 3 | 1 | 2 | 7 | 5 | 8 | 6 |   |
|---|---|---|---|---|---|---|---|---|---|---|
|   |   | 1 | 2 | 3 | 4 | 5 | 6 | 7 | 8 |   |
| 1 | 0 | 82 | ∞ | 6 | 1 | -5 | 7 | 82 | 1 |
| 2 | 71 | 0 | 41 | ∞ | 51 | 13 | 26 | 98 | 2 |
| 3 | 18 | ∞ | 0 | 83 | 43 | 52 | 33 | 67 | 3 |
| 4 | ∞ | 4 | 48 | 0 | 38 | -1 | 10 | 30 | 4 |
| 5 | 82 | 88 | 77 | 56 | 0 | ∞ | 44 | 30 | 5 |
| 6 | 21 | 69 | 32 | 33 | 23 | 0 | 11 | ∞ | 6 |
| 7 | 43 | 28 | 25 | -4 | ∞ | 14 | 0 | 1 | 7 |
| 8 | 20 | 27 | -18 | 8 | 58 | 47 | ∞ | 0 | 8 |
|   |   | 1 | 2 | 3 | 4 | 5 | 6 | 7 | 8 |   |

$j = 1, 2$.

There are no negative entries in the first and second columns.

$j = 3$.

The only negative entry in the third column is in row 8.

$d(8, 3) = -18$. MIN(M)(3, 1) = 1. $\sigma_4^{-1}(1) = 3$. (3, 3) is a loop.

MIN(M)(3, 2) = 4. $\sigma_4^{-1}(4) = 1$. $d(3, 1) = 18$. $d(8, 3) + d(3, 1) \geq 0$.

MIN(M)(3, $i$), $i = 3,4,5,6,7,8$ are each greater than 18. Thus,

$d(8, 3) + d(3, i) > 0$, $i = 3,4,5,6,7,8$.

$j = 4$.

The only negative entry in column four occurs in row 7.

$d(7, 4) = -4$. MIN(M)(4, 1) = 5. $\sigma_4^{-1}(5) = 6$. $d(4, 6) = -1$.

$d(7, 4) + d(4, 6) = -5. < d(7, 6)$. Thus, we place $-5$ in $(7, 6)$.

$MIN(M)(4, 2) = 2$. $\sigma_4^{-1}(2) = 4$. $(4, 4)$ is a loop.

$MIN(M)(4, 3) = 3$. $\sigma_4^{-1}(3) = 2$. $d(4, 2) = 4$. $d(7, 4) + d(4, 6) \geq 0$. The $MIN(M)(4, i)$ for $i = 4,5,6,7,8$ yield values in $\sigma_4^{-1} M^-$ greater than 4. Thus, $d(7, 4) + d(4, i) > 0$ for $i = 4,5,6,7,8$.

$$P_4(4)$$

|   | 1 | 2 | 3 | 4 | 5 | 6 | 7 | 8 |   |
|---|---|---|---|---|---|---|---|---|---|
| 1 |   |   |   |   |   |   |   |   | 1 |
| 2 |   |   |   |   |   |   |   |   | 2 |
| 3 |   |   |   |   |   |   |   |   | 3 |
| 4 |   |   |   |   |   |   |   |   | 4 |
| 5 |   |   |   |   |   |   |   |   | 5 |
| 6 |   |   |   |   |   |   |   |   | 6 |
| 7 |   |   |   |   |   | 1 |   |   | 7 |
| 8 |   |   |   |   |   |   |   |   | 8 |
|   | 1 | 2 | 3 | 4 | 5 | 6 | 7 | 8 |   |

$j = 5$.

There are no negative entries in column 5.

$j = 6$.

There are three negative entries in column 6: rows 1, 4 and 7.

$d(1, 6) = -5$. $MIN(M)(6, 1) = 5$. $\sigma_4^{-1}(5) = 6$. $(6, 6)$ is a loop.

$MIN(M)(6, 2) = 8$. $\sigma_4^{-1}(8) = 7$. $d(6, 7) = 11$. $d(1, 6) + d(6, 7) > 0$.

For $i = 3,4,5,6,7,8$, $MIN(6, i) > MIN(M)(6, 2)$. Thus,

d(1, 6) + d(6, *i*) > 0.

d(4, 6) = - 1. MIN(M)(6, 2) = 8. $\sigma_4^{-1}(8) = 7$. d(6, 7) = 11.

d(4, 6) + d(6, 7) > 0. The same is true with

MIN(M)(6, *i*), *i* = 3,4,5,6,7,8.

d(7, 6) = - 5. MIN(M)(6, 2) = 8. $\sigma_4^{-1}(8) = 7$. d(6, 7) = 11.

d(4, 6) + d(6, 7) > 0. The same is true with

MIN(M)(6, *i*), *i* = 3,4,5,6,7,8.

*j* = 7, 8.

There are no negative entries in columns 7 and 8.

We have only one negative path in $D_3^{-1} M^-(8)$: [7, 4, 6] with a value of -5. This path can't be extended in $D_3^{-1} M^-(16)$. It follows that $\sigma_{APOPT} = D_3$. The total sum of its values is 155. We now apply Karp's algorithm to obtain $T_{UPPERBOUND}$, the 2-cycle (3 8) yields

$T_{UPPERBOUND}$ = $D_3(3\ 8)$ = (1 5 7 8 6 4 2 3) : 171. We now check to see if we can use cycles obtained earlier to obtain derangements that are tours are tours having values less than 171. We first check those cycles that have negative values: $s_{11}$ and $s_{21}$ both have values of -1.

$D_2 s_{11} = (1^{17} 4^{51} 8^{34} 6^{29} 5^6 7^{21} 2^2 3^1)$ : 161. This is the smallest-valued tour that we can obtain. We thus define $D_2 s_{11}$ as $T_{UPPERBOUND}$.

### PHASE 3.

We now apply Step 3a to obtain an approximation to $\sigma_{TSPOPT}$. Using only entries of $D_3^{-1} M^-$, we must be able to obtain a positive cycle of value no greater than 5. We use the modified F-W algorithm.

$$D_3^{-1}M^-$$

|   | 4 | 3 | 1 | 2 | 7 | 5 | 8 | 6 |   |
|---|---|---|---|---|---|---|---|---|---|
|   | 1 | 2 | 3 | 4 | 5 | 6 | 7 | 8 |   |
| 1 | 0 | 82 | ∞ | 6 | 1 | -5 | 7 | 82 | 1 |
| 2 | 71 | 0 | 41 | ∞ | 51 | 13 | 26 | 98 | 2 |
| 3 | 18 | ∞ | 0 | 83 | 43 | 52 | 33 | 67 | 3 |
| 4 | ∞ | 4 | 48 | 0 | 38 | -1 | 10 | 30 | 4 |
| 5 | 82 | 88 | 77 | 56 | 0 | ∞ | 44 | 30 | 5 |
| 6 | 21 | 69 | 32 | 33 | 23 | 0 | 11 | ∞ | 6 |
| 7 | 43 | 28 | 25 | -4 | ∞ | 14 | 0 | 1 | 7 |
| 8 | 20 | 27 | -18 | 8 | 58 | 47 | ∞ | 0 | 8 |
|   | 1 | 2 | 3 | 4 | 5 | 6 | 7 | 8 |   |

$j$=3.

(8, 3)(3, 1): 0.

$j$=4.

(7, 4)(4, 2): 0;  (7, 4)(4, 6): -5.

$j$=8.

(7, 8)(8, 3): -17.

$$D_3^{-1}M^-(8)$$

|   | 4 | 3 | 1 | 2 | 7 | 5 | 8 | 6 |   |
|---|---|---|---|---|---|---|---|---|---|
|   | 1 | 2 | 3 | 4 | 5 | 6 | 7 | 8 |   |
| 1 | 0 | 82 | ∞ | 6 | 1 | -5 | 7 | 82 | 1 |
| 2 | 71 | 0 | 41 | ∞ | 51 | 13 | 26 | 98 | 2 |
| 3 | 18 | ∞ | 0 | 83 | 43 | 52 | 33 | 67 | 3 |
| 4 | ∞ | 4 | 48 | 0 | 38 | -1 | 10 | 30 | 4 |
| 5 | 82 | 88 | 77 | 56 | 0 | ∞ | 44 | 30 | 5 |
| 6 | 21 | 69 | 32 | 33 | 23 | 0 | 11 | ∞ | 6 |
| 7 | 43 | <u>0</u> | <u>-17</u> | -4 | ∞ | <u>-5</u> | 0 | 1 | 7 |
| 8 | <u>0</u> | 27 | *-18* | 8 | 58 | 47 | ∞ | 0 | 8 |
|   | 1 | 2 | 3 | 4 | 5 | 6 | 7 | 8 |   |

$j=1$.

(8, 1)(1, 6): -5;  (8, 1)(1, 5): 1.

$j=3$.

(7, 3)(3, 1): 1.

$$D_3^{-1}M^-\ (16)$$

|   | 4 | 3 | 1 | 2 | 7 | 5 | 8 | 6 |   |
|---|---|---|---|---|---|---|---|---|---|
|   | 1 | 2 | 3 | 4 | 5 | 6 | 7 | 8 |   |
| 1 | 0 | 82 | ∞ | 6 | 1 | -5 | 7 | 82 | 1 |
| 2 | 71 | 0 | 41 | ∞ | 51 | 13 | 26 | 98 | 2 |
| 3 | 18 | ∞ | 0 | 83 | 43 | 52 | 33 | 67 | 3 |
| 4 | ∞ | 4 | 48 | 0 | 38 | -1 | 10 | 30 | 4 |
| 5 | 82 | 88 | 77 | 56 | 0 | ∞ | 44 | 30 | 5 |
| 6 | 21 | 69 | 32 | 33 | 23 | 0 | 11 | ∞ | 6 |
| 7 | *1* | 0 | *-17* | *-4* | ∞ | -5 | 0 | 1 | 7 |
| 8 | *0* | 27 | *-18* | 8 | *1* | *-5* | ∞ | 0 | 8 |
|   | 1 | 2 | 3 | 4 | 5 | 6 | 7 | 8 |   |

$j=1$.

(7, 1)(1, 5): 2

(7, 5) cannot be extended to a path of value less than 6.

It follows that $\sigma_{FWTSPOPT} = T_{UPPERBOUND} = (1^{17}4^{51}8^{34}6^{29}5^6 7^{21}2^2 3^1)$.

Furthermore, from the corollary to theorem 3.6, since we are unable to obtain a positive cycle of value less than 6, $\sigma_{TSPOPT} = \sigma_{FWTSPOPT}$.

**Example 3.5**

M

|    | 1  | 2  | 3  | 4  | 5  | 6  | 7  | 8  | 9  | 10 | 11 | 12 | 13 | 14 | 15 | 16 | 17 | 18 | 19 | 20 |    |
|----|----|----|----|----|----|----|----|----|----|----|----|----|----|----|----|----|----|----|----|----|----|
| 1  | ∞  | 88 | 72 | 97 | 14 | 38 | 9  | 59 | 39 | 46 | 52 | 50 | 29 | 17 | 48 | 65 | 2  | 72 | 86 | 65 | 1  |
| 2  | 80 | ∞  | 28 | 72 | 67 | 18 | 99 | 20 | 58 | 66 | 24 | 76 | 32 | 45 | 11 | 62 | 54 | 62 | 25 | 45 | 2  |
| 3  | 55 | 40 | ∞  | 82 | 34 | 26 | 73 | 76 | 97 | 17 | 13 | 33 | 23 | 94 | 76 | 87 | 56 | 32 | 74 | 81 | 3  |
| 4  | 29 | 76 | 77 | ∞  | 92 | 63 | 94 | 88 | 87 | 18 | 38 | 59 | 94 | 62 | 33 | 18 | 9  | 67 | 93 | 31 | 4  |
| 5  | 17 | 32 | 60 | 80 | ∞  | 49 | 21 | 64 | 77 | 54 | 41 | 18 | 91 | 4  | 35 | 29 | 10 | 19 | 99 | 35 | 5  |
| 6  | 21 | 19 | 31 | 75 | 19 | ∞  | 98 | 94 | 72 | 40 | 30 | 43 | 48 | 18 | 94 | 82 | 69 | 70 | 22 | 71 | 6  |
| 7  | 1  | 99 | 82 | 68 | 4  | 35 | ∞  | 74 | 28 | 44 | 81 | 59 | 24 | 33 | 9  | 10 | 91 | 95 | 58 | 89 | 7  |
| 8  | 1  | 49 | 71 | 49 | 52 | 30 | 52 | ∞  | 59 | 92 | 25 | 48 | 56 | 65 | 34 | 59 | 24 | 78 | 67 | 70 | 8  |
| 9  | 92 | 99 | 49 | 6  | 63 | 79 | 66 | 19 | ∞  | 63 | 90 | 24 | 92 | 21 | 4  | 43 | 77 | 68 | 84 | 66 | 9  |
| 10 | 46 | 65 | 2  | 4  | 80 | 54 | 92 | 92 | 72 | ∞  | 34 | 10 | 86 | 63 | 63 | 40 | 73 | 99 | 85 | 5  | 10 |
| 11 | 11 | 14 | 81 | 18 | 91 | 37 | 46 | 51 | 32 | 23 | ∞  | 90 | 30 | 85 | 18 | 66 | 54 | 85 | 31 | 19 | 11 |
| 12 | 43 | 1  | 52 | 64 | 6  | 39 | 79 | 89 | 39 | 44 | 77 | ∞  | 78 | 42 | 47 | 62 | 68 | 65 | 25 | 7  | 12 |
| 13 | 19 | 82 | 93 | 76 | 20 | 80 | 15 | 81 | 15 | 87 | 45 | 67 | ∞  | 54 | 9  | 92 | 16 | 67 | 7  | 4  | 13 |
| 14 | 84 | 73 | 89 | 90 | 56 | 96 | 31 | 52 | 28 | 26 | 23 | 54 | 19 | ∞  | 91 | 37 | 6  | 95 | 59 | 26 | 14 |
| 15 | 84 | 53 | 65 | 42 | 65 | 54 | 62 | 81 | 90 | 80 | 98 | 52 | 59 | 44 | ∞  | 18 | 79 | 39 | 50 | 91 | 15 |
| 16 | 13 | 5  | 77 | 60 | 81 | 5  | 88 | 17 | 58 | 48 | 62 | 12 | 59 | 20 | 48 | ∞  | 69 | 61 | 20 | 57 | 16 |
| 17 | 46 | 48 | 25 | 59 | 8  | 83 | 83 | 24 | 28 | 1  | 19 | 75 | 17 | 28 | 82 | 75 | ∞  | 71 | 31 | 6  | 17 |
| 18 | 31 | 50 | 84 | 98 | 26 | 80 | 67 | 51 | 83 | 80 | 82 | 90 | 42 | 9  | 3  | 26 | 68 | ∞  | 51 | 41 | 18 |
| 19 | 53 | 89 | 6  | 44 | 58 | 48 | 26 | 17 | 64 | 88 | 63 | 63 | 87 | 61 | 42 | 57 | 32 | 59 | ∞  | 68 | 19 |
| 20 | 43 | 6  | 73 | 51 | 49 | 52 | 14 | 56 | 35 | 18 | 24 | 80 | 65 | 47 | 6  | 55 | 91 | 85 | 84 | ∞  | 20 |
|    | 1  | 2  | 3  | 4  | 5  | 6  | 7  | 8  | 9  | 10 | 11 | 12 | 13 | 14 | 15 | 16 | 17 | 18 | 19 | 20 |    |

MIN(M)

|    | 1  | 2  | 3  | 4  | 5  | 6  | 7  | 8  | 9  | 10 | 11 | 12 | 13 | 14 | 15 | 16 | 17 | 18 | 19 | 20 |    |
|----|----|----|----|----|----|----|----|----|----|----|----|----|----|----|----|----|----|----|----|----|----|
| 1  | 17 | 7  | 5  | 14 | 13 | 6  | 9  | 10 | 15 | 12 | 11 | 8  | 16 | 20 | 3  | 18 | 19 | 2  | 4  | 1  | 1  |
| 2  | 15 | 6  | 8  | 11 | 19 | 3  | 13 | 14 | 20 | 17 | 9  | 16 | 18 | 10 | 5  | 4  | 12 | 1  | 7  | 2  | 2  |
| 3  | 11 | 10 | 13 | 6  | 18 | 12 | 5  | 2  | 1  | 17 | 7  | 19 | 8  | 15 | 20 | 4  | 16 | 14 | 9  | 3  | 3  |
| 4  | 17 | 10 | 16 | 1  | 20 | 15 | 11 | 12 | 14 | 6  | 18 | 2  | 3  | 9  | 8  | 5  | 7  | 13 | 19 | 4  | 4  |
| 5  | 14 | 17 | 1  | 12 | 18 | 7  | 16 | 2  | 15 | 20 | 11 | 6  | 10 | 3  | 8  | 9  | 4  | 13 | 19 | 5  | 5  |
| 6  | 14 | 2  | 5  | 1  | 19 | 11 | 3  | 10 | 12 | 13 | 17 | 18 | 20 | 9  | 4  | 16 | 8  | 15 | 7  | 6  | 6  |
| 7  | 1  | 5  | 9  | 16 | 13 | 9  | 6  | 14 | 10 | 19 | 12 | 4  | 8  | 11 | 3  | 20 | 17 | 18 | 2  | 7  | 7  |
| 8  | 1  | 17 | 11 | 6  | 15 | 12 | 2  | 4  | 5  | 7  | 13 | 9  | 16 | 14 | 19 | 20 | 3  | 18 | 10 | 8  | 8  |
| 9  | 15 | 4  | 8  | 14 | 12 | 16 | 3  | 5  | 10 | 7  | 20 | 18 | 17 | 6  | 19 | 11 | 1  | 13 | 2  | 9  | 9  |
| 10 | 3  | 4  | 20 | 12 | 11 | 16 | 1  | 6  | 14 | 15 | 2  | 9  | 17 | 5  | 19 | 13 | 7  | 8  | 18 | 10 | 10 |
| 11 | 1  | 2  | 4  | 15 | 20 | 10 | 13 | 19 | 9  | 6  | 7  | 8  | 17 | 16 | 3  | 14 | 18 | 12 | 5  | 11 | 11 |
| 12 | 2  | 5  | 20 | 19 | 6  | 9  | 14 | 1  | 10 | 15 | 3  | 16 | 4  | 17 | 18 | 11 | 13 | 7  | 8  | 12 | 12 |
| 13 | 20 | 19 | 15 | 7  | 9  | 17 | 1  | 5  | 11 | 14 | 12 | 18 | 4  | 6  | 8  | 2  | 10 | 16 | 3  | 13 | 13 |
| 14 | 17 | 13 | 11 | 10 | 20 | 9  | 7  | 16 | 8  | 12 | 5  | 19 | 2  | 1  | 3  | 4  | 15 | 18 | 6  | 14 | 14 |
| 15 | 16 | 18 | 4  | 14 | 12 | 2  | 6  | 13 | 19 | 7  | 3  | 5  | 17 | 10 | 8  | 1  | 9  | 20 | 11 | 15 | 15 |
| 16 | 2  | 6  | 12 | 1  | 8  | 14 | 19 | 10 | 15 | 20 | 9  | 13 | 4  | 18 | 11 | 17 | 3  | 5  | 7  | 16 | 16 |
| 17 | 10 | 5  | 20 | 13 | 11 | 8  | 3  | 9  | 14 | 19 | 1  | 2  | 4  | 18 | 12 | 16 | 15 | 6  | 7  | 17 | 17 |
| 18 | 15 | 14 | 5  | 16 | 1  | 20 | 13 | 2  | 8  | 19 | 7  | 17 | 6  | 10 | 11 | 9  | 3  | 12 | 4  | 18 | 18 |
| 19 | 3  | 8  | 7  | 17 | 15 | 4  | 6  | 1  | 16 | 5  | 18 | 14 | 11 | 12 | 9  | 20 | 13 | 10 | 2  | 19 | 19 |
| 20 | 2  | 15 | 7  | 10 | 11 | 9  | 1  | 14 | 5  | 4  | 6  | 16 | 8  | 13 | 3  | 12 | 19 | 18 | 17 | 20 | 20 |
|    | 1  | 2  | 3  | 4  | 5  | 6  | 7  | 8  | 9  | 10 | 11 | 12 | 13 | 14 | 15 | 16 | 17 | 18 | 19 | 20 |    |

$$D = (1\ 2\ 3\ 4\ 5\ 6\ 7\ 8\ 9\ 10\ 11\ 12\ 13\ 14\ 15\ 16\ 17\ 18\ 19\ 20).$$

From M, the values of the arcs of $D$ are

|  | 88 | 28 | 82 | 92 | 49 | 98 | 74 | 59 | 63 | 34 |  |
|---|---|---|---|---|---|---|---|---|---|---|---|
| (1 | 2 | 3 | 4 | 5 | 6 | 7 | 8 | 9 | 10 | | |
| 90 | 78 | 54 | 91 | 18 | 69 | 71 | 51 | 68 | 43 | | |
| 11 | 12 | 13 | 14 | 15 | 16 | 17 | 18 | 19 | 20). | | |

We now use MIN(M) to obtain the smallest value in each row of M.

MIN(M)(1, 1) = 17: d(1, 17) = 2. MIN(M)(2, 1) = 15: d(2, 15) = 11.

MIN(M)(3, 1) = 11: d(3, 11) = 13. MIN(M)(4, 1) = 17: d(4, 17) = 9. IN(M)(5, 1) = 14: d(5, 14) = 18.

MIN(M)(6, 1) = 14: d(6, 14) = 18.

MIN(M)(7, 1) = 1: d(7, 1) = 1. MIN(M)(8, 1) = 1: d(8, 1) = 1.

MIN(M)(9, 1) = 15: d(9, 15) = 4. MIN(M)(10, 1) = 3: d(10, 3) = 2.

MIN(M)(11, 1) = 1: d(11, 1) = 11. MIN(M)(12, 1) = 2: d(12, 2) = 1.

MIN(M)(13, 1) = 20: d(13, 20) = 4. MIN(M)(14, 1) = 17: d(14, 17) = 6.

MIN(M)(15, 1) = 16: d(15, 16) = 18. MIN(M)(16, 1) = 2: d(16, 2) = 5.

MIN(M)(17, 1) = 10: d(17, 10) = 1. MIN(M)(18, 1) = 15: d(18, 15) = 3.

MIN(M)(19, 1) = 3: d(19, 3) = 6. MIN(M)(20, 1) = 2: d(20, 2) = 6.

For $j = 1, 2, \ldots, 19$, we now obtain d($j$, MIN(M)($j$,1)) - d($j$, $j$+1);

also, d(20, 2) - d(20, 1). We then replace the values given the arcs of $D$

by the corresponding differences we've just constructed. We now have

$$D = \begin{pmatrix} -86 & -17 & -69 & -83 & -45 & -80 & -73 & -58 & -59 & -32 \\ 1 & 2 & 3 & 4 & 5 & 6 & 7 & 8 & 9 & 10 \\ -79 & -77 & -50 & -85 & 0 & -64 & -70 & -48 & -62 & -37 \\ 11 & 12 & 13 & 14 & 15 & 16 & 17 & 18 & 19 & 20 \end{pmatrix}$$

We next construct a permutation consisting of at most two cycles which when applied to $D$ will reduce the total sum of its arcs in M. To start the construction, we choose an arc where the absolute value of its difference is greatest. We therefore start our first cycle with 1. Before going on, we introduce an economical way to write permutations, their inverses and their differences as we proceed from derangement to derangement.

We first write two ordered lists of the numbers from 1 through n.

We then first enter the number corresponding to $D_i(a)$ underneath a in the first list. We then write the number preceding $a$ underneath $a$ in the second list. Simply speaking, if have

$$\begin{matrix} a \\ b \end{matrix}$$

occurring in a column in the first list, we have

$$\begin{matrix} b \\ a \end{matrix}$$

occurring in a column in the second list. If we don't know the predecessor, $p$, of the first point $a$ in a cycle of a permutation, it will always occur in the second row of the permutation when we reach the arc $(p, D_i(p))$ where $D_i(p) = a$. Finally, we place the value of $DIFF(a)$ above $a$ in the second list. When we follow this practice for $D$, we must write all points of both $D$ and $D^{-1}$ because the permutation $D$ is itself a derangement. However, when constructing further permutations, we need only change those columns in which the point $a$ is contained in a permutation.

We now rewrite $D$ and $D^{-1}$ in the form mentioned earlier. For simplicity, we call it the *row form* of D.

1 2 3 4 5 6 7 8  9 10 11 12 13 14 15 16 17 18 19 20

2 3 4 5 6 7 8 9 10 11 12 13 14 15 16 17 18 19 20  1

$1^{-86}$  $2^{-17}$  $3^{-69}$  $4^{-83}$  $5^{-45}$  $6^{-80}$  $7^{-73}$  $8^{-58}$  $9^{-59}$  $10^{-32}$  $11^{-79}$  $12^{\mathbf{-77}}$  $13^{-50}$  $14^{-85}$

20   1   2   3   4   5   6   7    8    9    10    11    12    13

$15^{0}$  $16^{-64}$  $17^{-70}$  $18^{-48}$  $19^{-62}$  $20^{-37}$

14   15   16    17    18    19

As we construct our lists, we keep track of the size of the values of *DIFF* obtained. We start off with *DIFF*(1). If and when we find a smaller value of *DIFF*, say at point *a*, we replace *DIFF*(1) by *DIFF(a)*, replacing the preceding value of *DIFF*. Thus, when we are finished we have obtained both the point at which the smallest value of *DIFF* occurs as well as the value of *DIFF* there.

Using this procedure, we see that the smallest value of *DIFF* occurs when $a = 1$. *DIFF(1)* = -86.

The rule that we adopted in Phase 1 was that we try the first $|log\ n|$ entries MIN(M)($i, j$) in row *i* where $|log\ n|$ represents the positive integer closest to $log\ n$ ($j = 1, 2, ... ,|log\ n|$). log(20) = 2.996. Thus, in this case, we will pick the first three possibilities {MIN(M)($i, j$) | $j = 1,2,3$}. There are two more things of which we must be aware. The first is that if an arc (*a, b*) of $D_i$ is chosen when constructing a permutation ,we will obtain (*a, a*) for the corresponding arc of the permutation. If (*a, b*) comes from MIN(M)(*a, 1*), then we must choose MIN(M)(*a, 2*). Otherwise, we will choose MIN(M)(*a, 1*). Secondly, in testing for cycles, the last arc of a permutation cycle is often chosen without knowing the arc of the value matrix from which it is obtained. Occasionally, it may occur that the arc chosen actually comes from a loop. i.e., an arc of the form (*a, a*). Thus, the current derangement

would lead to a permutation which is not a derangement. To avoid this occurrence, we must check the terminal arc of any possible which we have not obtained directly from an arc obtained from MIN(M). To do this, given the permutation arc, say $(a, b)$, we obtain from the first two rows of the row form, $(a, D_i(b))$. If $D_i(b) \neq a$, then we have no difficulty. Otherwise, we go on to test the next possible terminal arc of a cycle.

TRIAL 1.

MIN(M)(1,1) = 17. d(1, 17) = 2. In order to construct a corresponding arc of the cycle, we must obtain $D^{-1}(17) = 16$. Thus, the first arc of the first cycle is (1, 16). We now continue with the vertex 16.

MIN(M)(16, 1) = 2. Thus, the second arc of M chosen is (16, 2) where d(16, 2) = 5. $D^{-1}(2) = 1$. Thus, the first cycle we've obtained is (1 16). Adding up the values corresponding to 1 and 16, respectively, if we applied the 2-cycle (1 16) to $D$, we would reduce the total sum of arcs of $D$ in M by 86 + 64 = 150.

TRIAL 2.

MIN(M)(1,2) = 7. d(1, 7) = 9. Thus, instead of - 86, we now have - 79. $D^{-1}(7) = 6$. Thus, the first arc of the first cycle of a permutation is (1 6). We will now proceed with the construction of the first cycle. We must always keep in mind that the path we construct must always have a negative sum of values. Otherwise, the permutation won't reduce the total value of the derangement to which it will be applied.

    (1  7)    (1 6)   - 79

    (6  14)   (6 13)   - 80

    (13  20) (13 19) - 50

    (19  3)   (19 2)   - 62

    (2  15)  (2 14)   - 17

    (14  17) (14  16) - 85

    (16  2)   (16 1)   - 64.

The cycle (1 6 13 19 2 14 16) reduces the sum of the values of the arcs of $D$ in M by 437.

TRIAL 3.

MIN(M)(1, 3) = 5. d(1, 5) = 14. Thus, the value of the arc of $D$ out of 1 is reduced by 74.

$\quad\quad\quad\quad$ (1  5)  (1  4) - 74

$\quad\quad\quad\quad$ (4  17) (4  16) - 83

$\quad\quad\quad\quad$ (16  2) (16  1) - 64.

The cycle (1  4  16) reduces the sum of the values of the arcs of $D$ by 221.

It follows that our best choice is the second one. Let

$s_1$ = (1 6 13 19 2 14 16). Next, we give the row form of the new permutation, $D_1$, obtained by applying

$s_1$ to $D$. Asterisks have been placed above each point of $D$ which has been moved.

*  *$\quad\quad\quad$*$\quad\quad\quad\quad\quad\quad\quad$*  *$\quad\quad$*$\quad\quad\quad$*

1  2  3  4  5$\quad$6  7  8$\quad$9  10  11  12  13  14  15  16  17  18  19  20

7 15  4  5  6  14  8  9  10  11  12  13  20  17  16   2  18  19   3   1

1$\quad$2$\quad$3$^{-69}$ 4$^{-83}$ 5$^{-45}$ 6$\quad$7$^{-73}$ 8$^{-58}$ 9$^{-59}$ 10$^{-32}$ 11$^{-79}$ 12$^{-77}$ 13  14  15$^0$  16

20$\quad$16  19$\quad$3$\quad$4$\quad$5  1$\quad$7$\quad$8$\quad$9$\quad$10$\quad$11$\quad$12  6  2  15

17$^{-70}$ 18$^{-48}$  19  20$^{-37}$

14$\quad$17$\quad$18  13

We now must obtain the $DIFF$ values for points under asterisks.

MIN(M)(1, 1) = 17 $\Rightarrow$ $DIFF$(1) = d(1, 17) - d(1, 17) = 0.

MIN(M)(2, 1) = 15 $\Rightarrow$ $DIFF$(2) = d(2, 15) - d(2, 15) = 0.

MIN(M)(6, 1) = 14 $\Rightarrow$ $DIFF$(6) = d(6, 14) - d(6, 14) = 0 .

MIN(M)(13, 1) = 20 $\Rightarrow$ $DIFF$(13) = d(13, 20) - d(13, 20) = 0.

MIN(M)(14, 1) = 17 $\Rightarrow$ $DIFF$(14) = d(14, 17) - d(14, 17) = 0.

MIN(M)(16, 1) = 2 ⇒ $DIFF(16) = d(16, 2) - d(16, 2) = 0$.

MIN(M)(19, 1) = 3 ⇒ $DIFF(19) = d(19, 3) - d(19, 3) = 0$.

Thus, we can complete the row form above.

1 2 3 4 5  6 7 8  9 10 11 12 13 14 15 16 17 18 19 20

$D_1$

 7 15 4 5 6 14 8 9 10 11 12 13 20 17 16  2 18 19  3  1

$1^0$ $2^0$ $3^{-69}$ $4^{-83}$ $5^{-45}$ $6^0$ $7^{-73}$ $8^{-58}$ $9^{-59}$ $10^{-32}$ $11^{-79}$ $12^{-77}$ $13^0$ $14^0$

$D_1^{-1}$

20 16 19   3   4   5   1   7   8      9     10    11    12    6

$15^0$ $16^0$ $17^{-70}$ $18^{-48}$ $19^0$ $20^{-37}$

   2   15   14     17    18   13

The smallest value of $DIFF$ occurs when $a = 4$.

TRIAL 1.

MIN(M)(4, 1) = 17.

$$(4, 17) \to (4, 14)$$

$$(14, 17) \to (14, 14)$$

(14, 17) is an arc of $D_1$. We next try to obtain MIN(M)(14, 2) = 13.

$$(14, 13) \to (14, 12)$$

$$(12, 2) \to (12, 16)$$

$$(16, 2) \to (16, 16)$$

(16, 2) is an arc of $D_1$. MIN(M)(16, 2) = 6.

$$(16, 6) \rightarrow (16, 5)$$

$$(5, 14) \rightarrow (5, 6)$$

$$(6, 14) \rightarrow (6, 6)$$

$(6, 14)$ is an arc of $D_1$.

$$(6, 2) \rightarrow (6, 16).$$

We thus have a path $P = [4, 14, 12, 16, 5, 6, 16]$.

Let $s_{11} = (4\ 14\ 12\ 16\ 5\ 6)$.

We check to see the arc of $D_1$ associated with (6 4). Going to the first two rows of *ROW FORM*, we obtain (6 5). In those cases in which we didn't encounter an arc of $D_1$ in constructing an arc of $s_{11}$, we can use the values associated with the points of the cycle in *ROW FORM*.

Thus, $s_{11} = (4^{-83}\ 14\ 12^{-77}\ 16\ 5^{-45})$.

We now find the values of 14, 16 and 6. For each vertex $a$, we have the arc in M, obtained from MIN(M)$(a, 2)$. We then subtract the arc in $D_1$ with the same initial vertex.

$d(14, 13) - d(14, 17) = 19 - 6 = 13$.

$d(16, 6) - d(16, 2) = 5 - 5 = 0$.

$d(6, 5) - d(6, 14) = 19 - 18 = 1$

Thus, $s_{11} = (4^{-83}\ 14^{13}\ 12^{-77}\ 16^0\ 5^{-45}\ 6^1)$.

$s_{11}$ has the value -191. Thus, when applied to $D_1$, it reduces the value of $D_1$ by 191.

We now consider the permutations obtained by including all partial paths of $s_{11}$ which begin with 4.

Let

$s_{111} = (4\ 14)$, $s_{112} = (4\ 14\ 12)$, $s_{113} = (4\ 14\ 12\ 16)$,

$s_{114} = (4\ 14\ 12\ 16\ 5)$.

Since 4 is mapped into 5 by $D_1$, in order to obtain the value of each of these permutations, we have to obtain the respective values of (14, 5), (12, 5), (16, 5). Since (5,5) is a loop, we need not consider $s_{114}$.

d(14, 5) = 56, d(12, 5) = 6, d(16, 5) = 81.

We next subtract the value of each arc of form $(a, D_1(a))$

($a$ = 14, 12, 5) from the respective values obtained above to obtain the value of each arc in its respective permutation.

d(14, 5) - d(14, 17) = 56 - 6 = 50.

d(12, 5) - d(14, 17) = 6 - 6 = 0.

d(16, 5) - d(16, 2) = 81 - 5 = 76.

Thus, $s_{111} = (4^{-83}\ 14^{50})$, $s_{112} = (4^{-83}\ 14^{13}\ 12^{0})$,

$s_{113} = (4^{-83}\ 14^{13}\ 12^{-77}\ 16^{76})$. Thus, the respective values of these permutations are -33, -70 and -71.

We now consider the permutation obtained by using the path $P$ to construct a permutation of two cycles.

$s_{12} = (4\ 14\ 12)(16\ 5\ 6)$.

We already have the values of all of the points of $s_{12}$ except for the arc (6 16). We now obtain the value associated with 6 in s $s_{12}$.

d(6, 16) → d(6, 15) = 94

d(6, 15) - d(6, 14) = 94 - 18 = 76.

$s_{12} = (4^{-83}\ 14^{13}\ 12^{0})(16^{0}\ 5^{-45}\ 6^{76})$.

The first cycle of $s_{12}$ is $s_{112}$. The value of the second cycle is non-negative so we needn't include it.

TRIAL 2. MIN(M)(4, 2) = 10.

$$(4, 10) \rightarrow (4, 9)$$

$$(9, 15) \rightarrow (9, 2)$$

$$(2, 15) \rightarrow (2, 2)$$

Thus, (2, 15) is an arc of $D_1$. MIN(M)(2,2) = 6

$$(2, 6) \rightarrow (2, 5)$$

$$(5, 14) \rightarrow (5, 6)$$

$$(6, 14) \rightarrow (6, 6).$$

$(6, 14)$ is an arc of $D_1$. MIN(M)(6, 2) = 2.

$$(6, 2) \rightarrow (6, 16)$$

$$(16, 2) \rightarrow (16, 16).$$

$(16, 2)$ is an arc of $D_1$.

$$(16, 6) \rightarrow (16, 5).$$

Let $P = [4, 9, 2, 5, 6\ 16, 5]$.

$s_{21} = (4\ 9\ 2\ 5\ 6\ 16)$.

We must obtain the values associated with the initial vertices 2, 6, 16. The other arcs were all obtained using terminal vertices in the first column of MIN(M). Thus, we can read them off directly from *ROW FORM*. We must compute $d(16, D_1(4)) = d(16, 5)$ rather than using $d(16, 6)$.

$d(2, 6) - d(2, 15) = 18 - 11 = 7$.

$d(6, 2) - d(6, 14) = 19 - 18 = 1$.

$d(16, 5) - d(16, 2) = 81 - 5 = 76$.

Thus, $s_{21} = (4^{-83}\ 9^{-59}\ 2^7\ 5^{-45}\ 6^1\ 16^{76})$. The value of $s_{21}$ is -103.

$s_{211} = (4\ 9)$, $s_{212} = (4\ 9\ 2)$, $s_{213} = (4\ 9\ 2\ 5)$, $s_{214} = (4\ 9\ 2\ 5\ 6)$.

We must thus obtain $d(9, 5), d(2, 5), d(6, 5)$. $s_{213}$ need not be considered since the arc (5, 4) of the permutation maps into the loop (5, 5).

$d(9, 5) = 63$, $d(2, 5) = 67$, $d(6, 5) = 19$.

The respective values in $s_{211}, s_{212}, s_{214}$ are

$d(9, 5) - d(9, 10) = 63 - 63 = 0$, $d(2, 5) - d(2, 15) = 67 - 11 = 56$,

$d(6, 5) - d(6, 14) = 19 - 18 = 1$.

$s_{211} = (4^{-83}\ 9^0)$, $(4^{-83}\ 9^{-59}\ 2^{56})$, $s_{214} = (4^{-83}\ 9^{-59}\ 2^7\ 5^{-45}\ 6^1)$.

The respective values of these permutations are -83, -86, -179.

We now construct $s_{22}$ from $P$.

$s_{22} = (4\ 9\ 2)(5\ 6\ 16)$. We have already obtained the values for all vertices except 16 previously.

$d(16, D_1(5)) = d(16, 6)$. $d(16, 6) - d(16, 2) = 5 - 5 = 0$.

Thus, we have

$s_{22} = (4^{-93}\ 9^{-59}\ 2^{56})(5^{-45}\ 6^1\ 16^0)$.

The value of $s_{22}$ is -140.

TRIAL 3. MIN(M)(4, 3) = 16.

$$(4, 16) \rightarrow (4, 15)$$

$$(15, 16) \rightarrow (15, 15)$$

(15, 16) is an arc of $D_1$.

$$(15, 18) \rightarrow (15, 17)$$

$$(17, 10) \rightarrow (17, 9)$$

$$(9, 15) \rightarrow (9, 2)$$

$$(2, 15) \rightarrow (2, 2)$$

(2, 15) is an arc of $D_1$.

$$(2, 6) \rightarrow (2, 5)$$

$$(5, 14) \rightarrow (5, 6)$$

$$(6, 14) \rightarrow (6, 6)$$

(6, 14) is an arc of $D_1$.

$$(6, 2) \rightarrow (6, 16)$$

$$(16, 2) \rightarrow (16, 16)$$

(16, 2) is an arc of $D_1$.

$$(16, 6) \rightarrow (16, 5)$$

$P = [4, 15, 17, 9, 2, 5, 6, 16, 5]$.

15, 2, 6 and 16 were the initial vertices of arcs in $D_1$. Thus, we had to use arcs obtained from the second column of MIN(M). We can immediately read off the values of the remaining vertices from *ROW FORM*..

$s_{31} = (4^{-83}\ 15\ 17^{-70}\ 9^{-59}\ 2\ 5^{-45}\ 6\ 16)$.

d(15, 18) - d(15, 16) = 39 - 18 = 21.

d(2, 6) - d(2, 15) = 18 - 11 = 7.

d(6, 2) - d(6, 14) = 19 - 18 = 1.

d(16, 5) - d(16, 2) = 81 - 5 = 76.

We thus obtain $s_{31} = (4^{-83}\ 15^{21}\ 17^{-70}\ 9^{-59}\ 2^7\ 5^{-45}\ 6^1\ 16^{76})$. The value of $s_{31}$ is -152.

$s_{311} = (4\ 15)$. $s_{312} = (4\ 15\ 17)$. $s_{313} = (4\ 15\ 17\ 9)$. $s_{314} = (4\ 15\ 17\ 9\ 2)$. Since (5, 4) maps into (5, 5), we need not consider (4 15 17 9 2 5).

$s_{315} = (4\ 15\ 17\ 9\ 2\ 5\ 6)$.

d(15, 5) - d(15, 16) = 65 - 18 = 47.

d(17, 5) - d(17, 18) = 8 - 71 = -63.

d(9, 5) - d(9, 10) = 63 - 63 = 0.

d(2, 5) - d(2, 15) = 67 - 11 = 56.

d(6, 5) - d(6, 14) = 19 - 18 = 1.

d(16, 5) - d(16, 2) = 81 - 5 = 76.

Thus, $s_{311} = (4^{-83}\ 15^{47})$, $s_{312} = (4^{-73}15^{21}17^{-73})$, $s_{313} = (4^{-73}15^{21}17^{-70}9^0)$,

$s_{314} = (4^{-83}15^{21}17^{-70}9^{-59}2^{56})$, $s_{315} = (4^{-83}15^{21}17^{-70}9^{-59}2^7 5^{-45}6^1)$.

The respective values of these permutations are: -36, -125, -122, -135,

From P, $s_{32} = (4\ 15\ 17\ 9\ 2)(5\ 6\ 16)$. The first cycle is $s_{314}$. The second was obtained in trial 2 as the second factor of $s_{22}$. It follows that the value of $s_{32}$ is (-135) + (-44) = -179.

After checking all cycles obtained, $s_{315}$ has the smallest value: -228.

Thus, $s_2 = (4\ 15\ 17\ 9\ 2\ 5\ 6)$. In order to obtain $D_2$ in row form, we need only change the points 4, 15, 17, 9, 2, 5, 6. In constructing permutations, we have obtained the changes in the arcs of $D_1$. However, if we wish, we can use the row form of $D_1$ to obtain new columns in $D_2$.

Our new arcs are: (4 16), (15 18), (17 10), (9 15), (2 6), (5 14), (6 5).

$D_2$

| 1 | 2 | 3 | 4 | 5 | 6 | 7 | 8 | 9 | 10 | 11 | 12 | 13 | 14 | 15 | 16 | 17 | 18 | 19 | 20 |
|---|---|---|---|---|---|---|---|---|----|----|----|----|----|----|----|----|----|----|----|
| 7 | 6 | 4 | 16 | 14 | 5 | 8 | 9 | 15 | 11 | 12 | 13 | 20 | 17 | 18 | 2 | 10 | 19 | 3 | 1 |
| 0 | -7 | -69 | -9 | 0 | -1 | -73 | -58 | 0 | -32 | -79 | -77 | 0 | 0 | -21 | 0 | 0 | -48 | 0 | -37 |

| 1 | 2 | 3 | 4 | 5 | 6 | 7 | 8 | 9 | 10 | 11 | 12 | 13 | 14 | 15 | 16 | 17 | 18 | 19 | 20 |
|---|---|---|---|---|---|---|---|---|----|----|----|----|----|----|----|----|----|----|----|

$D_2^{-1}$

| 20 | 16 | 19 | 3 | 6 | 2 | 1 | 7 | 8 | 17 | 10 | 11 | 12 | 5 | 9 | 4 | 14 | 15 | 18 | 13 |
|----|----|----|---|---|---|---|---|---|----|----|----|----|---|---|---|----|----|----|----|

$d(2, 15) - d(2, 6) = 11 - 18 = -7$.

$d(4, 17) - d(4, 16) = 9 - 18 = -9$.

$d(5, 14) - d(5, 14) = 0$.

$d(6, 14) - d(6, 5) = 18 - 19 = -1$

$d(9, 15) - d(9, 15) = 0$.

$d(15, 16) - d(15, 18) = 18 - 39 = -21$.

$d(17, 10) - d(17, 10) = 0$.

We choose 11 as an initial vertex for our trials.

TRIAL 1. MIN(M)(11, 1) = 1.

$$(11, 1) \rightarrow (11, 20)$$
$$(20, 2) \rightarrow (20, 16)$$
$$(16, 2) \rightarrow (16, 16)$$

(16, 2) is an arc of $D_2$.

$$(16, 6) \to (16, 2)$$
$$(2, 15) \to (2, 9)$$
$$(9, 15) \to (9, 9).$$

(9, 15) is an arc of $D_2$.

$$(9, 4) \to (9, 3).$$
$$(3, 11) \to (3, 10)$$
$$(10, 3) \to (10, 19)$$
$$(19, 3) \to (19, 19)$$

(19, 3) is an arc of $D_2$.

$$(19, 8) \to (19, 7)$$
$$(7, 1) \to (7, 20)$$

P = [11, 20, 16, 2 9 3 10 19, 7, 20].

$s_{11} = (11^{-79}\ 20^{-37}\ 16\ 2\ 9\ 3^{-69}\ 10^{-32}\ 19\ 7)$

d(16, 6) - d(16, 2) = 5 - 5 = 0.

d(2, 6) - d(2, 15) = 18 - 11 = 7.

d(9, 4) - d(9, 15) = 6 - 4 = 2.

d(19, 8) - d(19, 3) = 17 - 6 = 11.

d(7, 12) - d(7, 1) = 59 - 1 = 58

$s_{11} = (11^{-79} 20^{-37} 16^0 2^7 9^2 3^{-69} 10^{-32} 19^{11} 7^{58})$.

The value of $s_{11}$ is -139.

$s_{111} = (11\ 20)$, $s_{112} = (11\ 20\ 16)$, $s_{113} = (11\ 20\ 16\ 2)$, $s_{114} = (11\ 20\ 16\ 2\ 9)$, $s_{115} = (11\ 20\ 16\ 2\ 9\ 3)$, $s_{116} = (11\ 20\ 16\ 2\ 9\ 3\ 10)$, $s_{117} = (11\ 20\ 16\ 2\ 9\ 3\ 10\ 19)$.

d(20, 12) - d(20, 1) = 80 - 43 = 37.

d(16, 12) - d(16, 2) = 12 - 5 = 7.

d(2, 12) - d(2, 6) = 76 - 18 = 58.

d(9, 12) - d(9, 15) = 24 - 4 = 20.

d(3, 12) - d(3, 4) = 33 - 82 = -49.

d(10, 12) - d(10, 11) = 10 - 34 = -24.

d(19, 12) - d(19, 3) = 63 - 6 = 57.

$s_{111} = (11^{-79} 20^{37})$, $s_{112} = (11^{-79} 20^{-37} 16^{7})$, $s_{113} = (11^{-79} 20^{-37} 16^{0} 2^{58})$,

$s_{114} = (11^{-79} 20^{-37} 16^{0} 2^{7} 9^{20})$, $s_{115} = (11^{-79} 20^{-37} 16^{0} 2^{7} 9^{2} 3^{-49})$.

$s_{116} = (11^{-79} 20^{-37} 16^{0} 2^{7} 9^{2} 3^{-69} 10^{-24})$, $s_{117} = (11^{-79} 20^{-37} 16^{0} 2^{7} 9^{2} 3^{-69} 10^{-32} 19^{57})$.

The respective values of the above permutations are: -42, -109, -58, -89, -156, -200, -151.

$s_{12} = (20^{-37} 16^{0} 2^{7} 9^{2} 3^{-69} 10^{-32} 19^{11} 7^{0})$.

The value of $s_{12}$ is -118.

TRIAL 2. MIN(M)(11, 2) = 2.

$$(11, 2) \rightarrow (11, 16)$$

From this point on, the path is the same as that of $P$ in trial 1 up until we reach 7. Thus, let us continue on from 7 allowing $P'$ to be

$P' = [11, 16, 2, 9, 3, 10, 19, 7]$

$$(7, 1) \rightarrow (7, 20)$$
$$(20, 2) \rightarrow (20, 16).$$

Thus, $P = [11, 16, 2, 9, 3, 10, 19, 7, 20, 16]$

$s_{21} = (11\ 16\ 2\ 9\ 3\ 10\ 19\ 7\ 20)$

d(11, 2) - d(11, 12) = 14 - 90 = -76.

d(20, 12) - d(20, 1) = 80 - 43 = 37.

$s_{21} = (11^{-76}16^0 2^7 9^2 3^{-69} 10^{-32} 19^{11} 7^0 20^{37})$.

The value of $s_{21}$ is -120.

$s_{211} = (11^{-76}16)$, $s_{212} = (11^{-76}16^0 2)$, $s_{213} = (11^{-76}16^0 2^7 9)$,

$s_{214} = (11^{-76}16^0 2^7 9^2 3)$, $s_{215} = (11^{-76}16^0 2^7 9^2 3^{-69})$, $s_{216} = (11^{-76}16^0 2^7 9^2 3^{-69} 10^{-32} 19)$,

$s_{217} = (11^{-76}16^0 2^7 9^2 3^{-69} 10^{-32} 19^{11} 7)$

d(16, 12) - d(16, 2) = 12 - 5 = 7.

d(2, 12) - d(2, 6) = 76 - 18 = 58.

d(9, 12) - d(9, 15) = 24 - 4 = 20.

d(3, 12) - d(3, 4) = 33 - 82 = -49.

d(10, 12) - d(10, 11) = 10 - 34 = -24.

d(19, 12) - d(19, 3) = 80 - 6 = 74.

d(7, 12) - d(7, 8) = 59 - 74 = -15.

$s_{211} = (11^{-76}16^7)$, $s_{212} = (11^{-76}16^0 2^{58})$, $s_{213} = (11^{-76}16^0 2^7 9^{20})$,

$s_{214} = (11^{-76}16^0 2^7 9^2 3^{-49})$, $s_{215} = (11^{-76}16^0 2^7 9^2 3^{-69} 10^{-24})$,

$s_{216} = (11^{-76}16^0 2^7 9^2 3^{-69} 10^{-32} 19^{74})$, $s_{217} = (11^{-76}16^0 2^7 9^2 3^{-69} 10^{-32} 19^{11} 7^{-15})$.

The respective values of the above permutations are: -69, -18, -49,

-116, -172, -94, -172.

$s_{22} = (16^0 2^7 9^2 3^{-69} 10^{-32} 19^{11} 7^0 20)$

d(20, 4) - d(20, 1) = 51 - 43 = 8.

$s_{22} = (16^0 2^7 9^2 3^{-69} 10^{-32} 19^{11} 7^0 20^8)$.

The value of $s_{22}$ is -73.

TRIAL 3. MIN(M)(11, 3) = 4.

$$(11, 4) \to (11, 3)$$

$$(3, 11) \to (3, 10)$$

Using $s_{22}$, we can assume that $P'$, a subpath of $P$, is

$P' = [11, 3, 10, 19, 7, 20]$.

Continuing from 20, we obtain

$$(20, 2) \rightarrow (20, 16)$$

$$(16, 2) \rightarrow (16, 16)$$

$(16, 2)$ is an arc of $D_2$.

$$(16, 6) \rightarrow (16, 2)$$

$$(2, 15) \rightarrow (2, 9)$$

$$(9, 15) \rightarrow (9, 9)$$

$(9, 15)$ is an arc of $D_2$.

$$(9, 4) \rightarrow (9, 3).$$

$P = [11, 3, 10\ 19, 7, 20, 16, 2, 9, 3]$.

$s_{31} = (11\ 3\ 10\ 19\ 7\ 20\ 16\ 2\ 9)$.

$d(11, 4) - d(11, 12) = 18 - 90 = -72$.

$d(20, 2) - d(20, 1) = 6 - 43 = -37$.

$d(16, 6) - d(16, 2) = 5 - 5 = 0$.

$d(2, 15) - d(2, 6) = 11 - 18 = -7$.

$d(9, 12) - d(9, 15) = 24 - 4 = 20$.

$s_{31} = ((11^{-72}3^{-69}10^{-32}19^{11}7^{0}20^{-37}16^{0}2^{-7}9^{20}))$.

The value of $s_{31}$ is -186.

$d(3, 12) - d(3, 4) = 33 - 82 = -59$.

$d(10, 12) - d(10, 11) = 10 - 34 = -24$.

$d(19, 12) - d(19, 3) = 63 - 6 = 57$.

$d(7, 12) - d(7, 8) = 59 - 74 = 15$.

$d(20, 12) - d(20, 1) = 80 - 43 = 37$.

d(16, 12) - d(16, 2) = 12 - 5 = 7.

d(2, 12) - d(2, 6) = 76 - 18 = 58.

$s_{311} = (11^{-72}3^{-59})$, $s_{312} = (11^{-72}3^{-69}10^{-24})$, $s_{313} = (11^{-72}3^{-69}10^{-32}19^{57})$,

$s_{314} = (11^{-72}3^{-59}10^{-32}19^{11}7^{15})$, $s_{315} = (11^{-72}3^{-69}10^{-32}19^{11}7^{0}20^{37})$,

$s_{316} = (11^{-72}3^{-69}10^{-32}19^{11}7^{0}20^{-37}16^{7})$, $s_{317} = (11^{-72}3^{-69}10^{-32}19^{11}20^{-37}16^{0}2^{58})$.

The respective values of the above permutations are: -131, -173, -116,

-147, -141, -192, -141.

d(9, 4) - d(9, 15) = 6 - 4 = 2.

$s_{32} = (3^{-69}10^{-32}19^{11}7^{0}20^{-37}16^{0}2^{-7}9^{2})$.

The value of $s_{32}$ is -132. It follows that our best possibility for $s_3$ is

$s_{116} = (11^{-79}20^{-37}16^{0}2^{7}9^{2}3^{-69}10^{-24})$.

(11, 20) → (11, 1), → (20, 16) → (20, 2), (16, 2) → (16, 6),

(2, 9) → (2, 15), (9, 3) → (9, 4), (3, 10) → (3, 11), (10, 11) → (10, 12).

Our new arcs in $D_3$ are: (11, 1), (20, 2), (16, 6), (2, 15), (9, 4), (3, 11) and (10, 12).

```
        * *              * * *             *           *
     1  2  3  4  5  6  7  8  9 10 11 12 13 14 15 16 17 18 19 20
D₃
     7 15 11 16 14  5  8  9  4 12  1 13 20 17 18  6 10 19  3  2
```

d(2, 15) - d(2, 15) = 0.

d(3, 11) - d(3, 11) = 0.

d(9, 15) - d(9, 4) = 4 - 6 = -2.

d(10, 3) - d(10, 12) = 2 - 10 = -8.

d(11, 1) - d(11, 1) = 0.

d(16, 2) - d(16, 6) = 5 - 5 = 0.

d(20, 2) - d(20, 2) = 0.

| 0 | 0 | 0 | -9 | 0 | 0 | -73 | -58 | -2 | -8 | 0 | -77 | 0 | 0 | -21 | 0 | 0 | -48 | 0 | 0 |
|---|---|---|---|---|---|---|---|---|---|---|---|---|---|---|---|---|---|---|---|
| 1 | 2 | 3 | 4 | 5 | 6 | 7 | 8 | 9 | 10 | 11 | 12 | 13 | 14 | 15 | 16 | 17 | 18 | 19 | 20 |
| 11 | 20 | 19 | 9 | 6 | 16 | 1 | 7 | 8 | 17 | 3 | 10 | 12 | 5 | 2 | 4 | 14 | 15 | 18 | 13 |

The smallest value of *DIFF* occurs out of 12.

TRIAL 1. MIN(M)(12, 1) = 2.

$$(12, 2) \rightarrow (12, 20)$$

$$(20, 2) \rightarrow (20, 20).$$

(20, 2) is an arc of $D_3$.

$$(20, 15) \rightarrow (20, 2)$$

$$(2, 15) \rightarrow (2, 2).$$

(2, 15) is an arc of $D_3$.

$$(2, 6) \rightarrow (2, 16)$$

$$(16, 2) \rightarrow (16, 20)$$

P = [12, 20, 2, 16, 20].

$s_{11}$ = (12 20 2 16)

d(20, 2) - d(20, 15) = 0.

d(2, 6) - d(2, 15) = 18 - 11 = 7.

d(16, 2) - d(16, 13) = 59 - 5 = 54.

$s_{11} = (12^{-77} \, 20^0 \, 2^7 \, 16^{54})$.

The value of $s_{11}$ is -16.

$s_{12}$ = (20 2 16).

d(16, 2) - d(16, 2) = 0.

$s_{12} = (20^0 \, 2^7 \, 16^0)$.

The value of $s_{12}$ is non-negative.

TRIAL 2. MIN(M)(12, 2) = 5.

d(12, 5) - d(12, 13) = 6 - 78 = -72.

$$(12, 5) \rightarrow (12, 6)$$

$$(6, 14) \rightarrow (6, 5)$$

$$(5, 14) \rightarrow (5, 5)$$

(5, 14) is an arc of $D_3$.

$$(5, 17) \rightarrow (5, 14)$$

$$(14, 17) \rightarrow (14, 14)$$

(14, 17) is an arc of $D_3$.

$$(14, 13) \rightarrow (14, 12)$$

P = [12, 6, 5, 14, 12]

$s_{21} = (12 \, 6 \, 5 \, 14)$

d(6, 14) - d(6, 5) = 18 - 19 = -1.

d(5, 17) - d(5, 14) = 10 - 4 = 6.

d(14, 13) - d(14, 17) = 19 - 6 = 13.

$s_{21} = (12^{-72} \, 6^{-1} \, 5^6 \, 14^{13})$.

The value of $s_{21}$ is -54.

TRIAL 3. MIN(M)(12, 3) = 20.

d(12, 20) - d(12, 13) = 7 - 78 = -71.

$$(12, 20) \rightarrow (12, 13)$$

$$(13, 20) \rightarrow (13, 13)$$

(13, 20) is an arc of $D_3$.

$$(13, 19) \to (13, 18)$$

$$(18, 15) \to (18, 2)$$

$$(2, 15) \to (2, 2)$$

(2, 15) is an arc of $D_3$.

$$(2, 6) \to (2, 16)$$

$$(16, 2) \to (16, 20)$$

$$(20, 2) \to (20, 20)$$

(20, 2) is an arc of $D_3$.

$$(20, 15) \to (20, 2)$$

$P = [12, 13, 18, 2, 16, 20, 2]$.

$s_{31} = (12\ 13\ 18\ 2\ 16\ 20)$

$(12, 13) \to (12, 20), (13, 18) \to (13, 19), (18, 2) \to (18, 15),$

$(2, 16) \to (2, 6), (16, 20) \to (16, 2), (20, 12) \to (20, 13)$.

d(12, 20) - d(12, 13) = 7 - 78 = -71.

d(13, 19) - d(13, 20) = 7 - 4 = 3.

d(18, 15) - d(18, 19) = 3 - 51 = -48.

d(2, 6) - d(2, 15) = 18 - 11 = 7.

d(16, 2) - d(16, 6) = 5 - 5 = 0.

d(20, 13) - d(20, 2) = 65 - 6 = 59.

$(13, 12) \to (13, 13)$. Thus, $s_{311}$ need not be considered.

$s_{312} = (12^{-71} 13^3 18^{-9})$, $s_{313} = (12^{-71} 13^3 18^{-48} 2^{21})$, $s_{314} = (12^{-71} 13^3 18^{-48} 2^7 16^{54})$.

The respective values of the above permutations are: -77, -95, -55.

d(20, 15) – d(20, 2) = 6 – 6 = 0.

$s_{32} = (12^{-71}13^{3}18^{-9})(2^{7}16^{0}20^{0})$.

The value of the second cycle of $s_{32}$ is non-negative, while the first cycle is precisely $s_{312}$.

The smallest-valued permutation is $s_{313}$. Thus, $s_{3} = (12\ 13\ 18\ 2)$.

$(12, 13) \rightarrow (12, 20), (13, 18) \rightarrow (13, 19), (18, 2) \rightarrow (18, 15)$,

$(2, 12) \rightarrow (2, 13)$. Thus, the new arcs of $D_4$ are: $(12, 20), (13, 19)$,

$(18, 15)$ and $(2, 13)$.

$D_4$

|   |   |   |   |   |   |   |   |   |   |   | * |   |   | * | * |   |   | * |   |   |
|---|---|---|---|---|---|---|---|---|---|---|---|---|---|---|---|---|---|---|---|---|
|   | 1 | 2 | 3 | 4 | 5 | 6 | 7 | 8 | 9 | 10 | 11 | 12 | 13 | 14 | 15 | 16 | 17 | 18 | 19 | 20 |
|   | 7 | 13 | 11 | 16 | 14 | 5 | 8 | 9 | 4 | 12 | 1 | 20 | 19 | 17 | 18 | 6 | 10 | 15 | 3 | 2 |

In order to create $D_4^{-1}$, we need only invert the four new arcs and use all of the remaining arcs from $D_3^{-1}$.

Before going on, we obtain $DIFF(a)$ for $a = 2, 12, 13$ and $18$.

$d(2, 15) - d(2, 13) = 11 - 32 = -21$.

$d(12, 2) - d(12, 20) = 1 - 7 = -6$.

$d(13, 20) - d(13, 19) = 4 - 7 = -3$

$d(18, 15) - d(18, 15) = 0$.

| 0 | –21 | 0 | -9 | 0 | 0 | –73 | -58 | -2 | -8 | 0 | -6 | -3 | 0 | -21 | 0 | 0 | 0 | 0 | 0 |
|---|---|---|---|---|---|---|---|---|---|---|---|---|---|---|---|---|---|---|---|
| 1 | 2 | 3 | 4 | 5 | 6 | 7 | 8 | 9 | 10 | 11 | 12 | 13 | 14 | 15 | 16 | 17 | 18 | 19 | 20 |
| 11 | 20 | 19 | 9 | 6 | 16 | 1 | 7 | 8 | 17 | 3 | 10 | 2 | 5 | 18 | 4 | 14 | 15 | 13 | 12 |

The smallest $DIFF$ value occurs when $a = 7$.

TRIAL 1.

$(7, 1) \rightarrow (7, 11)$

$$(11, 1) \rightarrow (11, 11)$$

$(11, 1)$ is an arc of $D_4$.

$$(11, 2) \rightarrow (11, 20)$$

$$(20, 2) \rightarrow (20, 20)$$

$(20, 2)$ is an arc of $D_4$.

$$(20, 15) \rightarrow (20, 18)$$

$$(18, 15) \rightarrow (18, 18)$$

$(18, 15)$ is an arc of $D_4$.

$$(18, 14) \rightarrow (18, 5)$$

$$(5, 14) \rightarrow (5, 5)$$

$(5, 14)$ is an arc of $D_4$.

$$(5, 17) \rightarrow (5, 14)$$

$$(14, 17) \rightarrow (14, 14)$$

$(14, 17)$ is an arc of $D_4$.

$$(14, 13) \rightarrow (14, 2)$$

$$(2, 15) \rightarrow (2, 18)$$

$P = [7, 11, 20, 18, 5, 14, 2, 18]$.

$s_{11} = (7\ 11\ \ 20\ 18\ 5\ 14\ 2)$

$d(7, 1) - d(7, 8) = 1 - 74 = -73$

$d(11, 2) - d(11, 1) = 14 - 11 = 3$

$d(20, 15) - d(20, 2) = 6 - 6 = 0.$

$d(18, 14) - d(18, 15) = 9 - 3 = 6.$

$d(5, 17) - d(5, 14) = 10 - 4 = 6.$

$d(14, 13) - d(14, 17) = 19 - 6 = 13.$

$d(2, 8) - d(2, 13) = 20 - 32 = -12.$

$s_{11} = (7^{-73} 11^3 20^0 18^6 5^6 14^{13} 2^{-12})$.

The value of $s_{11}$ is –57.

d(11, 8) – d(11, 1) = 51 – 11 = 40.

d(20, 8) – d(20, 2) = 56 – 6 = 50.

d(18, 8) – d(18, 15) = 51 – 3 = 48.

d(5, 8) – d(5, 14) = 64 – 4 = 60.

d(14, 8) – d(14, 17) = 52 – 6 = 46.

$s_{114} = (7^{-73} 11^3 20^0 18^6 5^{60})$, $s_{115} = (7^{-73} 11^3 20^0 18^6 5^6 14^{46})$.

The respective values of the above permutations are: -33, -20, -21, -14, -12.

$s_{12} = (7\ 11\ 20)(18\ 5\ 14) = (7^{-73} 11^3 20^{50})(18^6 5^6 14^{85})$.

d(14, 15) – d(14, 17) = 91 – 6 = 85.

The first cycle is $s_{112}$; the second cycle is non-negative.

TRIAL 2. MIN(M)(7, 2) = 5.

d(7, 5) – d(7,8) = 4 – 74 = -70

$$(7, 5) \rightarrow (7, 6)$$

$$(6, 14) \rightarrow (6, 5)$$

In Trial 1, 5 was followed by 14, 2 and 18. We continue with 18.

$$(18, 15) \rightarrow (18, 18)$$

(18, 15) is an arc of $D_4$.

$$(18, 14) \rightarrow (18, 5)$$

P = [7, 6, 5, 14, 2,18, 5].

$s_{21} = (7\ 6\ 5\ 14\ 2\ 18)$.

d(6, 14) – d(6, 5) = 18 – 19 = -1.

d(5, 17) – d(5, 14) = 6.

$d(14, 13) - d(14, 17) = 13$.

$d(2, 15) - d(2, 13) = 11 - 32 = -21$.

$d(18, 8) - d(18, 15) = 51 - 3 = 48$.

$s_{21} = (7^{-70} 6^{-1} 5^6 14^{13} 2^{-21} 18^{48})$

The value of $s_{21}$ is -25.

$d(6, 8) - d(6, 5) = 94 - 19 = 75$.

$d(5, 8) - d(5, 14) = 64 - 4 = 60$.

$d(14, 8) - d(14, 17) = 52 - 6 = 46$.

$d(2, 8) - d(2, 13) = 20 - 32 = -12$.

$s_{211} = (7^{-70} 6^{75})$, $s_{212} = (7^{-70} 6^0 5^{60})$, $s_{213} = (7^{-70} 6^0 5^6 14^{46})$,

$s_{214} = (7^{-70} 6^0 5^6 14^{13} 2^{-12})$.

The respective values of the above permutations are: 5, -10, -18, -63.

$d(2, 15) - d(2, 13) = 11 - 32 = -21$.

$d(18, 14) - d(18, 15) = 9 - 3 = 6$.

$s_{22} = (7^{-70} 6^{76})(5^6 14^{13} 2^{-21} 18^6)$.

Each cycle has a non-negative value.

TRIAL 3. MIN(M)(7, 3) = 9.

$d(7, 9) - d(7,8) = 28 - 74 = -46$.

$$(7, 9) \rightarrow (7, 4)$$

$$(4, 17) \rightarrow (4, 14)$$

In Trial 1, we obtained 14, 2, 18, 5.

We continue with 5.

$$(5, 14) \rightarrow (5, 5)$$

(5, 14) is an arc of $D_4$.

$$(5, 17) \rightarrow (5, 14)$$

$P = [7, 4, 14, 2, 18, 5, 14]$.

$s_{31} = (7\ 4\ 14\ 2\ 18\ 5)$.

$d(4, 17) - d(4, 16) = 9 - 18 = -9$.

$d(14, 13) - d(14, 17) = 19 - 6 = 13$.

$d(2, 15) - d(2, 13) = 11 - 32 = -21$.

$d(18, 14) - d(18, 15) = 9 - 3 = 6$.

$d(5, 8) - d(5, 14) = 64 - 4 = 60$.

$s_{31} = (7^{-46} 4^{-9} 14^{13} 2^{-21} 18^{6} 5^{60})$.

The value of $s_{31}$ is 3.

$d(4, 8) - d(4, 16) = 88 - 18 = 70$.

$d(14, 8) - d(14, 17) = 52 - 6 = 46$.

$d(2, 8) - d(2, 13) = 20 - 32 = -12$.

$d(18, 8) - d(18, 15) = 51 - 3 = 48$.

$s_{311} = (7^{-46} 4^{70})$, $s_{312} = (7^{-46} 4^{-9} 14^{46})$, $s_{313} = (7^{-46} 4^{-9} 14^{13} 2^{-12})$,

$s_{314} = (7^{-46} 4^{-9} 14^{13} 2^{-21} 18^{48})$.

The respective values of the above permutations are: 24, 1, -54, -15.

$d(5, 17) - d(5, 14) = 10 - 4 = 6$.

$s_{32} = (7^{-46} 4^{70})(14^{13} 2^{-21} 18^{6} 5^{6})$

Both cycles are non-negative.

The smallest value of all permutations is $s_{214} = (7^{-70} 6^{0} 5^{6} 14^{13} 2^{-12})$. We obtain the following new arcs for

$D_5$: (7, 5), (6, 14), (5, 17), (14, 13), (2, 8).

        *        * * *                      *

$D_5$

| 1 | 2 | 3 | 4 | 5 | 6 | 7 | 8 | 9 | 10 | 11 | 12 | 13 | 14 | 15 | 16 | 17 | 18 | 19 | 20 |
|---|---|---|---|---|---|---|---|---|----|----|----|----|----|----|----|----|----|----|----|
| 7 | 8 | 11 | 16 | 17 | 14 | 5 | 9 | 4 | 12 | 1 | 20 | 19 | 13 | 18 | 6 | 10 | 15 | 3 | 2 |

d(2, 15) – d(2, 8) = 11 – 20 = -9.

d(5, 14) – d(5, 17) = 4 – 10 = -6.

d(6, 14) – d(6,14) = 0.

d(7, 1) - d(7, 5) = 1 – 4 = -3.

d(14, 17) – d(14, 13) = 6 – 19 = -13.

| 0 | –9 | 0 | -9 | -6 | 0 | -3 | –58 | -2 | -8 | 0 | -6 | -3 | –13 | -21 | 0 | 0 | 0 | 0 | 0 |
|---|----|---|----|----|---|----|-----|----|----|---|----|----|-----|-----|---|---|---|---|---|
| 1 | 2 | 3 | 4 | 5 | 6 | 7 | 8 | 9 | 10 | 11 | 12 | 13 | 14 | 15 | 16 | 17 | 18 | 19 | 20 |

$D_5^{-1}$

| 11 | 20 | 19 | 9 | 7 | 16 | 1 | 2 | 8 | 17 | 3 | 10 | 14 | 6 | 18 | 4 | 5 | 15 | 13 | 12 |
|----|----|----|---|---|----|---|---|---|----|---|----|----|---|----|---|---|----|----|----|

We start with $a = 8$.

TRIAL 1.

d(8, 1) – d(8, 9) = 1 – 59 = -58.

(8, 1) → (8, 11)

(11, 1) → (11, 11)

(11, 1) is an arc of $D_5$.

(11, 2) → (11, 20)

(20, 2) → (20, 20)

(20, 2) is an arc of $D_5$.

(20, 15) → (20, 18)

$$(18, 15) \rightarrow (18, 18)$$

(18, 15) is an arc of $D_5$.

$$(18, 14) \rightarrow (18, 6)$$

$$(6, 14) \rightarrow (6, 6)$$

(6, 14) is an arc of $D_5$

$$(6, 2) \rightarrow (6, 20)$$

$P = [8, 11, 20, 18, 6, 20]$.

$s_{11} = (8\ 11\ 20\ 18\ 6)$

$d(8, 1) - d(8, 9) = -58$.

$d(11, 2) - d(11, 1) = 14 - 11 = 3$.

$d(20, 15) - d(20, 2) = 6 - 6 = 0$.

$d(18, 14) - d(18, 15) = 9 - 3 = 6$.

$d(6, 9) - d(6, 14) = 72 - 18 = 54$.

$s_{11} = (8^{-58} 11^3 20^0 18^6 6^{54})$.

$d(11, 9) - d(11, 1) = 32 - 11 = 21$.

$d(20, 9) - d(20, 2) = 35 - 6 = 29$.

$d(18, 9) - d(18, 15) = 83 - 3 = 80$.

$s_{111} = (8^{-58} 11^{21})$, $s_{112} = (8^{-58} 11^3 20^{29})$, $s_{113} = (8^{-58} 11^3 20^0 18^{80})$.

The respective values of the above permutations are: -37, -27, 25.

$d(6, 2) - d(6, 14) = 19 - 18 = 1$.

$s_{12} = (8^{-58} 11^{21})(20^0 18^6 6^1)$.

The first cycle is $s_{111}$. The second cycle has a non-negative value.

TRIAL 2. MIN(M)(8, 2) = 17.

$d(8, 17) - d(8, 9) = 24 - 58 = -34$.

$$(8, 17) \rightarrow (8, 5)$$

$$(5, 14) \rightarrow (5, 13)$$

$$(13, 20) \rightarrow (13, 2)$$

$$(2, 15) \rightarrow (2, 18)$$

$$(18, 15) \rightarrow (18, 18)$$

$(18, 15)$ is an arc of $D_5$.

$$(18, 14) \rightarrow (18, 6)$$

$$(6, 14) \rightarrow (6, 6)$$

$(6, 14)$ is an arc of $D_5$.

$$(6, 2) \rightarrow (6, 20)$$

$$(20, 2) \rightarrow (20, 20)$$

$(20, 2)$ is an arc of $D_5$.

$$(20, 15) \rightarrow (20, 18)$$

$P = [8, 5, 13, 2, 18, 6, 20, 18]$.

$s_{21} = (8\ 5\ 13\ 2\ 18\ 6\ 20)$.

$d(8, 17) - d(8, 9) = -34$.

$d(5, 14) - d(5, 17) = 4 - 10 = -6$.

$d(13, 20) - d(13, 19) = 4 - 7 = -3$.

$d(2, 15) - d(2, 8) = 11 - 20 = -9$.

$d(18, 14) - d(18, 15) = 9 - 3 = 6$.

$d(6, 2) - d(6, 14) = 19 - 18 = 1$.

$d(20, 15) - d(20, 2) = 6 - 6 = 0$.

$s_{21} = (8^{-34} 5^{-6} 13^{-3} 2^{-9} 18^6 6^1 20^{29})$.

The value of $s_{21}$ is $-16$.

$d(5, 9) - d(5, 17) = 77 - 10 = 67$.

d(13, 9) – d(13, 19) = 15 – 7 = 8.

d(2, 9) – d(2, 8) = 58 – 20 = 38.

d(18, 9) – d(18, 15) = 83 – 3 = 80.

d(6, 9) – d(6, 14) = 72 – 18 = 54.

$s_{211} = (8^{-34}5^{67})$, $s_{212} = (8^{-34}5^{-6}13^{8})$, $s_{213} = (8^{-34}5^{-6}13^{-3}2^{38})$,

$s_{214} = (8^{-34}5^{-6}13^{-3}2^{-9}18^{80})$, $s_{215} = (8^{-34}5^{-6}13^{-3}2^{-9}18^{6}6^{54})$.

The respective values of the above permutations are: 33, -32, -5, 28, 8.

$s_{22} = (8\ 5\ 13\ 2)(18\ 6\ 20)$.

The first cycle is $s_{213}$.

d(20, 15) – d(20, 2) = 6 – 6 = 0. Thus, the second cycle is $(18^{6}6^{1}20^{0})$. Thus, it has a non-negative value.

TRIAL 3. MIN(M)(8, 3) = 11.

$$d(8, 11) – d(8,9) = 25 – 59 = -34.$$

$$(8, 11) \to (8, 3)$$

$$(3, 11) \to (3, 3)$$

(3, 11) is an arc of $D_5$.

$$(3, 10) \to (3, 17)$$

$$(17, 10) \to (17, 17)$$

(17, 10) is an arc of $D_5$.

$$(17, 5) \to (17, 7)$$

$$(7, 1) \to (7, 11)$$

$$(11, 1) \to (11, 11)$$

(11, 1) is an arc of $D_5$.

$$(11, 2) \to (11, 20)$$

$$(20, 2) \rightarrow (20, 20)$$

$(20, 2)$ is an arc of $D_5$.

$$(20, 15) \rightarrow (20, 18)$$
$$(18, 15) \rightarrow (18, 18)$$

$(18, 15)$ is an arc of $D_5$.

$$(18, 14) \rightarrow (18, 6)$$
$$(6, 14) \rightarrow (6, 6)$$

$(6, 14)$ is an arc of $D_5$.

$$(6, 2) \rightarrow (6, 20)$$
$$(20, 2) \rightarrow (20, 20)$$

$(20, 2)$ is an arc of $D_5$.

$$(20, 15) \rightarrow (20, 18)$$

$P = [8, 3, 17, 7, 11, 20, 18, 6, 20]$.

$d(3, 10) - d(3, 11) = 17 - 13 = 4$.

$d(17, 5) - d(17, 10) = 8 - 1 = 7$.

$d(7, 1) - d(7, 5) = 1 - 4 = -3$.

$d(11, 2) - d(11, 1) = 14 - 11 = 3$.

$d(20, 15) - d(20, 2) = 6 - 6 = 0$.

$d(18, 14) - d(18, 15) = 9 - 3 = 6$.

$d(6, 9) - d(6, 14) = 72 - 18 = 54$.

$s_{31} = (8^{-34} 3^4 17^7 7^{-3} 11^3 20^0 18^6 6^{54})$.

$d(3, 9) - d(3, 11) = 97 - 13 = 84$.

$d(17, 9) - d(17, 10) = 28 - 1 = 27$.

$d(7, 9) - d(7, 5) = 28 - 4 = 24$.

$d(11, 9) - d(11, 2) = 32 - 14 = 18$.

d(20, 9) − d(20, 2) = 35 − 6 = 29.

d(18, 9) − d(18, 15) = 83 − 3 = 80.

$s_{311} = (8^{-34}3^{84})$, $s_{312} = (8^{-34}3^{4}17^{27})$, $s_{313} = (8^{-34}3^{4}17^{7}7^{24})$,

$s_{314} = (8^{-34}3^{4}17^{7}7^{-3}11^{18})$, $s_{315} = (8^{-34}3^{4}17^{7}7^{-3}11^{3}20^{29})$,

$s_{316} = (8^{-34}3^{4}17^{7}7^{-3}11^{3}20^{0}18^{80})$.

$s_{32} = (8^{-34}3^{4}17^{7}7^{-3}11^{18})(20^{0}18^{6}6^{1})$.

The first cycle is $s_{314}$. The second is non-negative.

The smallest negative cycle is $s_{111} = (8^{-58}\,11^{21})$. Thus, $s_4 = (8\ 11)$.

|  | 1 | 2 | 3 | 4 | 5 | 6 | 7 | 8 | 9 | 10 | 11 | 12 | 13 | 14 | 15 | 16 | 17 | 18 | 19 | 20 |
|---|---|---|---|---|---|---|---|---|---|---|---|---|---|---|---|---|---|---|---|---|
| $D_6$ | 7 | 8 | 11 | 16 | 17 | 14 | 5 | 1 | 4 | 12 | 9 | 20 | 19 | 13 | 18 | 6 | 10 | 15 | 3 | 2 |
|  | −7 | −9 | 0 | −9 | −6 | 0 | −3 | 0 | −2 | −8 | −21 | −6 | −3 | −13 | −21 | 0 | 0 | 0 | 0 | 0 |
|  | 1 | 2 | 3 | 4 | 5 | 6 | 7 | 8 | 9 | 10 | 11 | 12 | 13 | 14 | 15 | 16 | 17 | 18 | 19 | 20 |
| $D_6^{-1}$ | 8 | 20 | 19 | 9 | 7 | 16 | 1 | 2 | 11 | 17 | 3 | 10 | 14 | 6 | 18 | 4 | 5 | 15 | 13 | 12 |

TRIAL 1. MIN(M)(15, 1) = 16.

$$d(15,\ 16) - d(15,\ 18) = 18 - 39 = -21.$$

$$(15, 16) \rightarrow (15, 4)$$

$$(4, 17) \rightarrow (4, 5)$$

$$(5, 14) \rightarrow (5, 6)$$

$$(6, 14) \rightarrow (6, 6)$$

(6, 14) is an arc of $D_6$.

$$(6, 2) \rightarrow (6, 20)$$

$$(20, 2) \rightarrow (20, 20)$$

(20, 2) is an arc of $D_6$.

$$(20, 15) \rightarrow (20, 18)$$

$$(18, 15) \rightarrow (18, 18)$$

(18, 15) is an arc of $D_6$.

$$(18, 14) \rightarrow (18, 6)$$

$P = [15, 4, 5, 6, 20, 18, 6]$.

$d(15, 16) - d(15, 18) = 18 - 39 = -21$.

$d(4, 17) - d(4, 16) = 9 - 18 = -9$

$d(5, 14) - d(5, 17) = 4 - 9 = -5$.

$d(6, 2) - d(6, 14) = 19 - 18 = 1$.

$d(20, 15) - d(20, 2) = 6 - 6 = 0$.

$s_{11} = (15\ 4\ 5\ 6\ 20\ 18)$.

We can't use $s_{11}$ since $D_6(18) = 18$. Thus, $D_6\ s_{11}$ wouldn't be a derangement.

$d(4, 18) - d(4, 16) = 67 - 18 = 49$.

$d(5, 18) - d(5, 17) = 19 - 10 = 9$.

$d(6, 18) - d(6, 14) = 70 - 18 = 52$.

$d(20, 18) - d(20, 2) = 85 - 6 = 79$.

$s_{111} = (15^{-21}4^{49})$, $s_{112} = (15^{-21}4^{-9}5^9)$, $s_{113} = (15^{-21}4^{-9}5^{-5}6^{52})$,

$s_{114} = (15^{-21}4^{-9}5^{-5}6^{1}20^{79})$ .

The respective values of the above permutations are: 28, -21, 17, 45.

$s_{22} = (15^{-21}4^{-9}5^9)(6^{1}20^{0}18^{6})$.

The first cycle is $s_{112}$. The second cycle is non-negative.

TRIAL 2. MIN(M)(15, 2) = 18.

$$(15, 18) - (15, 18) = 0.$$

Thus, we can proceed no further in trials 2 and 3.

The smallest negative permutation obtained is $s_{112} = (15^{-21}4^{-9}5^9)$.

Thus, $s_6 = (15^{-21}4^{-9}5^9)$. The new arcs of $D_7$ are: (15, 16), (4, 17), (5, 18).

|   |   |   |   |   |   |   |   |   |   |   |   |   |   | * |   |   | * |   |   |
|---|---|---|---|---|---|---|---|---|---|---|---|---|---|---|---|---|---|---|---|
| 1 | 2 | 3 | 4 | 5 | 6 | 7 | 8 | 9 | 10 | 11 | 12 | 13 | 14 | 15 | 16 | 17 | 18 | 19 | 20 |

$D_7$

| 7 | 20 | 11 | 17 | 18 | 14 | 5 | 1 | 4 | 12 | 9 | 20 | 19 | 13 | 16 | 6 | 10 | 15 | 3 | 2 |

$d(15, 16) - d(15, 16) = 0.$

$d(4, 17) - d(4, 17) = 0.$

$d(5, 14) - d(5, 18) = 4 - 19 = -15.$

| -7 | -9 | 0 | 0 | -15 | 0 | -3 | 0 | -2 | -8 | -21 | -6 | -3 | -13 | 0 | 0 | 0 | 0 | 0 | 0 |
|---|---|---|---|---|---|---|---|---|---|---|---|---|---|---|---|---|---|---|---|
| 1 | 2 | 3 | 4 | 5 | 6 | 7 | 8 | 9 | 10 | 11 | 12 | 13 | 14 | 15 | 16 | 17 | 18 | 19 | 20 |

$D_7^{-1}$

| 8 | 20 | 19 | 9 | 7 | 16 | 1 | 2 | 11 | 17 | 3 | 10 | 14 | 6 | 18 | 15 | 4 | 5 | 13 | 12 |

We let $a = 11$.

TRIAL 1. MIN(M)(11, 1) = 1.

$$(11, 1) \rightarrow (11, 8)$$

$$(8, 1) \rightarrow (8, 8)$$

(8, 1) is an arc of $D_7$.

$$(8, 17) \rightarrow (8, 10)$$

$$(10, 3) \rightarrow (10, 19)$$

$$(19, 3) \rightarrow (19, 19)$$

$(19, 3)$ is an arc of $D_7$.

$$(19, 8) \rightarrow (19, 2)$$

$$(2, 15) \rightarrow (2, 18)$$

$$(18, 15) \rightarrow (18, 18)$$

$(18, 15)$ is an arc of $D_7$.

$$(18, 14) \rightarrow (18, 6)$$

$$(6, 14) \rightarrow (6, 6)$$

$(6, 14)$ is an arc of $D_7$.

$$(6, 2) \rightarrow (6, 20)$$

$$(20, 2) \rightarrow (20, 20)$$

$(20, 2)$ is an arc of $D_7$.

$$(20, 15) \rightarrow (20, 18)$$

$P = [11, 8, 10, 19, 2, 18. 6, 20, 18]$

$d(8, 17) - d(8, 1) = 24 - 1 = 23$.

The path $[11, 8, 10]$ has a positive value $-21 + 23 = 2$.

We need proceed no further in Trial 1. The only possible permutation is $(11\ 8)$.

$d(11, 1) - d(11, 9) = -21$.

$d(8, 9) - d(8, 1) = 59 - 1 = 58$.

Thus, the value of $(11\ 8)$ is non-negative.

TRIAL 2. MIN(M )(11, 2) = 2.

$$d(11, 2) - d(11, 9) = 14 - 32 = -18.$$

$$(11, 2) \rightarrow (11, 20)$$

$$(20, 2) \rightarrow (20, 20)$$

$(20, 2)$ is an arc of $D_7$.

$$(20, 15) \rightarrow (20, 18)$$

$d(20, 15) - d(20, 2) = 6 - 6 = 0$.

$$(18, 15) \rightarrow (18, 18)$$

$(18, 15)$ is an arc of $D_7$.

$$(18, 14) \rightarrow (18, 6)$$

$d(18, 14) - d(18, 15) = 9 - 3 = 6$.

$$(6, 14) \rightarrow (6, 6)$$

$(6, 14)$ is an arc of $D_7$.

$$(6, 2) \rightarrow (6, 20)$$

$d(6, 2) - d(6, 14) = 19 - 18 = 1$.

$P = [11, 20, 18, 6, 20]$

$d(6, 9) - d(6, 14) = 72 - 18 = 54$.

$s_{21} = (11^{-18} 20^0 18^6 6^{54})$.

The value of $s_{21}$ is non-negative.

$d(20, 9) - d(20, 2) = 35 - 6 = 29$.

$d(18, 9) - d(18, 15) = 83 - 3 = 80$.

$s_{211} = (11^{-21} 20^{29})$, $s_{212} = (11^{-21} 20^0 18^{80})$.

The value of both of the above permutations are non-negative.

$s_{22} = (20^0 18^6 6^1)$.

The value of $s_{22}$ is non-negative.

TRIAL 3. MIN(M)(11, 3) = 4.

$d(11, 4) - d(11, 9) = 18 - 32 = -14$.

$$(11, 4) \to (11, 9)$$

$$(9, 15) \to (9, 18)$$

$d(9, 15) - d(9, 4) = 4 - 6 = -2$.

$$(18, 15) \to (18, 18)$$

$(18, 15)$ is an arc of $D_7$.

$$(18, 14) \to (18, 6)$$

$d(18, 14) - d(18, 15) = 9 - 3 = 6$.

$$(6, 14) \to (6, 6)$$

$(6, 14)$ is an arc of $D_7$.

$$(6, 2) \to (6, 20)$$

$d(6, 2) - d(6, 14) = 19 - 18 = 1$.

$$(20, 2) \to (20, 20)$$

$(20, 2)$ is an arc of $D_7$.

$$(20, 15) \to (20, 18)$$

$d(20, 15) - d(20, 2) = 6 - 6 = 0$.

$P = [11, 9, 18, 6, 20, 18]$.

$d(20, 9) - d(20, 2) = 35 - 6 = 29$.

$s_{31} = (11^{-14} 9^{-2} 18^6 6^1 20^{29})$.

The value of $s_{31}$ is non-negative.

$d(9, 9)$ is not permissible.

$d(18, 9) - d(18, 15) = 83 - 3 = 80$.

$d(6, 9) - d(6, 14) = 72 - 18 = 54$.

$s_{312} = (11^{-14} 9^{-2} 18^{80})$, $s_{313} = (11^{-14} 9^{-2} 18^6 6^{54})$.

Each of the above permutations has a non-negative value.

$s_{32} = (18^6 6^1 20^0)$.

The value of $s_{32}$ is non-negative.

Since we were unable to obtain a negatively-valued permutation in all three trials, we can go no further using Phase 1 of our algorithm.

Before going on to PHASE 2, we note that we can still use MIN(M) to obtain smallest values. For example, if $(a\ i)$ is the $i-th$ smallest value in row $a$ of $M$, then $(a\ D_j^{-1}(i))$ is the $i-th$ smallest value in row $a$ of $D_j^{-1}(i)$. In order to simplify this, we can write $D_j^{-1}(i)$ in row form before we start working on the next matrix.

$$D_7^{-1}M^-$$

| | 7 | 8 | 11 | 17 | 18 | 14 | 5 | 1 | 4 | 12 | 9 | 20 | 19 | 13 | 16 | 6 | 10 | 15 | 3 | 2 | |
|---|---|---|---|---|---|---|---|---|---|---|---|---|---|---|---|---|---|---|---|---|---|
| | 1 | 2 | 3 | 4 | 5 | 6 | 7 | 8 | 9 | 10 | 11 | 12 | 13 | 14 | 15 | 16 | 17 | 18 | 19 | 20 | |
| 1 | 0 | 50 | 43 | -7 | 63 | 8 | 5 | ∞ | 88 | 42 | 30 | 56 | 77 | 20 | 56 | 29 | 37 | 39 | 63 | 79 | 1 |
| 2 | 79 | 0 | 4 | 34 | 42 | 25 | 47 | 60 | 52 | 56 | 38 | 25 | 5 | 12 | 42 | -2 | 46 | -9 | 8 | ∞ | 2 |
| 3 | 60 | 63 | 0 | 43 | 19 | 81 | 21 | 42 | 69 | 20 | 64 | 68 | 61 | 10 | 74 | 13 | 4 | 63 | ∞ | 27 | 3 |
| 4 | 85 | 79 | 29 | 0 | 58 | 53 | 83 | 20 | ∞ | 50 | 78 | 22 | 84 | 85 | 9 | 54 | 9 | 24 | 68 | 67 | 4 |
| 5 | 2 | 45 | 22 | -9 | 0 | -15 | ∞ | -2 | 61 | -1 | 58 | 16 | 80 | 72 | 10 | 30 | 35 | 16 | 41 | 13 | 5 |
| 6 | 80 | 76 | 12 | 51 | 52 | 0 | 1 | 3 | 57 | 25 | 54 | 53 | 4 | 30 | 64 | ∞ | 22 | 76 | 13 | 1 | 6 |
| 7 | ∞ | 70 | 77 | 87 | 91 | 29 | 0 | -3 | 64 | 55 | 24 | 85 | 54 | 20 | 6 | 31 | 40 | 5 | 78 | 95 | 7 |
| 8 | 53 | ∞ | 24 | 23 | 77 | 64 | 51 | 0 | 48 | 47 | 58 | 69 | 66 | 55 | 58 | 29 | 91 | 33 | 70 | 48 | 8 |
| 9 | 60 | 13 | 84 | 71 | 62 | 15 | 57 | 86 | 0 | 18 | ∞ | 60 | 78 | 86 | 37 | 73 | 57 | -2 | 43 | 93 | 9 |
| 10 | 82 | 82 | 24 | 63 | 89 | 53 | 70 | 36 | -6 | 0 | 62 | -5 | 75 | 76 | 30 | 44 | ∞ | 53 | -8 | 55 | 10 |
| 11 | 14 | 19 | ∞ | 22 | 53 | 53 | 59 | -21 | -14 | 58 | 0 | -13 | -1 | -2 | 34 | 5 | -9 | -14 | 49 | -18 | 11 |
| 12 | 72 | 82 | 70 | 61 | 58 | 35 | -1 | 37 | 57 | ∞ | 32 | 0 | 18 | 71 | 55 | 32 | 37 | 40 | 43 | -6 | 12 |
| 13 | 8 | 74 | 38 | 9 | 60 | 47 | 13 | 12 | 69 | 60 | 8 | -3 | 0 | ∞ | 85 | 73 | 80 | 2 | 86 | 75 | 13 |
| 14 | 12 | 33 | 4 | -13 | 76 | ∞ | 37 | 65 | 71 | 35 | 9 | 7 | 40 | 0 | 18 | 77 | 7 | 72 | 70 | 54 | 14 |
| 15 | 44 | 63 | 80 | 61 | 21 | 26 | 47 | 66 | 24 | 34 | 72 | 73 | 32 | 41 | 0 | 36 | 62 | ∞ | 47 | 35 | 15 |
| 16 | 83 | 12 | 57 | 64 | 56 | 15 | 76 | 8 | 55 | 7 | 53 | 52 | 15 | 54 | ∞ | 0 | 43 | 43 | 72 | 0 | 16 |
| 17 | 82 | 23 | 18 | ∞ | 70 | 27 | 7 | 45 | 58 | 74 | 27 | 5 | 30 | 16 | 74 | 82 | 0 | 81 | 24 | 47 | 17 |
| 18 | 64 | 63 | 79 | 65 | ∞ | 6 | 23 | 28 | 95 | 87 | 80 | 38 | 48 | 39 | 23 | 77 | 77 | 0 | 81 | 47 | 18 |
| 19 | 20 | 48 | 57 | 26 | 59 | 55 | 52 | 47 | 38 | 57 | 58 | 62 | ∞ | 81 | 51 | 42 | 82 | 36 | 0 | 83 | 19 |
| 20 | 8 | 11 | 18 | 85 | 79 | 41 | 43 | 37 | 45 | 74 | 29 | ∞ | 78 | 59 | 49 | 46 | 12 | 0 | 67 | 0 | 20 |
| | 1 | 2 | 3 | 4 | 5 | 6 | 7 | 8 | 9 | 10 | 11 | 12 | 13 | 14 | 15 | 16 | 17 | 18 | 19 | 20 | |

j = 4

  (14 4)(4 15) = (14 15): -4

  (14 4)(4 17) = (14 17): -4

j = 6

  (5 6)(6 3) = (5 3): -3

  (5 6)(6 7) = (5 7): -14

  (5 6)(6 8) = (5 8): -12

  (5 6)(6 13) = (5 13): -11

  (5 6)(6 19) = (5 19): -2

  (5 6)(6 20) = (5 20): -14

j = 7

  (5 7)(7 8) = (5 8): -17

  (5 7)(7 15) = (5 15): -8

  (5 7)(7 18) = (5 18): -9

j = 10

  (5 10)(10 9) = (5 9): -7

j = 12

  (10 12)(12 7) = (10 7): -6

  (10 12)(12 20) = (10 20): -11

  (11 12)(12 7) = (11 7): -14

  (11 12)(12 20) = (11 20): -19

  (13 12)(12 7) = (13 7): -4

  (13 12)(12 20) = (13 20): -9

j = 13

  (5 13)(13 1) = (5 1): -3

j = 14

   (11 14)(14 4) = (11 4): -15

j = 16

   (2 16)(16 20) = (2 20): -2

j = 18

   (2 18)(18 6) = (2 6): -3

   (10 18)(18 6) = (10 6): -2

   (11 18)(18 6) = (11 6): -10

j = 20

   (5 20)(20 1) = (5 1): -6

   (5 20)(20 2) = (5 2): -3

   (5 20)(20 17) = (5 17): -2

   (5 20)(20 18) = (5 18): -14

   (11 20)(20 18) = (11 18): -19

   (12 20)(20 18) = (13 18): -9

$$D_7^{-1}M^-(20)$$

| | 7 | 8 | 11 | 17 | 18 | 14 | 5 | 1 | 4 | 12 | 9 | 20 | 19 | 13 | 16 | 6 | 10 | 15 | 3 | 2 | |
|---|---|---|---|---|---|---|---|---|---|---|---|---|---|---|---|---|---|---|---|---|---|
| | 1 | 2 | 3 | 4 | 5 | 6 | 7 | 8 | 9 | 10 | 11 | 12 | 13 | 14 | 15 | 16 | 17 | 18 | 19 | 20 | |
| 1 | 0 | 50 | 43 | *-7* | 63 | 8 | 5 | ∞ | 88 | 42 | 30 | 56 | 77 | 20 | 56 | 29 | 37 | 39 | 63 | 79 | 1 |
| 2 | 79 | 0 | 4 | 34 | 42 | <u>-3</u> | 47 | 60 | 52 | 56 | 38 | 25 | 5 | 12 | 42 | *-2* | 46 | -9 | 8 | *-2* | 2 |
| 3 | 60 | 63 | 0 | 43 | 19 | 81 | 21 | 42 | 69 | 20 | 64 | 68 | 61 | 10 | 74 | 13 | 4 | 63 | ∞ | 27 | 3 |
| 4 | 85 | 79 | 29 | 0 | 58 | 53 | 83 | 20 | ∞ | 50 | 78 | 22 | 84 | 85 | 9 | 54 | 9 | 24 | 68 | 67 | 4 |
| 5 | <u>-6</u> | <u>-3</u> | <u>-3</u> | -9 | 0 | *-15* | *-14* | *-17* | <u>-7</u> | *-1* | 58 | 16 | *-11* | 72 | *-8* | 30 | <u>-2</u> | *-14* | *-2* | *-14* | 5 |
| 6 | 80 | 76 | 12 | 51 | 52 | 0 | 1 | 3 | 57 | 25 | 54 | 53 | 4 | 30 | 64 | ∞ | 22 | 76 | 13 | 1 | 6 |
| 7 | ∞ | 70 | 77 | 87 | 91 | 29 | 0 | *-3* | 64 | 55 | 24 | 85 | 54 | 20 | 6 | 31 | 40 | 5 | 78 | 95 | 7 |
| 8 | 53 | ∞ | 24 | 23 | 77 | 64 | 51 | 0 | 48 | 47 | 58 | 69 | 66 | 55 | 58 | 29 | 91 | 33 | 70 | 48 | 8 |
| 9 | 60 | 13 | 84 | 71 | 62 | 15 | 57 | 86 | 0 | 18 | ∞ | 60 | 78 | 86 | 37 | 73 | 57 | *-2* | 43 | 93 | 9 |
| 10 | 82 | 82 | 24 | 63 | 89 | <u>-2</u> | <u>-6</u> | 36 | -6 | 0 | 62 | -5 | 75 | 76 | 30 | 44 | ∞ | *-8* | *-8* | *-11* | 10 |
| 11 | 14 | <u>-1</u> | ∞ | <u>-15</u> | 53 | <u>-10</u> | <u>-14</u> | *-21* | *-14* | 58 | 0 | *-13* | *-1* | -2 | 34 | 5 | *-9* | <u>-19</u> | 49 | *-19* | 11 |
| 12 | 72 | 82 | 70 | 61 | 58 | 35 | - | *-4* | 57 | ∞ | 32 | 0 | 18 | 71 | 55 | 32 | 37 | <u>-6</u> | 43 | -6 | 12 |
| 13 | 8 | 74 | 38 | 9 | 60 | 47 | <u>-4</u> | 12 | 69 | 60 | 8 | -3 | 0 | ∞ | *-4* | 73 | *-4* | <u>-9</u> | 86 | *-9* | 13 |
| 14 | 12 | 33 | 4 | *-13* | 76 | ∞ | 37 | 65 | 71 | 35 | 9 | 7 | 40 | 0 | *-4* | 77 | 7 | *-4* | 70 | 54 | 14 |
| 15 | 44 | 63 | 80 | 61 | 21 | 26 | 47 | 66 | 24 | 34 | 72 | 73 | 32 | 41 | 0 | 36 | 62 | ∞ | 47 | 35 | 15 |
| 16 | 83 | 12 | 57 | 64 | 56 | 15 | 76 | 8 | 55 | 7 | 53 | 52 | 15 | 54 | ∞ | 0 | 43 | 43 | 72 | 0 | 16 |
| 17 | 82 | 23 | 18 | ∞ | 70 | 27 | 7 | 45 | 58 | 74 | 27 | 5 | 30 | 16 | 74 | 82 | 0 | 81 | 24 | 47 | 17 |
| 18 | 64 | 63 | 79 | 65 | ∞ | 6 | 23 | 28 | 95 | 87 | 80 | 38 | 48 | 39 | 23 | 77 | 77 | 0 | 81 | 47 | 18 |
| 19 | 20 | 48 | 57 | 26 | 59 | 55 | 52 | 47 | 38 | 57 | 58 | 62 | ∞ | 81 | 51 | 42 | 82 | 36 | 0 | 83 | 19 |
| 20 | 8 | 11 | 18 | 85 | 79 | 41 | 43 | 37 | 45 | 74 | 29 | ∞ | 78 | 59 | 49 | 46 | 12 | 0 | 67 | 0 | 20 |
| | 1 | 2 | 3 | 4 | 5 | 6 | 7 | 8 | 9 | 10 | 11 | 12 | 13 | 14 | 15 | 16 | 17 | 18 | 19 | 20 | |

$P_{20}$

|    | 1  | 2  | 3 | 4  | 5 | 6  | 7  | 8 | 9  | 10 | 11 | 12 | 13 | 14 | 15 | 16 | 17 | 18 | 19 | 20 |    |
|----|----|----|---|----|---|----|----|---|----|----|----|----|----|----|----|----|----|----|----|----|----|
| 1  |    |    |   |    |   |    |    |   |    |    |    |    |    |    |    |    |    |    |    |    | 1  |
| 2  |    |    |   |    |   | 18 |    |   |    |    |    |    |    |    |    | 2  |    | 2  |    | 16 | 2  |
| 3  |    |    |   |    |   |    |    |   |    |    |    |    |    |    |    |    |    |    |    |    | 3  |
| 4  |    |    |   |    |   |    |    |   |    |    |    |    |    |    |    |    |    |    |    |    | 4  |
| 5  | 20 | 20 | 6 |    |   | 5  | 6  | 7 | 10 | 5  |    |    | 6  |    | 7  |    | 20 | 20 | 6  | 6  | 5  |
| 6  |    |    |   |    |   |    |    |   |    |    |    |    |    |    |    |    |    |    |    |    | 6  |
| 7  |    |    |   |    |   |    |    |   |    |    |    |    |    |    |    |    |    |    |    |    | 7  |
| 8  |    |    |   |    |   |    |    |   |    |    |    |    |    |    |    |    |    |    |    |    | 8  |
| 9  |    |    |   |    |   |    |    |   |    |    |    |    |    |    |    |    |    |    |    |    | 9  |
| 10 |    |    |   |    |   | 18 | 12 |   | 10 |    |    | 10 |    |    |    |    |    | 9  |    | 12 | 10 |
| 11 |    | 9  |   | 14 |   | 18 | 12 |   | 11 |    |    | 11 |    | 11 |    |    |    | 20 |    | 12 | 11 |
| 12 |    |    |   |    |   |    |    |   |    |    |    |    |    |    |    |    |    |    |    |    | 12 |
| 13 |    |    |   |    |   |    | 12 |   |    |    |    |    | 13 |    |    |    |    | 20 |    | 12 | 13 |
| 14 |    |    |   | 14 |   |    |    |   |    |    |    |    |    |    | 4  |    | 4  |    |    |    | 14 |
| 15 |    |    |   |    |   |    |    |   |    |    |    |    |    |    |    |    |    |    |    |    | 15 |
| 16 |    |    |   |    |   |    |    |   |    |    |    |    |    |    |    |    |    |    |    |    | 16 |
| 17 |    |    |   |    |   |    |    |   |    |    |    |    |    |    |    |    |    |    |    |    | 17 |
| 18 |    |    |   |    |   |    |    |   |    |    |    |    |    |    |    |    |    |    |    |    | 18 |
| 19 |    |    |   |    |   |    |    |   |    |    |    |    |    |    |    |    |    |    |    |    | 19 |
| 20 |    |    |   |    |   |    |    |   |    |    |    |    |    |    |    |    |    |    |    |    | 20 |
|    | 1  | 2  | 3 | 4  | 5 | 6  | 7  | 8 | 9  | 10 | 11 | 12 | 13 | 14 | 15 | 16 | 17 | 18 | 19 | 20 |    |

*NEGPATHS(20)*

(1) (5 1);  (2) (5 2), (7 2), (11 2);  (3) (5 3);  (4) (11 4);  (6) (2 6), (10 6), (11 6);
(7) (10 7), (11 7), (13 7);  (9) (5 9);  (17) (5 17);  (18) (5 18), (11 18), (12 18), (13 18).

j = 1

  (5 1)(1 4) = (5 4): -13

j = 2

  (5 2)(2 16) = (5 16): -5

  (12 2)(2 16) = (11 16): -3

j = 4

  (5 4)(4 17) = (5 17): -4

j = 6

  (2 6)(6 7) = (2 7): -2

  (11 6)(6 13) = (11 13): -9

j = 7

  (2 7)(7 8) = (2 8): -5

  (10 7)(7 8) = (10 8): -9

j = 13

  (11 13)(13 1) = (11 1): -1

  **(11 13)(13 11) = (11 11): -1**

We thus have obtained a negative cycle using column 13 and need continue no further in the construction of $D_7^{-1}M^-(40)$.

$$D_7^{-1}M^-(40)$$

| | 7 | 8 | 11 | 17 | 18 | 14 | 5 | 1 | 4 | 12 | 9 | 20 | 19 | 13 | 16 | 6 | 10 | 15 | 3 | 2 | |
|---|---|---|---|---|---|---|---|---|---|---|---|---|---|---|---|---|---|---|---|---|---|
| | 1 | 2 | 3 | 4 | 5 | 6 | 7 | 8 | 9 | 10 | 11 | 12 | 13 | 14 | 15 | 16 | 17 | 18 | 19 | 20 | |
| 1 | 0 | 50 | 43 | *-7* | 63 | 8 | 5 | ∞ | 88 | 42 | 30 | 56 | 77 | 20 | 56 | 29 | 37 | 39 | 63 | 79 | 1 |
| 2 | 79 | 0 | 4 | 34 | 42 | *-3* | *-2* | 60 | 52 | 56 | 38 | 25 | 5 | 12 | 42 | *-2* | 46 | -9 | 8 | *-2* | 2 |
| 3 | 60 | 63 | 0 | 43 | 19 | 81 | 21 | 42 | 69 | 20 | 64 | 68 | 61 | 10 | 74 | 13 | 4 | 63 | ∞ | 27 | 3 |
| 4 | 85 | 79 | 29 | 0 | 58 | 53 | 83 | 20 | ∞ | 50 | 78 | 22 | 84 | 85 | 9 | 54 | 9 | 24 | 68 | 67 | 4 |
| 5 | *-6* | *-3* | *-3* | *-13* | 0 | *-15* | *-14* | *-17* | *-7* | *-1* | 58 | 16 | *-11* | 72 | *-8* | <u>*-5*</u> | <u>*-4*</u> | *-14* | *-2* | *-14* | 5 |
| 6 | 80 | 76 | 12 | 51 | 52 | 0 | 1 | 3 | 57 | 25 | 54 | 53 | 4 | 30 | 64 | ∞ | 22 | 76 | 13 | 1 | 6 |
| 7 | ∞ | 70 | 77 | 87 | 91 | 29 | 0 | *-3* | 64 | 55 | 24 | 85 | 54 | 20 | 6 | 31 | 40 | 5 | 78 | 95 | 7 |
| 8 | 53 | ∞ | 24 | 23 | 77 | 64 | 51 | 0 | 48 | 47 | 58 | 69 | 66 | 55 | 58 | 29 | 91 | 33 | 70 | 48 | 8 |
| 9 | 60 | 13 | 84 | 71 | 62 | 15 | 57 | 86 | 0 | 18 | ∞ | 60 | 78 | 86 | 37 | 73 | 57 | -2 | 43 | 93 | 9 |
| 10 | 82 | 82 | 24 | 63 | 89 | *-2* | *-6* | *-9* | *-6* | 0 | 62 | *-5* | 75 | 76 | 30 | 44 | ∞ | *-8* | *-8* | *-11* | 10 |
| 11 | <u>*-1*</u> | *-1* | ∞ | *-15* | 53 | *-13* | *-14* | *-21* | *-14* | 58 | **-1** | *-13* | *-9* | *-2* | <u>*-6*</u> | <u>*-3*</u> | <u>*-11*</u> | <u>*-19*</u> | 49 | *-19* | 11 |
| 12 | 72 | 82 | 70 | 61 | 58 | 35 | - | *-4* | 57 | ∞ | 32 | 0 | 18 | 71 | 55 | 32 | 37 | <u>*-6*</u> | 43 | -6 | 12 |
| 13 | 8 | 74 | 38 | 9 | 60 | 47 | *-4* | *-7* | 69 | 60 | 8 | *-3* | 0 | ∞ | *-4* | 73 | *-4* | <u>*-9*</u> | 86 | *-9* | 13 |
| 14 | 12 | 33 | 4 | *-13* | 76 | ∞ | 37 | 65 | 71 | 35 | 9 | 7 | 40 | 0 | *-4* | 77 | 7 | *-4* | 70 | 54 | 14 |
| 15 | 44 | 63 | 80 | 61 | 21 | 26 | 47 | 66 | 24 | 34 | 72 | 73 | 32 | 41 | 0 | 36 | 62 | ∞ | 47 | 35 | 15 |
| 16 | 83 | 12 | 57 | 64 | 56 | 15 | 76 | 8 | 55 | 7 | 53 | 52 | 15 | 54 | ∞ | 0 | 43 | 43 | 72 | 0 | 16 |
| 17 | 82 | 23 | 18 | ∞ | 70 | 27 | 7 | 45 | 58 | 74 | 27 | 5 | 30 | 16 | 74 | 82 | 0 | 81 | 24 | 47 | 17 |
| 18 | 64 | 63 | 79 | 65 | ∞ | 6 | 23 | 28 | 95 | 87 | 80 | 38 | 48 | 39 | 23 | 77 | 77 | 0 | 81 | 47 | 18 |
| 19 | 20 | 48 | 57 | 26 | 59 | 55 | 52 | 47 | 38 | 57 | 58 | 62 | ∞ | 81 | 51 | 42 | 82 | 36 | 0 | 83 | 19 |
| 20 | 8 | 11 | 18 | 85 | 79 | 41 | 43 | 37 | 45 | 74 | 29 | ∞ | 78 | 59 | 49 | 46 | 12 | 0 | 67 | 0 | 20 |
| | 1 | 2 | 3 | 4 | 5 | 6 | 7 | 8 | 9 | 10 | 11 | 12 | 13 | 14 | 15 | 16 | 17 | 18 | 19 | 20 | |

$P_{40}$

|    | 1  | 2  | 3 | 4 | 5 | 6  | 7  | 8 | 9 | 10 | 11 | 12 | 13 | 14 | 15 | 16 | 17 | 18 | 19 | 20 |
|----|----|----|---|---|---|----|----|---|---|----|----|----|----|----|----|----|----|----|----|----|
| 1  |    |    |   |   |   |    |    |   |   |    |    |    |    |    |    |    |    |    |    |    |
| 2  |    |    |   |   |   | 2  | 6  | 7 |   |    |    |    |    |    |    |    |    |    |    |    |
| 3  |    |    |   |   |   |    |    |   |   |    |    |    |    |    |    |    |    |    |    |    |
| 4  |    |    |   |   |   |    |    |   |   |    |    |    |    |    |    |    |    |    |    |    |
| 5  | 5  | 5  |   | 1 |   |    |    |   |   |    |    |    |    |    |    | 2  | 4  |    |    |    |
| 6  |    |    |   |   |   |    |    |   |   |    |    |    |    |    |    |    |    |    |    |    |
| 7  |    |    |   |   |   |    |    |   |   |    |    |    |    |    |    |    |    |    |    |    |
| 8  |    |    |   |   |   |    |    |   |   |    |    |    |    |    |    |    |    |    |    |    |
| 9  |    |    |   |   |   |    |    |   |   |    |    |    |    |    |    |    |    |    |    |    |
| 10 |    |    |   |   |   | 10 | 7  |   |   |    |    |    |    |    |    |    |    |    |    | 12 |
| 11 | 13 | 11 |   |   |   | 11 |    |   |   |    | **13** |    | 6  |    |    | 2  |    |    |    | 12 |
| 12 |    |    |   |   |   |    |    |   |   |    |    |    |    |    |    |    |    |    |    | 12 |
| 13 |    |    |   |   |   |    |    |   |   |    |    |    |    |    |    |    |    |    |    | 12 |
| 14 |    |    |   |   |   |    |    |   |   |    |    |    |    |    |    |    |    |    |    |    |
| 15 |    |    |   |   |   |    |    |   |   |    |    |    |    |    |    |    |    |    |    |    |
| 16 |    |    |   |   |   |    |    |   |   |    |    |    |    |    |    |    |    |    |    |    |
| 17 |    |    |   |   |   |    |    |   |   |    |    |    |    |    |    |    |    |    |    |    |
| 18 |    |    |   |   |   |    |    |   |   |    |    |    |    |    |    |    |    |    |    |    |
| 19 |    |    |   |   |   |    |    |   |   |    |    |    |    |    |    |    |    |    |    |    |
| 20 |    |    |   |   |   |    |    |   |   |    |    |    |    |    |    |    |    |    |    |    |
|    | 1  | 2  | 3 | 4 | 5 | 6  | 7  | 8 | 9 | 10 | 11 | 12 | 13 | 14 | 15 | 16 | 17 | 18 | 19 | 20 |

From the entry –1 at (11 13) of $D_7^{-1}M^-(40)$, the value of s is –1. Using $P_{40}$, checking for further points of s, we obtain (11 6 13). Going back to $P_{20}$, $18 \to 6, 20 \to 18, 12 \to 20, 11 \to 12$. Thus,

$s$ = (11 12 20 18 6 13). $D_8 = D_7 s \Rightarrow D_8^{-1} = s^{-1}D_7^{-1}$. Applying $s^{-1}$ to the second row of $D_7^{-1}M$, we obtain $D_8^{-1}M$.

$$D_8^{-1}M^-$$

|   | 7 | 8 | 11 | 17 | 18 | 19 | 5 | 1 | 4 | 12 | 20 | 2 | 9 | 13 | 16 | 6 | 10 | 14 | 3 | 15 |   |
|---|---|---|----|----|----|----|---|---|---|----|----|---|---|----|----|---|----|----|---|----|---|
|   | 1 | 2 | 3 | 4 | 5 | 6 | 7 | 8 | 9 | 10 | 11 | 12 | 13 | 14 | 15 | 16 | 17 | 18 | 19 | 20 |   |
| 1 | 0 | 50 | 43 | -7 | 63 | 77 | 5 | ∞ | 88 | 41 | 56 | 79 | 30 | 20 | 56 | 29 | 37 | 8 | 63 | 39 | 1 |
| 2 | 79 | 0 | 4 | 34 | 42 | 5 | 47 | 60 | 52 | 56 | 25 | ∞ | 38 | 12 | 42 | -2 | 46 | 25 | 8 | -9 | 2 |
| 3 | 60 | 63 | 0 | 43 | 19 | 61 | 21 | 42 | 49 | 20 | 68 | 27 | 84 | 10 | 74 | 13 | 4 | 81 | ∞ | 63 | 3 |
| 4 | 85 | 79 | 29 | 0 | 58 | 84 | 83 | 20 | ∞ | 50 | 22 | 67 | 78 | 85 | 9 | 54 | 9 | 53 | 68 | 24 | 4 |
| 5 | 2 | 45 | 22 | -9 | 0 | 80 | ∞ | -2 | 61 | -1 | 16 | 13 | 58 | 72 | 10 | 30 | 35 | -15 | 41 | 16 | 5 |
| 6 | 76 | 72 | 8 | 47 | 48 | 0 | -3 | -1 | 53 | 21 | 49 | -3 | 50 | 26 | 60 | ∞ | 18 | -2 | 9 | 72 | 6 |
| 7 | ∞ | 70 | 77 | 87 | 91 | 54 | 0 | -3 | 64 | 55 | 85 | 95 | 24 | 20 | 6 | 31 | 40 | 29 | 78 | 5 | 7 |
| 8 | 51 | ∞ | 24 | 23 | 77 | 66 | 51 | 0 | 48 | 47 | 69 | 48 | 58 | 55 | 58 | 29 | 91 | 64 | 70 | 33 | 8 |
| 9 | 60 | 13 | 84 | 71 | 62 | 78 | 57 | 86 | 0 | 18 | 60 | 93 | ∞ | 86 | 37 | 73 | 57 | 15 | 43 | -2 | 9 |
| 10 | 82 | 82 | 24 | 63 | 89 | 75 | 70 | 36 | -6 | 0 | -5 | 55 | 62 | 76 | 30 | 44 | ∞ | 53 | -8 | 53 | 10 |
| 11 | 27 | 32 | ∞ | 35 | 66 | 12 | 72 | -8 | -1 | 71 | 0 | -5 | 13 | 11 | 47 | 18 | 4 | 66 | 62 | -1 | 11 |
| 12 | 78 | 88 | 76 | 67 | 64 | 24 | 5 | 42 | 63 | ∞ | 6 | 0 | 38 | 77 | 61 | 38 | 43 | 41 | 51 | 44 | 12 |
| 13 | 0 | 66 | 30 | 1 | 52 | -8 | 5 | 4 | 61 | 52 | -11 | 67 | 0 | ∞ | 77 | 65 | 72 | 39 | 78 | -6 | 13 |
| 14 | 12 | 33 | 4 | -15 | 76 | 40 | 37 | 65 | 71 | 35 | 7 | 54 | 9 | 0 | 18 | 77 | 7 | ∞ | 70 | 72 | 14 |
| 15 | 44 | 63 | 80 | 61 | 21 | 32 | 47 | 66 | 24 | 34 | 73 | 35 | 72 | 41 | 0 | 36 | 62 | 26 | 47 | ∞ | 15 |
| 16 | 83 | 12 | 57 | 64 | 56 | 15 | 76 | 8 | 55 | 7 | 52 | 0 | 53 | 54 | ∞ | 0 | 43 | 15 | 72 | 43 | 16 |
| 17 | 82 | 23 | 18 | ∞ | 70 | 30 | 7 | 45 | 58 | 74 | 5 | 47 | 27 | 16 | 74 | 82 | 0 | 27 | 24 | 81 | 17 |
| 18 | 58 | 42 | 73 | 59 | ∞ | 42 | 17 | 22 | 89 | 81 | 32 | 41 | 74 | 33 | 17 | 71 | 71 | 0 | 75 | -6 | 18 |
| 19 | 20 | 11 | 57 | 26 | 53 | ∞ | 52 | 47 | 38 | 57 | 62 | 83 | 58 | 81 | 51 | 42 | 82 | 55 | 0 | -6 | 19 |
| 20 | 8 | 50 | 18 | 85 | 79 | 78 | 43 | 37 | 45 | 74 | ∞ | 0 | 29 | 59 | 49 | 46 | 12 | 41 | 67 | 0 | 20 |
|   | 1 | 2 | 3 | 4 | 5 | 6 | 7 | 8 | 9 | 10 | 11 | 12 | 13 | 14 | 15 | 16 | 17 | 18 | 19 | 20 |   |

j = 4

   (14 4)(4 15) = (14 15): -6

   (14 4)(4 17) = (14 17): -6

j = 6

   (13 6)(6 7) = (13 7): -17

   (13 6)(6 8) = (13 8): -9

   (13 6)(6 12) = (13 12): -11

   (13 6)(6 18) = (13 18): -10

j = 7

   (6 7)(7 8) = (6 8): -6

   (13 7)(7 8) = (13 8): -20

   (13 7)(7 15) = (13 15): -11

   (13 7)(7 20) = (13 20): -12

j = 9

   (10 9)(9 20) = (10 20): -8

   (11 9)(9 20) = (11 20): -3

   (10 9)(9 20) = (10 20): -8

j = 10

   (5 10)(10 9) = (5 9): -7

   (5 10)(10 11) = (5 11): -6

   (5 10)(10 19) = (5 19): -9

j = 11

   (5 11)(11 8) = (5 8): -14

   (5 11)(11 9) = (5 9): -7

   (5 11)(11 12) = (5 12): -11

   (5 11)(11 17) = (5 17): -2

   (5 11)(11 20) = (5 20): -7

   (10 11)(11 12) = (10 12) : -10

(10  11)(11  17) = (10  17): -1

(13  11)(11  9) = (13  9): -12

(13  11)(11 12) = (13  12): -16

(13  11)(11  17) = (13  17): -7

(13  11)(11  20) = (13  20): -12

j = 12

(5  12)(12  7) = (5  7): -6

(10 12)(12  7) = (10  7): -5

j = 18

(5  18)(18  20) = (5  20): -21

(13  18)(18  20) = (13  20): -16

j = 20

(2  20)(20  1) = (2  1): -1

(2  20)(20  12) = (2  12): -9

(5  20)(20  1) = (5  1): -13

(5  20)(20  3) = (5  3): -3

(5  20)(20  12) = (5  12): -21

(5  20)(20  17) = (5  17): -9

(10  20)(20  1) = (10  1): -6

(10  20)(20  12) = (10  12): -8

(11  20)(20  12) = (11  12): -3

(13  20)(20  1) = (13  1): -6

(13  20)(20  12) = (13  12): -14

(13  20)(20  17) = (13  17): -2

(18  20)(20  12) = (18  12): -6

(19  20)(20  12) = (19  12): -6

$$D_8^{-1} M^- (20)$$

| | 7 | 8 | 11 | 17 | 18 | 19 | 5 | 1 | 4 | 12 | 20 | 2 | 9 | 13 | 16 | 6 | 10 | 14 | 3 | 15 | |
|---|---|---|---|---|---|---|---|---|---|---|---|---|---|---|---|---|---|---|---|---|---|
| | 1 | 2 | 3 | 4 | 5 | 6 | 7 | 8 | 9 | 10 | 11 | 12 | 13 | 14 | 15 | 16 | 17 | 18 | 19 | 20 | |
| 1 | 0 | 50 | 43 | -7 | 63 | 77 | 5 | ∞ | 88 | 41 | 56 | 79 | 30 | 20 | 56 | 29 | 37 | 8 | 63 | 39 | 1 |
| 2 | <u>-1</u> | 0 | 4 | 34 | 42 | 5 | 47 | 60 | 52 | 56 | 25 | <u>-9</u> | 38 | 12 | 42 | -2 | 46 | 25 | 8 | *-9* | 2 |
| 3 | 60 | 63 | 0 | 43 | 19 | 61 | 21 | 42 | 49 | 20 | 68 | 27 | 84 | 10 | 74 | 13 | 4 | 81 | ∞ | 63 | 3 |
| 4 | 85 | 79 | 29 | 0 | 58 | 84 | 83 | 20 | ∞ | 50 | 22 | 67 | 78 | 85 | 9 | 54 | 9 | 53 | 68 | 24 | 4 |
| 5 | <u>-13</u> | 45 | <u>-3</u> | -9 | 0 | 80 | <u>-6</u> | <u>-14</u> | <u>-7</u> | <u>-7</u> | -6 | <u>-21</u> | 58 | 72 | 10 | 30 | <u>-9</u> | *-15* | *-9* | *-21* | 5 |
| 6 | 76 | 72 | 8 | 47 | 48 | 0 | -3 | <u>-6</u> | 53 | 21 | 49 | -3 | 50 | 26 | 60 | ∞ | 18 | -2 | 9 | 72 | 6 |
| 7 | ∞ | 70 | 77 | 87 | 91 | 54 | 0 | *-3* | 64 | 55 | 85 | 95 | 24 | 20 | 6 | 31 | 40 | 29 | 78 | 5 | 7 |
| 8 | 51 | ∞ | 24 | 23 | 77 | 66 | 51 | 0 | 48 | 47 | 69 | 48 | 58 | 55 | 58 | 29 | 91 | 64 | 70 | 33 | 8 |
| 9 | 60 | 13 | 84 | 71 | 62 | 78 | 57 | 86 | 0 | 18 | 60 | 93 | ∞ | 86 | 37 | 73 | 57 | 15 | 43 | -2 | 9 |
| 10 | <u>-6</u> | 82 | 24 | 63 | 89 | 75 | <u>-5</u> | 36 | -6 | 0 | -5 | <u>-14</u> | 62 | 76 | 30 | 44 | <u>-2</u> | 53 | *-8* | *-14* | 10 |
| 11 | 27 | 32 | ∞ | 35 | 66 | 12 | 72 | *-8* | *-1* | 71 | 0 | -5 | 13 | 11 | 47 | 18 | 4 | 66 | 62 | *-3* | 11 |
| 12 | 78 | 88 | 76 | 67 | 64 | 24 | 5 | 42 | 63 | ∞ | 6 | 0 | 38 | 77 | 61 | 38 | 43 | 41 | 51 | 44 | 12 |
| 13 | <u>-8</u> | 66 | 30 | 1 | 52 | -8 | *-17* | *-20* | <u>-12</u> | 52 | *-11* | <u>-16</u> | 0 | ∞ | *-11* | 65 | <u>-4</u> | *-10* | 78 | *-16* | 13 |
| 14 | 12 | 33 | 4 | *-15* | 76 | 40 | 37 | 65 | 71 | 35 | <u>-1</u> | 54 | 9 | 0 | -6 | 77 | -6 | ∞ | 70 | 72 | 14 |
| 15 | 44 | 63 | 80 | 61 | 21 | 32 | 47 | 66 | 24 | 34 | 73 | 35 | 72 | 41 | 0 | 36 | 62 | 26 | 47 | ∞ | 15 |
| 16 | 83 | 12 | 57 | 64 | 56 | 15 | 76 | 8 | 55 | 7 | 52 | 0 | 53 | 54 | ∞ | 0 | 43 | 15 | 72 | 43 | 16 |
| 17 | 82 | 23 | 18 | ∞ | 70 | 30 | 7 | 45 | 58 | 74 | 5 | 47 | 27 | 16 | 74 | 82 | 0 | 27 | 24 | 81 | 17 |
| 18 | 58 | 42 | 73 | 59 | ∞ | 42 | 17 | 22 | 89 | 81 | 32 | <u>-6</u> | 74 | 33 | 17 | 71 | 71 | 0 | 75 | -6 | 18 |
| 19 | 20 | 11 | 57 | 26 | 53 | ∞ | 52 | 47 | 38 | 57 | 62 | <u>-6</u> | 58 | 81 | 51 | 42 | 82 | 55 | 0 | -6 | 19 |
| 20 | 8 | 50 | 18 | 85 | 79 | 78 | 43 | 37 | 45 | 74 | ∞ | 0 | 29 | 59 | 49 | 46 | 12 | 41 | 67 | 0 | 20 |
| | 1 | 2 | 3 | 4 | 5 | 6 | 7 | 8 | 9 | 10 | 11 | 12 | 13 | 14 | 15 | 16 | 17 | 18 | 19 | 20 | |

$P_{20}$

| | 1 | 2 | 3 | 4 | 5 | 6 | 7 | 8 | 9 | 10 | 11 | 12 | 13 | 14 | 15 | 16 | 17 | 18 | 19 | 20 | |
|---|---|---|---|---|---|---|---|---|---|---|---|---|---|---|---|---|---|---|---|---|---|
| 1 | | | | | | | | | | | | | | | | | | | | | 1 |
| 2 | 2 | | | 1 | | | 12 | | | | 12 | 2 | | | | | | | | | 2 |
| 3 | | | | | | | | | | | | | | | | | | | | | 3 |
| 4 | | | | | | | | | | | | | | | | | | | | | 4 |
| 5 | 5 | 19 | | | | | 12 | 11 | 10 | 5 | 12 | 5 | | | 7 | | | | 10 | | 5 |
| 6 | | | | | | | | | | | | | | | | | | | | | 6 |
| 7 | | | | | | | | | | | | | | | | | | | | | 7 |
| 8 | | | | | | | | | | | | | | | | | | | | | 8 |
| 9 | | | | | | | | | | | | | | | | | | | | | 9 |
| 10 | 10 | | | 1 | | | 12 | 10 | 7 | | 12 | 10 | | | 4 | | 4 | | | | 10 |
| 11 | | | | | | | | | | | | | | | | | | | | | 11 |
| 12 | | | | | | | | | | | | | | | | | | | | | 12 |
| 13 | 13 | | | 1 | | | | | | | | | | | | | 4 | | | | 13 |
| 14 | | | | 14 | | | | 11 | | | | 14 | | | 4 | | 4 | | | | 14 |
| 15 | | | | | | | | | | | | | | | | | | | | | 15 |
| 16 | | | | | | | | | | | | | | | | | | | | | 16 |
| 17 | | | | | | | | | | | | | | | | | | | | | 17 |
| 18 | | | | | | | | | | | | | | | | | | | | | 18 |
| 19 | | | | | | | | | | | | | | | | | | | | | 19 |
| 20 | | | | | | | | | | | | | | | | | | | | | 20 |
| | 1 | 2 | 3 | 4 | 5 | 6 | 7 | 8 | 9 | 10 | 11 | 12 | 13 | 14 | 15 | 16 | 17 | 18 | 19 | 20 | |

*NEGPATHS(20)*

(1) (2 1), (5 1), (10 1), (13 1); (3) (5 3); (7) (5 7), (10 7); (8) (5 8), (6 8), (9) (5 9), (13 9); (10) (5 10); (11) (14 11); (12) (2 12), (5 12), (10 12), (13 12), (18 12), (19 12); (17) (5 17), (10 17), (13 17).

j = 1

   (2 1)(1 4) = (2 4): -8

   (5 1)(1 4) = (5 4): -20

   (5 1)(1 7) = (5 7): -7

   (10 1)(1 4) = (10 4): -13

   (13 1)(1 4) = (13 4): -15

j = 4

   (5 4)(4 15) = (5 15): -11

   (5 4)(4 17) = (5 17): -11

   (10 4)(4 15) = (10 15): -4

   (10 4)(4 17) = (10 17): -4

   (12 4)(4 17) = (13 17); -6

j = 7

   (10 7)(7 8) = (10 8): -8

j = 10

   (5 10)(10 9) = (5 9): -13

   (5 10)(10 11) = (5 11): -12

   (5 10)(10 19) = (5 19): -16

j = 11

   (5 11)(11 8) = (5 8): -20

   (5 11)(11 14) = (5 14): -1

   (13 11)(11 9) = (14 9): -2

j = 12

   (5 12)(12 7) = (5 7): -16

   (5 12)(12 11) = (5 11): -15

   (10 12)(12 7) = (10 7): -9

   (10 12)(12 11) = (10 11): -8

j = 19

   (5 19)(19 2) = (5 2): -5

   (5 19)(19 20) = (5 19)(19 20): -5

j = 20

   (5 20)(20 1) = (5 1): -14

   (5 20)(20 12) = (5 12): -22

$$D_8^{-1} M^- (40)$$

| | 7 | 8 | 11 | 17 | 18 | 19 | 5 | 1 | 4 | 12 | 20 | 2 | 9 | 13 | 16 | 6 | 10 | 14 | 3 | 15 | |
|---|---|---|---|---|---|---|---|---|---|---|---|---|---|---|---|---|---|---|---|---|---|
| | 1 | 2 | 3 | 4 | 5 | 6 | 7 | 8 | 9 | 10 | 11 | 12 | 13 | 14 | 15 | 16 | 17 | 18 | 19 | 20 | |
| 1 | 0 | 50 | 43 | -7 | 63 | 77 | 5 | ∞ | 88 | 41 | 56 | 79 | 30 | 20 | 56 | 29 | 37 | 8 | 63 | 39 | 1 |
| 2 | *-1* | 0 | 4 | *-8* | 42 | 5 | *-4* | 60 | 52 | 56 | *-3* | *-9* | 38 | 12 | 42 | *-2* | 46 | 25 | 8 | *-9* | 2 |
| 3 | 60 | 63 | 0 | 43 | 19 | 61 | 21 | 42 | 49 | 20 | 68 | 27 | 84 | 10 | 74 | 13 | 4 | 81 | ∞ | 63 | 3 |
| 4 | 85 | 79 | 29 | 0 | 58 | 84 | 83 | 20 | ∞ | 50 | 22 | 67 | 78 | 85 | 9 | 54 | 9 | 53 | 68 | 24 | 4 |
| 5 | *-14* | *-5* | *-3* | *-20* | 0 | 80 | *-16* | *-20* | *-13* | *-7* | *-16* | *-22* | 58 | *-1* | *-11* | 30 | *-11* | *-15* | *-16* | *-22* | 5 |
| 6 | 76 | 72 | 8 | 47 | 48 | 0 | *-3* | -6 | 53 | 21 | 49 | *-3* | 50 | 26 | 60 | ∞ | 18 | *-2* | 9 | 72 | 6 |
| 7 | ∞ | 70 | 77 | 87 | 91 | 54 | 0 | -3 | 64 | 55 | 85 | 95 | 24 | 20 | 6 | 31 | 40 | 29 | 78 | 5 | 7 |
| 8 | 51 | ∞ | 24 | 23 | 77 | 66 | 51 | 0 | 48 | 47 | 69 | 48 | 58 | 55 | 58 | 29 | 91 | 64 | 70 | 33 | 8 |
| 9 | 60 | 13 | 84 | 71 | 62 | 78 | 57 | 86 | 0 | 18 | 60 | 93 | ∞ | 86 | 37 | 73 | 57 | 15 | 43 | *-2* | 9 |
| 10 | *-6* | 82 | 24 | *-13* | 89 | 75 | *-9* | *-8* | -6 | 0 | *-8* | *-14* | 62 | 76 | *-4* | 44 | *-4* | 53 | *-8* | *-14* | 10 |
| 11 | 27 | 32 | ∞ | 35 | 66 | 12 | 72 | *-8* | *-1* | 71 | 0 | -5 | 13 | 11 | 47 | 18 | 4 | 66 | 62 | *-3* | 11 |
| 12 | 78 | 88 | 76 | 67 | 64 | 24 | 5 | 42 | 63 | ∞ | 6 | 0 | 38 | 77 | 61 | 38 | 43 | 41 | 51 | 44 | 12 |
| 13 | *-8* | 66 | 30 | *-15* | 52 | *-8* | *-17* | *-20* | *-12* | 52 | *-11* | *-16* | 0 | ∞ | *-11* | 65 | *-6* | *-10* | 78 | *-16* | 13 |
| 14 | 12 | 33 | 4 | *-15* | 76 | 40 | 37 | 65 | *-2* | 35 | *-1* | 54 | 9 | 0 | -6 | 77 | -6 | ∞ | 70 | 72 | 14 |
| 15 | 44 | 63 | 80 | 61 | 21 | 32 | 47 | 66 | 24 | 34 | 73 | 35 | 72 | 41 | 0 | 36 | 62 | 26 | 47 | ∞ | 15 |
| 16 | 83 | 12 | 57 | 64 | 56 | 15 | 76 | 8 | 55 | 7 | 52 | 0 | 53 | 54 | ∞ | 0 | 43 | 15 | 72 | 43 | 16 |
| 17 | 82 | 23 | 18 | ∞ | 70 | 30 | 7 | 45 | 58 | 74 | 5 | 47 | 27 | 16 | 74 | 82 | 0 | 27 | 24 | 81 | 17 |
| 18 | 58 | 42 | 73 | 59 | ∞ | 42 | 17 | 22 | 89 | 81 | 32 | *-6* | 74 | 33 | 17 | 71 | 71 | 0 | 75 | *-6* | 18 |
| 19 | 20 | 11 | 57 | 26 | 53 | ∞ | 52 | 47 | 38 | 57 | 62 | *-6* | 58 | 81 | 51 | 42 | 82 | 55 | 0 | *-6* | 19 |
| 20 | 8 | 50 | 18 | 85 | 79 | 78 | 43 | 37 | 45 | 74 | ∞ | 0 | 29 | 59 | 49 | 46 | 12 | 41 | 67 | 0 | 20 |
| | 1 | 2 | 3 | 4 | 5 | 6 | 7 | 8 | 9 | 10 | 11 | 12 | 13 | 14 | 15 | 16 | 17 | 18 | 19 | 20 | |

$P_{40}$

| | 1 | 2 | 3 | 4 | 5 | 6 | 7 | 8 | 9 | 10 | 11 | 12 | 13 | 14 | 15 | 16 | 17 | 18 | 19 | 20 | |
|---|---|---|---|---|---|---|---|---|---|---|---|---|---|---|---|---|---|---|---|---|---|
| 1 | | | | | | | | | | | | | | | | | | | | | 1 |
| 2 | 2 | | | 1 | | | | | | | | | | | | | | | | | 2 |
| 3 | | | | | | | | | | | | | | | | | | | | | 3 |
| 4 | | | | | | | | | | | | | | | | | | | | | 4 |
| 5 | 20 | 19 | | 1 | | | 12 | 11 | 10 | 5 | 12 | 20 | | 11 | 4 | | 4 | | 10 | 19 | 5 |
| 6 | | | | | | | | | | | | | | | | | | | | | 6 |
| 7 | | | | | | | | | | | | | | | | | | | | | 7 |
| 8 | | | | | | | | | | | | | | | | | | | | | 8 |
| 9 | | | | | | | | | | | | | | | | | | | | | 9 |
| 10 | 10 | | | 1 | | | 10 | 7 | | | 12 | 10 | | | 4 | | 4 | | | | 10 |
| 11 | | | | | | | 12 | | | | | 10 | | | | | | | | | 11 |
| 12 | | | | | | | | | | | | | | | | | | | | | 12 |
| 13 | 13 | | | 1 | | | | | | | | | | | | | 4 | | | | 13 |
| 14 | | | | | | | | | 11 | | 14 | | | | | | | | | | 14 |
| 15 | | | | | | | | | | | | | | | | | | | | | 15 |
| 16 | | | | | | | | | | | | | | | | | | | | | 16 |
| 17 | | | | | | | | | | | | | | | | | | | | | 17 |
| 18 | | | | | | | | | | | | | | | | | | | | | 18 |
| 19 | | | | | | | | | | | | | | | | | | | | | 19 |
| 20 | | | | | | | | | | | | | | | | | | | | | 20 |
| | 1 | 2 | 3 | 4 | 5 | 6 | 7 | 8 | 9 | 10 | 11 | 12 | 13 | 14 | 15 | 16 | 17 | 18 | 19 | 20 | |

## NEGPATHS(40)

(1) (5 1); (2) (5 2); (4) (2 4); (7) (5 7), (10 7); (8) (5 8); (9) (5 9), (14 9); (11) (5 11), (10 11); (12) (5 12).

j = 1

   (5 1)(1 4) = (5 4): -21

j = 2

   (5 2)(2 16) = (5 16): -7

j = 4

   (5 4)(4 15) = (5 15): -12

j = 7

   (10 7)(7 8) = (10 8): -12

j = 9

   (14 9)(9 20) = (14 20): -4

j = 11

   (5 11)(11 8) = (5 8): -24

   (5 11)(11 9) = (5 9): -17

   (10 11)(11 8) = (10 8): -16

   (10 11)(11 9) = (10 9): -9

j = 12

   (5 12)(12 7): -17

j = 20

   (14 20)(20 12): -4

$$D_8^{-1}M^-(60)$$

|     | 7 | 8 | 11 | 17 | 18 | 19 | 5 | 1 | 4 | 12 | 20 | 2 | 9 | 13 | 16 | 6 | 10 | 14 | 3 | 15 |     |
|-----|---|---|----|----|----|----|---|---|---|----|----|---|---|----|----|---|----|----|---|----|-----|
|     | 1 | 2 | 3  | 4  | 5  | 6  | 7 | 8 | 9 | 10 | 11 | 12| 13| 14 | 15 | 16| 17 | 18 | 19| 20 |     |
| 1   | 0 | 50 | 43 | -7 | 63 | 77 | 5 | ∞ | 88 | 41 | 56 | 79 | 30 | 20 | 56 | 29 | 37 | 8 | 63 | 39 | 1 |
| 2   | *-1* | 0 | 4 | *-8* | 42 | 5 | *-4* | 60 | 52 | 56 | *-3* | *-9* | 38 | 12 | 42 | *-2* | 46 | 25 | 8 | *-9* | 2 |
| 3   | 60 | 63 | 0 | 43 | 19 | 61 | 21 | 42 | 49 | 20 | 68 | 27 | 84 | 10 | 74 | 13 | 4 | 81 | ∞ | 63 | 3 |
| 4   | 85 | 79 | 29 | 0 | 58 | 84 | 83 | 20 | ∞ | 50 | 22 | 67 | 78 | 85 | 9 | 54 | 9 | 53 | 68 | 24 | 4 |
| 5   | *-14* | *-5* | *-3* | *-21* | 0 | 80 | *-17* | *-24* | *-17* | *-7* | *-16* | *-22* | 58 | *-1* | *-12* | *-7* | *-12* | *-15* | *-16* | *-22* | 5 |
| 6   | 76 | 72 | 8 | 47 | 48 | 0 | *-3* | *-6* | 53 | 21 | 49 | *-3* | 50 | 26 | 60 | ∞ | 18 | *-2* | 9 | 72 | 6 |
| 7   | ∞ | 70 | 77 | 87 | 91 | 54 | 0 | *-3* | 64 | 55 | 85 | 95 | 24 | 20 | 6 | 31 | 40 | 29 | 78 | 5 | 7 |
| 8   | 51 | ∞ | 24 | 23 | 77 | 66 | 51 | 0 | 48 | 47 | 69 | 48 | 58 | 55 | 58 | 29 | 91 | 64 | 70 | 33 | 8 |
| 9   | 60 | 13 | 84 | 71 | 62 | 78 | 57 | 86 | 0 | 18 | 60 | 93 | ∞ | 86 | 37 | 73 | 57 | 15 | 43 | *-2* | 9 |
| 10  | *-6* | 82 | 24 | *-13* | 89 | 75 | *-9* | *-12* | *-9* | 0 | *-8* | *-14* | 62 | 76 | *-4* | 44 | *-4* | 53 | *-8* | *-14* | 10 |
| 11  | 27 | 32 | ∞ | 35 | 66 | 12 | 72 | *-8* | *-1* | 71 | 0 | *-5* | 13 | 11 | 47 | 18 | 4 | 66 | 62 | *-3* | 11 |
| 12  | 78 | 88 | 76 | 67 | 64 | 24 | 5 | 42 | 63 | ∞ | 6 | 0 | 38 | 77 | 61 | 38 | 43 | 41 | 51 | 44 | 12 |
| 13  | *-8* | 66 | 30 | *-15* | 52 | *-8* | *-17* | *-20* | *-12* | 52 | *-11* | *-16* | 0 | ∞ | *-11* | 65 | *-6* | *-10* | 78 | *-16* | 13 |
| 14  | 12 | 33 | 4 | *-15* | 76 | 40 | 37 | 65 | *-2* | 35 | *-1* | *-4* | 9 | 0 | *-6* | 77 | *-6* | ∞ | 70 | *-4* | 14 |
| 15  | 44 | 63 | 80 | 61 | 21 | 32 | 47 | 66 | 24 | 34 | 73 | 35 | 72 | 41 | 0 | 36 | 62 | 26 | 47 | ∞ | 15 |
| 16  | 83 | 12 | 57 | 64 | 56 | 15 | 76 | 8 | 55 | 7 | 52 | 0 | 53 | 54 | ∞ | 0 | 43 | 15 | 72 | 43 | 16 |
| 17  | 82 | 23 | 18 | ∞ | 70 | 30 | 7 | 45 | 58 | 74 | 5 | 47 | 27 | 16 | 74 | 82 | 0 | 27 | 24 | 81 | 17 |
| 18  | 58 | 42 | 73 | 59 | ∞ | 42 | 17 | 22 | 89 | 81 | 32 | *-6* | 74 | 33 | 17 | 71 | 71 | 0 | 75 | *-6* | 18 |
| 19  | 20 | 11 | 57 | 26 | 53 | ∞ | 52 | 47 | 38 | 57 | 62 | *-6* | 58 | 81 | 51 | 42 | 82 | 55 | 0 | *-6* | 19 |
| 20  | 8 | 50 | 18 | 85 | 79 | 78 | 43 | 37 | 45 | 74 | ∞ | 0 | 29 | 59 | 49 | 46 | 12 | 41 | 67 | 0 | 20 |
|     | 1 | 2 | 3 | 4 | 5 | 6 | 7 | 8 | 9 | 10 | 11 | 12 | 13 | 14 | 15 | 16 | 17 | 18 | 19 | 20 |     |

We can't extend any path further as shown in $D_8^{-1}M^-(80)$. It follows that $\sigma_{APOPT} = D_8$.

$P_{60}$

|    | 1 | 2 | 3 | 4 | 5 | 6 | 7  | 8  | 9  | 10 | 11 | 12 | 13 | 14 | 15 | 16 | 17 | 18 | 19 | 20 |    |
|----|---|---|---|---|---|---|----|----|----|----|----|----|----|----|----|----|----|----|----|----|----|
| 1  |   |   |   |   |   |   |    |    |    |    |    |    |    |    |    |    |    |    |    |    | 1  |
| 2  |   |   |   |   |   |   |    |    |    |    |    |    |    |    |    |    |    |    |    |    | 2  |
| 3  |   |   |   |   |   |   |    |    |    |    |    |    |    |    |    |    |    |    |    |    | 3  |
| 4  |   |   |   |   |   |   |    |    |    |    |    |    |    |    |    |    |    |    |    |    | 4  |
| 5  | 5 | 5 |   | 1 |   |   | 12 | 11 | 11 |    | 5  | 5  |    |    | 4  | 2  | 4  |    |    |    | 5  |
| 6  |   |   |   |   |   |   |    |    |    |    |    |    |    |    |    |    |    |    |    |    | 6  |
| 7  |   |   |   |   |   |   |    |    |    |    |    |    |    |    |    |    |    |    |    |    | 7  |
| 8  |   |   |   |   |   |   |    |    |    |    |    |    |    |    |    |    |    |    |    |    | 8  |
| 9  |   |   |   |   |   |   |    |    |    |    |    |    |    |    |    |    |    |    |    |    | 9  |
| 10 |   |   |   |   |   |   | 10 | 11 | 11 |    | 10 |    |    |    |    |    |    |    |    |    | 10 |
| 11 |   |   |   |   |   |   |    |    |    |    |    |    |    |    |    |    |    |    |    |    | 11 |
| 12 |   |   |   |   |   |   |    |    |    |    |    |    |    |    |    |    |    |    |    |    | 12 |
| 13 |   |   |   |   |   |   |    |    |    |    |    |    |    |    |    |    |    |    |    |    | 13 |
| 14 |   |   |   |   |   |   |    |    | 14 |    |    | 20 |    |    |    |    |    |    |    | 9  | 14 |
| 15 |   |   |   |   |   |   |    |    |    |    |    |    |    |    |    |    |    |    |    |    | 15 |
| 16 |   |   |   |   |   |   |    |    |    |    |    |    |    |    |    |    |    |    |    |    | 16 |
| 17 |   |   |   |   |   |   |    |    |    |    |    |    |    |    |    |    |    |    |    |    | 17 |
| 18 |   |   |   |   |   |   |    |    |    |    |    |    |    |    |    |    |    |    |    |    | 18 |
| 19 |   |   |   |   |   |   |    |    |    |    |    |    |    |    |    |    |    |    |    |    | 19 |
| 20 |   |   |   |   |   |   |    |    |    |    |    |    |    |    |    |    |    |    |    |    | 20 |
|    | 1 | 2 | 3 | 4 | 5 | 6 | 7  | 8  | 9  | 10 | 11 | 12 | 13 | 14 | 15 | 16 | 17 | 18 | 19 | 20 |    |

$$D_8^{-1}M^-(80)$$

| | 7 | 8 | 11 | 17 | 18 | 19 | 5 | 1 | 4 | 12 | 20 | 2 | 9 | 13 | 16 | 6 | 10 | 14 | 3 | 15 | |
|---|---|---|---|---|---|---|---|---|---|---|---|---|---|---|---|---|---|---|---|---|---|
| | 1 | 2 | 3 | 4 | 5 | 6 | 7 | 8 | 9 | 10 | 11 | 12 | 13 | 14 | 15 | 16 | 17 | 18 | 19 | 20 | |
| 1 | 0 | 50 | 43 | -7 | 63 | 77 | 5 | ∞ | 88 | 41 | 56 | 79 | 30 | 20 | 56 | 29 | 37 | 8 | 63 | 39 | 1 |
| 2 | *-1* | 0 | 4 | *-8* | 42 | 5 | *-4* | 60 | 52 | 56 | *-3* | *-9* | 38 | 12 | 42 | *-2* | 46 | 25 | 8 | *-9* | 2 |
| 3 | 60 | 63 | 0 | 43 | 19 | 61 | 21 | 42 | 49 | 20 | 68 | 27 | 84 | 10 | 74 | 13 | 4 | 81 | ∞ | 63 | 3 |
| 4 | 85 | 79 | 29 | 0 | 58 | 84 | 83 | 20 | ∞ | 50 | 22 | 67 | 78 | 85 | 9 | 54 | 9 | 53 | 68 | 24 | 4 |
| 5 | *-14* | *-5* | *-3* | *-21* | 0 | 80 | *-17* | *-24* | *-17* | *-7* | *-16* | *-22* | 58 | *-1* | *-12* | *-7* | *-12* | *-15* | *-16* | *-22* | 5 |
| 6 | 76 | 72 | 8 | 47 | 48 | 0 | *-3* | *-6* | 53 | 21 | 49 | *-3* | 50 | 26 | 60 | ∞ | 18 | *-2* | 9 | 72 | 6 |
| 7 | ∞ | 70 | 77 | 87 | 91 | 54 | 0 | *-3* | 64 | 55 | 85 | 95 | 24 | 20 | 6 | 31 | 40 | 29 | 78 | 5 | 7 |
| 8 | 51 | ∞ | 24 | 23 | 77 | 66 | 51 | 0 | 48 | 47 | 69 | 48 | 58 | 55 | 58 | 29 | 91 | 64 | 70 | 33 | 8 |
| 9 | 60 | 13 | 84 | 71 | 62 | 78 | 57 | 86 | 0 | 18 | 60 | 93 | ∞ | 86 | 37 | 73 | 57 | 15 | 43 | *-2* | 9 |
| 10 | *-6* | 82 | 24 | *-13* | 89 | 75 | *-9* | *-12* | *-9* | 0 | *-8* | *-14* | 62 | 76 | *-4* | 44 | *-4* | 53 | *-8* | *-14* | 10 |
| 11 | 27 | 32 | ∞ | 35 | 66 | 12 | 72 | *-8* | *-1* | 71 | 0 | *-5* | 13 | 11 | 47 | 18 | 4 | 66 | 62 | *-3* | 11 |
| 12 | 78 | 88 | 76 | 67 | 64 | 24 | 5 | 42 | 63 | ∞ | 6 | 0 | 38 | 77 | 61 | 38 | 43 | 41 | 51 | 44 | 12 |
| 13 | *-8* | 66 | 30 | *-15* | 52 | *-8* | *-17* | *-20* | *-12* | 52 | *-11* | *-16* | 0 | ∞ | *-11* | 65 | *-6* | *-10* | 78 | *-16* | 13 |
| 14 | 12 | 33 | 4 | *-15* | 76 | 40 | 37 | 65 | *-2* | 35 | *-1* | *-4* | 9 | 0 | *-6* | 77 | *-6* | ∞ | 70 | *-4* | 14 |
| 15 | 44 | 63 | 80 | 61 | 21 | 32 | 47 | 66 | 24 | 34 | 73 | 35 | 72 | 41 | 0 | 36 | 62 | 26 | 47 | ∞ | 15 |
| 16 | 83 | 12 | 57 | 64 | 56 | 15 | 76 | 8 | 55 | 7 | 52 | 0 | 53 | 54 | ∞ | 0 | 43 | 15 | 72 | 43 | 16 |
| 17 | 82 | 23 | 18 | ∞ | 70 | 30 | 7 | 45 | 58 | 74 | 5 | 47 | 27 | 16 | 74 | 82 | 0 | 27 | 24 | 81 | 17 |
| 18 | 58 | 42 | 73 | 59 | ∞ | 42 | 17 | 22 | 89 | 81 | 32 | *-6* | 74 | 33 | 17 | 71 | 71 | 0 | 75 | *-6* | 18 |
| 19 | 20 | 11 | 57 | 26 | 53 | ∞ | 52 | 47 | 38 | 57 | 62 | *-6* | 58 | 81 | 51 | 42 | 82 | 55 | 0 | *-6* | 19 |
| 20 | 8 | 50 | 18 | 85 | 79 | 78 | 43 | 37 | 45 | 74 | ∞ | 0 | 29 | 59 | 49 | 46 | 12 | 41 | 67 | 0 | 20 |
| | 1 | 2 | 3 | 4 | 5 | 6 | 7 | 8 | 9 | 10 | 11 | 12 | 13 | 14 | 15 | 16 | 17 | 18 | 19 | 20 | |

Our next step is to construct $T_{UPPERBOUND}$ by patching the cycles of $D_8$.

We now construct $D_8^{-1}M$. Using the 2-cycle *(11 12)*, we obtain

$$T_{UPPERBOUND} = (1\ 7\ 5\ 18\ 14\ 13\ 9\ 4\ 17\ 10\ 12\ 20\ 15\ 16\ 6\ 19\ 3\ 11\ 2\ 8): 213$$

Its value is 1 more than that if $D_8$. It also has the same value as $D_7$. Our next step is to see if $\sigma_{APOPT} = D_8$ contains any positive cycles

$$D_8^{-1} M^-$$

|   | 7 | 8 | 11 | 17 | 18 | 19 | 5 | 1 | 4 | 12 | 20 | 2 | 9 | 13 | 16 | 6 | 10 | 14 | 3 | 15 |   |
|---|---|---|----|----|----|----|---|---|---|----|----|---|---|----|----|---|----|----|---|----|---|
|   | 1 | 2 | 3 | 4 | 5 | 6 | 7 | 8 | 9 | 10 | 11 | 12 | 13 | 14 | 15 | 16 | 17 | 18 | 19 | 20 |   |
| 1 | 0 | 50 | 43 | -7 | 63 | 77 | 5 | ∞ | 88 | 41 | 56 | 79 | 30 | 20 | 56 | 29 | 37 | 8 | 63 | 39 | 1 |
| 2 | 79 | 0 | 4 | 34 | 42 | 5 | 47 | 60 | 52 | 56 | 25 | ∞ | 38 | 12 | 42 | -2 | 46 | 25 | 8 | -9 | 2 |
| 3 | 60 | 63 | 0 | 43 | 19 | 61 | 21 | 42 | 49 | 20 | 68 | 27 | 84 | 10 | 74 | 13 | 4 | 81 | ∞ | 63 | 3 |
| 4 | 85 | 79 | 29 | 0 | 58 | 84 | 83 | 20 | ∞ | 50 | 22 | 67 | 78 | 85 | 9 | 54 | 9 | 53 | 68 | 24 | 4 |
| 5 | 2 | 45 | 22 | -9 | 0 | 80 | ∞ | -2 | 61 | -1 | 16 | 13 | 58 | 72 | 10 | 30 | 35 | -15 | 41 | 16 | 5 |
| 6 | 76 | 72 | 8 | 47 | 48 | 0 | -3 | -1 | 53 | 21 | 49 | -3 | 50 | 26 | 60 | ∞ | 18 | -2 | 9 | 72 | 6 |
| 7 | ∞ | 70 | 77 | 87 | 91 | 54 | 0 | -3 | 64 | 55 | 85 | 95 | 24 | 20 | 6 | 31 | 40 | 29 | 78 | 5 | 7 |
| 8 | 51 | ∞ | 24 | 23 | 77 | 66 | 51 | 0 | 48 | 47 | 69 | 48 | 58 | 55 | 58 | 29 | 91 | 64 | 70 | 33 | 8 |
| 9 | 60 | 13 | 84 | 71 | 62 | 78 | 57 | 86 | 0 | 18 | 60 | 93 | ∞ | 86 | 37 | 73 | 57 | 15 | 43 | -2 | 9 |
| 10 | 82 | 82 | 24 | 63 | 89 | 75 | 70 | 36 | -6 | 0 | -5 | 55 | 62 | 76 | 30 | 44 | ∞ | 53 | -8 | 53 | 10 |
| 11 | 27 | 32 | ∞ | 35 | 66 | 12 | 72 | -8 | -1 | 71 | 0 | -5 | 13 | 11 | 47 | 18 | 4 | 66 | 62 | -1 | 11 |
| 12 | 78 | 88 | 76 | 67 | 64 | 24 | 5 | 42 | 63 | ∞ | 6 | 0 | 38 | 77 | 61 | 38 | 43 | 41 | 51 | 44 | 12 |
| 13 | 0 | 66 | 30 | 1 | 52 | -8 | 5 | 4 | 61 | 52 | -11 | 67 | 0 | ∞ | 77 | 65 | 72 | 39 | 78 | -6 | 13 |
| 14 | 12 | 33 | 4 | -15 | 76 | 40 | 37 | 65 | 71 | 35 | 7 | 54 | 9 | 0 | 18 | 77 | 7 | ∞ | 70 | 72 | 14 |
| 15 | 44 | 63 | 80 | 61 | 21 | 32 | 47 | 66 | 24 | 34 | 73 | 35 | 72 | 41 | 0 | 36 | 62 | 26 | 47 | ∞ | 15 |
| 16 | 83 | 12 | 57 | 64 | 56 | 15 | 76 | 8 | 55 | 7 | 52 | 0 | 53 | 54 | ∞ | 0 | 43 | 15 | 72 | 43 | 16 |
| 17 | 82 | 23 | 18 | ∞ | 70 | 30 | 7 | 45 | 58 | 74 | 5 | 47 | 27 | 16 | 74 | 82 | 0 | 27 | 24 | 81 | 17 |
| 18 | 58 | 42 | 73 | 59 | ∞ | 42 | 17 | 22 | 89 | 81 | 32 | 41 | 74 | 33 | 17 | 71 | 71 | 0 | 75 | -6 | 18 |
| 19 | 20 | 11 | 57 | 26 | 53 | ∞ | 52 | 47 | 38 | 57 | 62 | 83 | 58 | 81 | 51 | 42 | 82 | 55 | 0 | -6 | 19 |
| 20 | 8 | 50 | 18 | 85 | 79 | 78 | 43 | 37 | 45 | 74 | ∞ | 0 | 29 | 59 | 49 | 46 | 12 | 41 | 67 | 0 | 20 |
|   | 1 | 2 | 3 | 4 | 5 | 6 | 7 | 8 | 9 | 10 | 11 | 12 | 13 | 14 | 15 | 16 | 17 | 18 | 19 | 20 |   |

j = 4

 (5 4)(4 15) = (5 15): 0

 (5 4)(4 17) = (5 17): 0

 (14 4)(4 15) = (14 15): -6

 (14 4)(4 17) = (14 17): -6

j = 6

 (13 6)(6 3) = (13 3): -14

 (13 6)(6 7) = (13 7): -11

 (13 6)(6 8) = (13 8): -9

j = 7

 (6 7)(7 8) = (6 8): -6

 (13 7)(7 8) = (13 8): -14

j = 9

 (11 9)(9 20) = (11 20): -3

j = 10

 (5 10)(10 9) = (5 9): -7

 (5 10)(10 11) = (5 11): -6

j = 11

 (5 11)(11 8) = (5 8): -14

 (5 11)(11 9) = (5 9): -7

 (10 11)(11 8) = (10 8): -13

 (10 11)(11 12) = (10 12): -10

 (10 11)(11 17) = (10 17): -1

 (13 11)(11 8) = (13 8): -19

 (13 11)(11 9) = (13 9): -12

 (13 11)(11 17) = (13 17): -7

 (13  11)(11  20) = (13  20): -12

j = 12

 (10  12)(12  7) = (10  7): -5

 (11  12)(12  7) = (11  7): 0

j = 18

 (5  18)(18  20) = (5  20): -21

 (6  18)(18  20) = (6  20): -8

j = 19

 (10  19)(19  20) = (10  20): -14

j = 20

 (2  20)(20  1) = (2  1): -1

 (2  20)(20  12) = (2  12): -9

 (5  20)(20  1) = (5  1): -13

 (5  20)(20  3) = (5  3): -3

 (5  20)(20  12) = (5  12): -21

 (5  20)(20  17) = (5  17): -9

 (6  20)(20  1) = (6  1): 0

 (6  20)(20  12) = (6  12): -8

 (9  20)(20  12) = (9  12): -2

 (10  20)(20  1) = (10  1): -6

 (10  20)(20  12) = (10  12): -14

 (10  20)(20  17) = (10  17): -2

 (13  20)(20  1) = (13  1): -4

 (13  20)(20  12) = (13  12): -12

 (18  20)(20  12) = (18  12): -6

(19 20)(20 12) = (19 12): -6

$$D_8^{-1}M^-(20)$$

|  | 7 | 8 | 11 | 17 | 18 | 19 | 5 | 1 | 4 | 12 | 20 | 2 | 9 | 13 | 16 | 6 | 10 | 14 | 3 | 15 |  |
|---|---|---|---|---|---|---|---|---|---|---|---|---|---|---|---|---|---|---|---|---|---|
|  | 1 | 2 | 3 | 4 | 5 | 6 | 7 | 8 | 9 | 10 | 11 | 12 | 13 | 14 | 15 | 16 | 17 | 18 | 19 | 20 |  |
| 1 | 0 | 50 | 43 | -7 | 63 | 77 | 5 | ∞ | 88 | 41 | 56 | 79 | 30 | 20 | 56 | 29 | 37 | 8 | 63 | 39 | 1 |
| 2 | -1 | 0 | 4 | 34 | 42 | 5 | 47 | 60 | 52 | 56 | 25 | -9 | 38 | 12 | 42 | -2 | 46 | 25 | 8 | *-9* | 2 |
| 3 | 60 | 63 | 0 | 43 | 19 | 61 | 21 | 42 | 49 | 20 | 68 | 27 | 84 | 10 | 74 | 13 | 4 | 81 | ∞ | 63 | 3 |
| 4 | 85 | 79 | 29 | 0 | 58 | 84 | 83 | 20 | ∞ | 50 | 22 | 67 | 78 | 85 | 9 | 54 | 9 | 53 | 68 | 24 | 4 |
| 5 | -13 | 45 | -3 | -9 | 0 | 80 | ∞ | -14 | -7 | -1 | -6 | -21 | 58 | 72 | *0* | 30 | -9 | *-15* | 41 | *-21* | 5 |
| 6 | 0 | 72 | 8 | 47 | 48 | 0 | -3 | -6 | 53 | 21 | 49 | -8 | 50 | 26 | 60 | ∞ | 18 | *-2* | 9 | *-8* | 6 |
| 7 | ∞ | 70 | 77 | 87 | 91 | 54 | 0 | -3 | 64 | 55 | 85 | 95 | 24 | 20 | 6 | 31 | 40 | 29 | 78 | 5 | 7 |
| 8 | 51 | ∞ | 24 | 23 | 77 | 66 | 51 | 0 | 48 | 47 | 69 | 48 | 58 | 55 | 58 | 29 | 91 | 64 | 70 | 33 | 8 |
| 9 | 60 | 13 | 84 | 71 | 62 | 78 | 57 | 86 | 0 | 18 | 60 | -2 | ∞ | 86 | 37 | 73 | 57 | 15 | 43 | *-2* | 9 |
| 10 | -6 | 82 | 24 | 63 | 89 | 75 | -5 | -13 | -6 | 0 | -5 | -14 | 62 | 76 | 30 | 44 | -2 | 53 | -8 | *-14* | 10 |
| 11 | 27 | 32 | ∞ | 35 | 66 | 12 | 0 | -8 | -1 | 71 | 0 | -5 | 13 | 11 | 47 | 18 | 4 | 66 | 62 | *-3* | 11 |
| 12 | 78 | 88 | 76 | 67 | 64 | 24 | 5 | 42 | 63 | ∞ | 6 | 0 | 38 | 77 | 61 | 38 | 43 | 41 | 51 | 44 | 12 |
| 13 | -4 | 66 | 0 | 1 | 52 | -8 | *-11* | *-19* | *-12* | 52 | *-11* | *-12* | 0 | ∞ | 77 | 65 | -7 | 39 | 78 | *-12* | 13 |
| 14 | 12 | 33 | 4 | *-15* | 76 | 40 | 37 | 65 | 71 | 35 | 7 | 54 | 9 | 0 | -6 | 77 | -6 | ∞ | 70 | 72 | 14 |
| 15 | 44 | 63 | 80 | 61 | 21 | 32 | 47 | 66 | 24 | 34 | 73 | 35 | 72 | 41 | 0 | 36 | 62 | 26 | 47 | ∞ | 15 |
| 16 | 83 | 12 | 57 | 64 | 56 | 15 | 76 | 8 | 55 | 7 | 52 | 0 | 53 | 54 | ∞ | 0 | 43 | 15 | 72 | 43 | 16 |
| 17 | 82 | 23 | 18 | ∞ | 70 | 30 | 7 | 45 | 58 | 74 | 5 | 47 | 27 | 16 | 74 | 82 | 0 | 27 | 24 | 81 | 17 |
| 18 | 58 | 42 | 73 | 59 | ∞ | 42 | 17 | 22 | 89 | 81 | 32 | -6 | 74 | 33 | 17 | 71 | 71 | 0 | 75 | -6 | 18 |
| 19 | 20 | 11 | 57 | 26 | 53 | ∞ | 52 | 47 | 38 | 57 | 62 | -6 | 58 | 81 | 51 | 42 | 82 | 55 | 0 | -6 | 19 |
| 20 | 8 | 50 | 18 | 85 | 79 | 78 | 43 | 37 | 45 | 74 | ∞ | 0 | 29 | 59 | 49 | 46 | 12 | 41 | 67 | 0 | 20 |
|  | 1 | 2 | 3 | 4 | 5 | 6 | 7 | 8 | 9 | 10 | 11 | 12 | 13 | 14 | 15 | 16 | 17 | 18 | 19 | 20 |  |

$P_{20}$

|    | 1  | 2 | 3  | 4  | 5  | 6  | 7  | 8  | 9  | 10 | 11 | 12 | 13 | 14 | 15 | 16 | 17 | 18 | 19 | 20 |    |
|----|----|---|----|----|----|----|----|----|----|----|----|----|----|----|----|----|----|----|----|----|----|
| 1  |    |   |    |    |    |    |    |    |    |    |    |    |    |    |    |    |    |    |    |    | 1  |
| 2  | 20 |   |    |    |    |    |    |    |    |    |    | 20 |    |    |    |    |    |    | 2  |    | 2  |
| 3  |    |   |    |    |    |    |    |    |    |    |    |    |    |    |    |    |    |    |    |    | 3  |
| 4  |    |   |    |    |    |    |    |    |    |    |    |    |    |    |    |    |    |    |    |    | 4  |
| 5  | 20 |   | 20 | 5  |    |    |    | 11 | 11 | 5  | 10 | 20 |    |    | 4  |    | 20 | 5  |    | 18 | 5  |
| 6  | 20 |   |    |    |    | 6  | 7  |    |    |    |    | 20 |    |    |    |    | 6  |    |    | 18 | 6  |
| 7  |    |   |    |    |    |    |    |    |    |    |    |    |    |    |    |    |    |    |    |    | 7  |
| 8  |    |   |    |    |    |    |    |    |    |    |    |    |    |    |    |    |    |    |    |    | 8  |
| 9  |    |   |    |    |    |    |    |    |    |    |    | 20 |    |    |    |    |    |    | 9  |    | 9  |
| 10 | 20 |   |    |    |    |    | 12 | 11 | 10 |    | 10 | 20 |    |    |    |    | 20 |    | 10 | 19 | 10 |
| 11 |    |   |    |    |    |    | 12 |    | 11 |    | 11 |    |    |    |    |    |    |    |    | 9  | 11 |
| 12 |    |   |    |    |    |    |    |    |    |    |    |    |    |    |    |    |    |    |    |    | 12 |
| 13 | 20 |   | 6  |    |    | 13 | 6  | 11 | 11 |    | 13 | 20 |    |    |    |    | 11 |    |    | 11 | 13 |
| 14 |    |   | 14 |    |    |    |    |    |    |    |    |    |    |    | 4  |    | 4  |    |    |    | 14 |
| 15 |    |   |    |    |    |    |    |    |    |    |    |    |    |    |    |    |    |    |    |    | 15 |
| 16 |    |   |    |    |    |    |    |    |    |    |    |    |    |    |    |    |    |    |    |    | 16 |
| 17 |    |   |    |    |    |    |    |    |    |    |    |    |    |    |    |    |    |    |    |    | 17 |
| 18 |    |   |    |    |    |    |    |    |    |    |    | 20 |    |    |    |    |    |    |    | 18 | 18 |
| 19 |    |   |    |    |    |    |    |    |    |    |    | 20 |    |    |    |    |    |    |    | 19 | 19 |
| 20 |    |   |    |    |    |    |    |    |    |    |    |    |    |    |    |    |    |    |    |    | 20 |
|    | 1  | 2 | 3  | 4  | 5  | 6  | 7  | 8  | 9  | 10 | 11 | 12 | 13 | 14 | 15 | 16 | 17 | 18 | 19 | 20 |    |

*NONNEGPATHS(20)*

(1) (2 1), (10 1), (13 1); (3) ( 3); (4) (6 4); (7) (7 5), (10 5); (8) (5 8), (6 8); (9) (5 9),

(13 9); (10) (5 10); (12) (2 12), (5 12), (10 12), (13 12), (18 12), (19 12); (17) (5 17), (10 17), (13 17).

$j = 1$

   (2 1)(1 4) = (2 4): -8

   (5 1)(1 4) = (5 4): -20

   (5 1)(1 7) = (5 7): -8

   (5 1)(1 18) = (5 18): -5

   (6 1)(1 4) = (6 4): -7

   (10 1)(1 4) = (10 4): -13

   (10 1)(1 7) = (10 7): -1

   (13 1)(1 4) = (13 4): -11

$j = 4$

   (5 4)(4 15) = (5 15): -11

   (5 4)(4 17) = (5 17): -11

   (10 4)(4 15) = (10 15): -4

   (10 4)(4 17) = (10 17): -4

   (13 4)(4 15) = (13 15): -2

$j = 9$

   (13 9)(9 20) = (13 20): -14

$j = 12$

   (2 12)(12 7) = (2 7): -4

   (2 12)(12 11) = (2 11): -3

   (5 12)(12 7) = (5 7): -16

   (5 12)(12 11) = (5 11): -15

   (10 12)(12 7) = (10 7): -9

   (10 12)(12 11) = (10 11): -8

    (18 12)(12 11) = (18 11): 0
    (19 12)(12 7) = (19 7): -1
j = 20

    (13 20)(20 1) = (13 1): -6

$$D_8^{-1}M^-(40)$$

| | 7 | 8 | 11 | 17 | 18 | 19 | 5 | 1 | 4 | 12 | 20 | 2 | 9 | 13 | 16 | 6 | 10 | 14 | 3 | 15 | |
|---|---|---|---|---|---|---|---|---|---|---|---|---|---|---|---|---|---|---|---|---|---|
| | 1 | 2 | 3 | 4 | 5 | 6 | 7 | 8 | 9 | 10 | 11 | 12 | 13 | 14 | 15 | 16 | 17 | 18 | 19 | 20 | |
| 1 | 0 | 50 | 43 | -7 | 63 | 77 | 5 | ∞ | 88 | 41 | 56 | 79 | 30 | 20 | 56 | 29 | 37 | 8 | 63 | 39 | 1 |
| 2 | *-1* | 0 | 4 | *-8* | 42 | 5 | *-4* | 60 | 52 | 56 | *-3* | -9 | 38 | 12 | 42 | -2 | 46 | 25 | 8 | *-9* | 2 |
| 3 | 60 | 63 | 0 | 43 | 19 | 61 | 21 | 42 | 49 | 20 | 68 | 27 | 84 | 10 | 74 | 13 | 4 | 81 | ∞ | 63 | 3 |
| 4 | 85 | 79 | 29 | 0 | 58 | 84 | 83 | 20 | ∞ | 50 | 22 | 67 | 78 | 85 | 9 | 54 | 9 | 53 | 68 | 24 | 4 |
| 5 | *-13* | 45 | *-3* | *-20* | 0 | 80 | *-16* | *-14* | -7 | *-1* | *-15* | *-21* | 58 | 72 | *-11* | 30 | *-11* | *-15* | 41 | *-21* | 5 |
| 6 | *0* | 72 | 8 | -7 | 48 | 0 | -3 | -6 | 53 | 21 | 49 | *-8* | 50 | 26 | 60 | ∞ | 18 | -2 | 9 | *-8* | 6 |
| 7 | ∞ | 70 | 77 | 87 | 91 | 54 | 0 | -3 | 64 | 55 | 85 | 95 | 24 | 20 | 6 | 31 | 40 | 29 | 78 | 5 | 7 |
| 8 | 51 | ∞ | 24 | 23 | 77 | 66 | 51 | 0 | 48 | 47 | 69 | 48 | 58 | 55 | 58 | 29 | 91 | 64 | 70 | 33 | 8 |
| 9 | 60 | 13 | 84 | 71 | 62 | 78 | 57 | 86 | 0 | 18 | 60 | *-2* | ∞ | 86 | 37 | 73 | 57 | 15 | 43 | *-2* | 9 |
| 10 | -6 | 82 | 24 | *-13* | 89 | 75 | -5 | *-13* | -6 | 0 | *-8* | *-14* | 62 | 76 | *-4* | 44 | *-4* | 53 | *-8* | *-14* | 10 |
| 11 | 27 | 32 | ∞ | 35 | 66 | 12 | *0* | -8 | -1 | 71 | *0* | -5 | 13 | 11 | 47 | 18 | 4 | 66 | 62 | *-3* | 11 |
| 12 | 78 | 88 | 76 | 67 | 64 | 24 | 5 | 42 | 63 | ∞ | 6 | 0 | 38 | 77 | 61 | 38 | 43 | 41 | 51 | 44 | 12 |
| 13 | *-6* | 66 | *0* | *-11* | 52 | -8 | *-11* | *-19* | *-12* | 52 | *-11* | *-12* | 0 | ∞ | -2 | 65 | -7 | 39 | 78 | *-14* | 13 |
| 14 | 12 | 33 | 4 | *-15* | 76 | 40 | 37 | 65 | 71 | 35 | 7 | 54 | 9 | 0 | -6 | 77 | -6 | ∞ | 70 | 72 | 14 |
| 15 | 44 | 63 | 80 | 61 | 21 | 32 | 47 | 66 | 24 | 34 | 73 | 35 | 72 | 41 | 0 | 36 | 62 | 26 | 47 | ∞ | 15 |
| 16 | 83 | 12 | 57 | 64 | 56 | 15 | 76 | 8 | 55 | 7 | 52 | 0 | 53 | 54 | ∞ | 0 | 43 | 15 | 72 | 43 | 16 |
| 17 | 82 | 23 | 18 | ∞ | 70 | 30 | 7 | 45 | 58 | 74 | 5 | 47 | 27 | 16 | 74 | 82 | 0 | 27 | 24 | 81 | 17 |
| 18 | 58 | 42 | 73 | 59 | ∞ | 42 | *-1* | 22 | 89 | 81 | *0* | -6 | 74 | 33 | 17 | 71 | 71 | 0 | 75 | -6 | 18 |
| 19 | 20 | 11 | 57 | 26 | 53 | ∞ | *-1* | 47 | 38 | 57 | *0* | -6 | 58 | 81 | 51 | 42 | 82 | 55 | 0 | -6 | 19 |
| 20 | 8 | 50 | 18 | 85 | 79 | 78 | 43 | 37 | 45 | 74 | ∞ | 0 | 29 | 59 | 49 | 46 | 12 | 41 | 67 | 0 | 20 |

$P_{40}$

|    | 1  | 2 | 3 | 4 | 5 | 6 | 7  | 8 | 9 | 10 | 11 | 12 | 13 | 14 | 15 | 16 | 17 | 18 | 19 | 20 |    |
|----|----|---|---|---|---|---|----|---|---|----|----|----|----|----|----|----|----|----|----|----|----|
| 1  |    |   |   |   |   |   |    |   |   |    |    |    |    |    |    |    |    |    |    |    | 1  |
| 2  | 2  |   |   |   |   |   | 1  |   |   |    | 12 | 2  |    |    |    |    |    |    |    |    | 2  |
| 3  |    |   |   |   |   |   |    |   |   |    |    |    |    |    |    |    |    |    |    |    | 3  |
| 4  |    |   |   |   |   |   |    |   |   |    |    |    |    |    |    |    |    |    |    |    | 4  |
| 5  |    |   |   |   |   |   | 12 |   |   |    | 12 | 5  |    |    |    |    |    |    |    |    | 5  |
| 6  |    |   |   |   |   |   |    |   |   |    |    |    |    |    |    |    |    |    |    |    | 6  |
| 7  |    |   |   |   |   |   |    |   |   |    |    |    |    |    |    |    |    |    |    |    | 7  |
| 8  |    |   |   |   |   |   |    |   |   |    |    |    |    |    |    |    |    |    |    |    | 8  |
| 9  |    |   |   |   |   |   |    |   |   |    |    |    |    |    |    |    |    |    |    |    | 9  |
| 10 |    |   |   |   |   |   |    |   |   |    | 12 | 10 |    |    |    |    |    |    |    |    | 10 |
| 11 |    |   |   |   |   |   |    |   |   |    |    |    |    |    |    |    |    |    |    |    | 11 |
| 12 |    |   |   |   |   |   |    |   |   |    |    |    |    |    |    |    |    |    |    |    | 12 |
| 13 | 20 |   |   |   |   |   |    |   |   |    |    |    |    |    |    |    |    |    | 13 |    | 13 |
| 14 |    |   |   |   |   |   |    |   |   |    |    |    |    |    |    |    |    |    |    |    | 14 |
| 15 |    |   |   |   |   |   |    |   |   |    |    |    |    |    |    |    |    |    |    |    | 15 |
| 16 |    |   |   |   |   |   |    |   |   |    |    |    |    |    |    |    |    |    |    |    | 16 |
| 17 |    |   |   |   |   |   |    |   |   |    |    |    |    |    |    |    |    |    |    |    | 17 |
| 18 |    |   |   |   |   |   | 12 |   |   |    | 12 | 18 |    |    |    |    |    |    |    |    | 18 |
| 19 |    |   |   |   |   |   | 12 |   |   |    | 12 | 19 |    |    |    |    |    |    |    |    | 19 |
| 20 |    |   |   |   |   |   |    |   |   |    |    |    |    |    |    |    |    |    |    |    | 20 |
|    | 1  | 2 | 3 | 4 | 5 | 6 | 7  | 8 | 9 | 10 | 11 | 12 | 13 | 14 | 15 | 16 | 17 | 18 | 19 | 20 |    |

*NONNEGPATHS(40)*

*NONNEGPATHS(40)*

(1) (13 1); (4) (13 4); (7) (2 7); (11) (2 11), (5 11), (10 11); (12) (6 12).

j = 1

   (13 1)(1 4) = (13 4): -13

   (13 1)(1 7) = (13 7): -1

j = 4

   (13 4)(4 15) = (13 15): -4

j = 7

   (2 7)(7 8) = (2 8): -7

j = 11

   (2 11)(11 8) = (2 8): -11

   (5 11)(11 6) = (5 6): -3

   (5 11)(11 8) = (5 8): -23

   (5 11)(11 9) = (5 9): -16

   (10 11)(11 8) = (10 8): -16

   (10 11)(11 9) = (10 9): -9

j = 12

   (6 12)(12 11) = (6 11): -2

*NONNEGPATHS(60)*

(6) (5 6); (8) (5 8), (10 8); (9) (5 9), (10 9); (11) (6 11).

j = 8

   (5 9)(9 2) = (5 2): -3

j = 11

   (6 11)(11 8) = (6 8): -10

   (6 11)(11 9) = (6 9): -3

$$D_8^{-1}M^-(60)$$

|  | 7 | 8 | 11 | 17 | 18 | 19 | 5 | 1 | 4 | 12 | 20 | 2 | 9 | 13 | 16 | 6 | 10 | 14 | 3 | 15 |  |
|---|---|---|---|---|---|---|---|---|---|---|---|---|---|---|---|---|---|---|---|---|---|
|  | 1 | 2 | 3 | 4 | 5 | 6 | 7 | 8 | 9 | 10 | 11 | 12 | 13 | 14 | 15 | 16 | 17 | 18 | 19 | 20 |  |
| 1 | 0 | 50 | 43 | *-7* | 63 | 77 | 5 | ∞ | 88 | 41 | 56 | 79 | 30 | 20 | 56 | 29 | 37 | 8 | 63 | 39 | 1 |
| 2 | *-1* | 0 | 4 | *-8* | 42 | 5 | *-4* | 60 | 52 | 56 | *-3* | *-9* | 38 | 12 | 42 | -2 | 46 | 25 | 8 | *-9* | 2 |
| 3 | 60 | 63 | 0 | 43 | 19 | 61 | 21 | 42 | 49 | 20 | 68 | 27 | 84 | 10 | 74 | 13 | 4 | 81 | ∞ | 63 | 3 |
| 4 | 85 | 79 | 29 | 0 | 58 | 84 | 83 | 20 | ∞ | 50 | 22 | 67 | 78 | 85 | 9 | 54 | 9 | 53 | 68 | 24 | 4 |
| 5 | *-13* | <u>*-3*</u> | *-3* | *-20* | 0 | *-3* | *-16* | *-14* | *-16* | *-1* | *-15* | *-21* | 58 | 72 | *-11* | 30 | *-11* | *-15* | 41 | *-21* | 5 |
| 6 | *0* | 72 | 8 | *-7* | 48 | 0 | *-3* | <u>*-10*</u> | <u>*-3*</u> | 21 | *-2* | *-8* | 50 | 26 | 60 | ∞ | 18 | *-2* | 9 | *-8* | 6 |
| 7 | ∞ | 70 | 77 | 87 | 91 | 54 | 0 | -3 | 64 | 55 | 85 | 95 | 24 | 20 | 6 | 31 | 40 | 29 | 78 | 5 | 7 |
| 8 | 51 | ∞ | 24 | 23 | 77 | 66 | 51 | 0 | 48 | 47 | 69 | 48 | 58 | 55 | 58 | 29 | 91 | 64 | 70 | 33 | 8 |
| 9 | 60 | 13 | 84 | 71 | 62 | 78 | 57 | 86 | 0 | 18 | 60 | *-2* | ∞ | 86 | 37 | 73 | 57 | 15 | 43 | *-2* | 9 |
| 10 | *-6* | 82 | 24 | *-13* | 89 | 75 | -5 | *-16* | -9 | 0 | *-8* | *-14* | 62 | 76 | *-4* | 44 | *-4* | 53 | *-8* | *-14* | 10 |
| 11 | 27 | 32 | ∞ | 35 | 66 | 12 | <u>0</u> | -8 | *-1* | 71 | *0* | -5 | 13 | 11 | 47 | 18 | 4 | 66 | 62 | *-3* | 11 |
| 12 | 78 | 88 | 76 | 67 | 64 | 24 | 5 | 42 | 63 | ∞ | 6 | 0 | 38 | 77 | 61 | 38 | 43 | 41 | 51 | 44 | 12 |
| 13 | -6 | 66 | *0* | *-13* | 52 | *-8* | *-11* | *-19* | *-12* | 52 | *-11* | *-12* | 0 | ∞ | *-4* | 65 | *-7* | 39 | 78 | *-14* | 13 |
| 14 | 12 | 33 | 4 | *-15* | 76 | 40 | 37 | 65 | 71 | 35 | 7 | 54 | 9 | 0 | -6 | 77 | -6 | ∞ | 70 | 72 | 14 |
| 15 | 44 | 63 | 80 | 61 | 21 | 32 | 47 | 66 | 24 | 34 | 73 | 35 | 72 | 41 | 0 | 36 | 62 | 26 | 47 | ∞ | 15 |
| 16 | 83 | 12 | 57 | 64 | 56 | 15 | 76 | 8 | 55 | 7 | 52 | 0 | 53 | 54 | ∞ | 0 | 43 | 15 | 72 | 43 | 16 |
| 17 | 82 | 23 | 18 | ∞ | 70 | 30 | 7 | 45 | 58 | 74 | 5 | 47 | 27 | 16 | 74 | 82 | 0 | 27 | 24 | 81 | 17 |
| 18 | 58 | 42 | 73 | 59 | ∞ | 42 | *-1* | 22 | 89 | 81 | <u>0</u> | -6 | 74 | 33 | 17 | 71 | 71 | 0 | 75 | *-6* | 18 |
| 19 | 20 | 11 | 57 | 26 | 53 | ∞ | *-1* | 47 | 38 | 57 | <u>0</u> | -6 | 58 | 81 | 51 | 42 | 82 | 55 | 0 | *-6* | 19 |
| 20 | 8 | 50 | 18 | 85 | 79 | 78 | 43 | 37 | 45 | 74 | ∞ | 0 | 29 | 59 | 49 | 46 | 12 | 41 | 67 | 0 | 20 |
|  | 1 | 2 | 3 | 4 | 5 | 6 | 7 | 8 | 9 | 10 | 11 | 12 | 13 | 14 | 15 | 16 | 17 | 18 | 19 | 20 |  |

$P_{60}$

|    | 1 | 2 | 3 | 4 | 5 | 6 | 7 | 8 | 9 | 10 | 11 | 12 | 13 | 14 | 15 | 16 | 17 | 18 | 19 | 20 |    |
|----|---|---|---|---|---|---|---|---|---|----|----|----|----|----|----|----|----|----|----|----|----|
| 1  |   |   |   |   |   |   |   |   |   |    |    |    |    |    |    |    |    |    |    |    | 1  |
| 2  |   |   |   |   |   |   |   |   |   |    |    |    |    |    |    |    |    |    |    |    | 2  |
| 3  |   |   |   |   |   |   |   |   |   |    |    |    |    |    |    |    |    |    |    |    | 3  |
| 4  |   |   |   |   |   |   |   |   |   |    |    |    |    |    |    |    |    |    |    |    | 4  |
| 5  |   | 9 |   |   |   |   |   |   | 5 |    |    |    |    |    |    |    |    |    |    |    | 5  |
| 6  |   |   |   |   |   |   |   | 11| 11|    | 6  |    |    |    |    |    |    |    |    |    | 6  |
| 7  |   |   |   |   |   |   |   |   |   |    |    |    |    |    |    |    |    |    |    |    | 7  |
| 8  |   |   |   |   |   |   |   |   |   |    |    |    |    |    |    |    |    |    |    |    | 8  |
| 9  |   |   |   |   |   |   |   |   |   |    |    |    |    |    |    |    |    |    |    |    | 9  |
| 10 |   |   |   |   |   |   |   |   |   |    |    |    |    |    |    |    |    |    |    |    | 10 |
| 11 |   |   |   |   |   |   |   |   |   |    |    |    |    |    |    |    |    |    |    |    | 11 |
| 12 |   |   |   |   |   |   |   |   |   |    |    |    |    |    |    |    |    |    |    |    | 12 |
| 13 |   |   |   |   |   |   |   |   |   |    |    |    |    |    |    |    |    |    |    |    | 13 |
| 14 |   |   |   |   |   |   |   |   |   |    |    |    |    |    |    |    |    |    |    |    | 14 |
| 15 |   |   |   |   |   |   |   |   |   |    |    |    |    |    |    |    |    |    |    |    | 15 |
| 16 |   |   |   |   |   |   |   |   |   |    |    |    |    |    |    |    |    |    |    |    | 16 |
| 17 |   |   |   |   |   |   |   |   |   |    |    |    |    |    |    |    |    |    |    |    | 17 |
| 18 |   |   |   |   |   |   |   |   |   |    |    |    |    |    |    |    |    |    |    |    | 18 |
| 19 |   |   |   |   |   |   |   |   |   |    |    |    |    |    |    |    |    |    |    |    | 19 |
| 20 |   |   |   |   |   |   |   |   |   |    |    |    |    |    |    |    |    |    |    |    | 20 |
|    | 1 | 2 | 3 | 4 | 5 | 6 | 7 | 8 | 9 | 10 | 11 | 12 | 13 | 14 | 15 | 16 | 17 | 18 | 19 | 20 |    |
104

(2) (5 2); (11) (6 8), (6 9).
j = 2

(5 2)(2 16) = (5 16): -5

Thus, the only path that can be extended is at (5 16). But, checking the possible extensions, none is non-negative. It follows that there exists no cycle, $s$, in $D_8^{-1}M^-$ such that $D_8 s$ is a tour. It follows that the $n$-cycle $D_7 = \sigma_{FWTSPOPT}$. Furthermore, from the corollary to theorem 3.6, $D_7 = \sigma_{TSPOPT}$.

As addenda, we include a full list of c-trees.

*Note*. In constructing the set of c-trees, we must include as a root of a tree every determining vertex of each cycle obtained while constructing $\sigma_{FWTSPOPT}$.